\documentclass[12pt]{article}

\usepackage[margin=1.17in]{geometry}
\usepackage{amssymb,amsmath}
\usepackage{amsthm,enumitem,mathtools,microtype}
\usepackage{hyperref}
\usepackage{mathrsfs}
\usepackage{textcomp}
\hypersetup{
  colorlinks,
  citecolor=black,
  filecolor=black,
  linkcolor=black,
  urlcolor=black,
}
\usepackage[abbrev,nobysame]{amsrefs}
\usepackage{todonotes}
\usepackage{verbatim}
\usepackage{sectsty}
\usepackage[T1]{fontenc}
\usepackage{lmodern}
\usepackage{needspace}

\numberwithin{equation}{section}
\newtheorem{theorem}{Theorem}[section]
\newtheorem{remark}[theorem]{Remark}
\newtheorem{proposition}[theorem]{Proposition}
\newtheorem{lemma}[theorem]{Lemma}
\newtheorem{definition}[theorem]{Definition}
\setlist[enumerate,1]{label=\textup{(\roman*)},leftmargin=2.4em}
\allowdisplaybreaks[1]
\title{\textbf{Classification and asymptotics for the boundary Yamabe equation with an isolated singularity}\bigskip}
\author{\medskip Tianling Jin\footnote{T. Jin is partially supported by NSFC grant 12122120, and Hong Kong RGC grants GRF 16304125, GRF 16303624  and GRF 16303822.},\ \ Yuxuan Liao,\ \ Jingang Xiong\footnote{J. Xiong is partially supported by NSFC grant 12325104.}}
\date{\today}

\begin{document}
\maketitle

\begin{abstract}

We classify positive singular solutions of the boundary Yamabe equation on the punctured half-space in every dimension $n\ge3$. Previous symmetry results reduce this problem to a nonlinear PDE in two variables, for which classical phase-plane methods are unavailable. We develop a Hamiltonian method that combines local bifurcation with compactness to prove that every nonstationary cylindrical profile of a singular solution is periodic and uniquely determined, up to translation, by its Hamiltonian energy. At each stationary or periodic profile, the space of solutions of the linearized equation is two-dimensional. Finally, every positive solution in a half-ball with a nonremovable isolated boundary singularity is asymptotic to a unique global singular solution, with an algebraically decaying relative error.

\medskip


\end{abstract}

\section{Introduction}

The classification of solutions is a central problem in nonlinear elliptic equations. Common approaches are to prove uniqueness, often for ground states or by exploiting an underlying convexity, or to establish radial symmetry and then apply ODE methods. These approaches alone do not resolve the classification of singular solutions of the boundary Yamabe equation, which we study in this paper.

Write $X=(x,y)\in\mathbb R^{n-1}\times\mathbb R$ with $n\ge3$, and set
$\mathbb R^n_+=\{X:y>0\}$. We consider positive solutions of the boundary
Yamabe equation
\begin{equation}\label{eq:halfspace}
\begin{cases}
\begin{aligned}
\Delta U&=0
&& \text{in }\mathbb R^n_+,\\
-\partial_yU&=U^{\frac{n}{n-2}}
&& \text{on }\mathbb R^{n-1}\setminus\{0\}.
\end{aligned}
\end{cases}
\end{equation}
Equation~\eqref{eq:halfspace} arises in blow-up analysis and the study
of isolated singularities for the boundary Yamabe problem; see
Escobar \cites{EscobarSF,EscobarYamabe} and Han--Li \cite{HanLi1999}.
If $U$ is locally bounded near $0$, the singularity is removable,
and the classification of Li--Zhu \cite{LiZhu1995} gives
\begin{equation}\label{eq:halfspace-bubble}
U(x,y)=U_{\mu,x_0}(x,y)
=\left(
\frac{(n-2)\mu}{|x-x_0|^2+(y+\mu)^2}
\right)^{\frac{n-2}{2}},
\quad \mu>0,\quad x_0\in\mathbb R^{n-1}.
\end{equation}
Otherwise, $0$ is a nonremovable singularity, and we call $U$ a
singular solution. The growth estimates and symmetry result of
Caffarelli--Jin--Sire--Xiong \cite{CJSX} give
\begin{equation}\label{eq:radial}
U(x,y)=U(|x|,y)\ \ \mbox{and}\ \  C^{-1}|X|^{-\frac{n-2}{2}} \le U(X)\le C|X|^{-\frac{n-2}{2}} \ \ \text{on }\overline{\mathbb{R}^n_+}\setminus\{0\}
\end{equation}
for some $C>1$. For the Yamabe equation
$-\Delta U=U^{\frac{n+2}{n-2}}$ in $\mathbb R^n\setminus\{0\}$,
radial symmetry and the Emden--Fowler transform reduce the problem to
an autonomous second-order ODE. Frank--K\"onig \cite{FrankKonig2019} also used ODE methods to classify positive singular solutions of the critical biharmonic equation. In the present problem \eqref{eq:halfspace}, axial symmetry in \eqref{eq:radial} 
leaves one angular variable, so the usual planar phase portrait does not
apply. As in the Yamabe problem, however, the cylindrical
equation admits a conserved Hamiltonian. In this paper, we develop a two-dimensional
linear theory and use this Hamiltonian to classify all positive
bounded cylindrical solutions by their energy.

We now describe this classification in cylindrical variables. For an
axially symmetric solution $U$, introduce the Emden--Fowler variables
\[
\begin{gathered}
X=\rho(\sin\theta\,\omega',\cos\theta),
\qquad t=-\log\rho,\qquad
U(X)=
\left(\frac{n-2}{2}\right)^{\frac{n-2}{2}}
\rho^{-\frac{n-2}{2}}v(t,\theta),
\end{gathered}
\]
with $\rho>0$, $0\le\theta\le\pi/2$, and $\omega'\in\mathbb S^{n-2}$.
Writing $m(\theta)=\sin^{n-2}\theta$ we obtain
\begin{equation}\label{eq:cylinder}
\begin{cases}
\begin{aligned}
v_{tt}+m^{-1}(mv_\theta)_\theta-\frac{(n-2)^2}{4}v&=0
&& \text{in }\mathbb R\times(0,\pi/2),\\
v_\theta(t,\pi/2)&=\frac{n-2}{2}v(t,\pi/2)^{\frac{n}{n-2}}
&& \text{for }t\in\mathbb R,\\
v_\theta(t,0)&=0
&& \text{for }t\in\mathbb R.
\end{aligned}
\end{cases}
\end{equation}
A nonnegative solution of \eqref{eq:cylinder} is either identically
zero or positive everywhere. By \eqref{eq:radial}, the cylindrical
profiles of global singular solutions are bounded. We therefore work
with \emph{positive bounded solutions} of \eqref{eq:cylinder}.
We identify each profile
$v(t,\theta)$ with its axially symmetric extension to
$\mathbb R\times\mathbb S^{n-1}_+$. The rotational extensions are understood to be smooth at $\theta=0$. We call a solution $v$ periodic if there exists $T>0$ such that
\[
v(t+T,\theta)=v(t,\theta)
\qquad
\text{for all }(t,\theta)\in\mathbb R\times[0,\pi/2].
\]

Multiplying the first equation in \eqref{eq:cylinder} by $mv_t$,
integrating in $\theta$, and using the boundary condition, one obtains the conserved
Hamiltonian
\begin{equation}\label{eq:hamiltonian}
\begin{aligned}
\mathcal H(v)
&=\frac12\int_0^{\pi/2}m\left(v_t^2-v_\theta^2-\frac{(n-2)^2}{4}v^2\right)\,\mathrm d\theta
+\frac{(n-2)^2}{4(n-1)}
v(t,\pi/2)^{\frac{2(n-1)}{n-2}}
\end{aligned}
\end{equation}
satisfying $$\frac{\mathrm d}{\mathrm dt}\mathcal H(v)=0.$$

There are two special types of positive solutions to \eqref{eq:cylinder}. One is the family of axial bubbles $U_{\mu,0}$ in
\eqref{eq:halfspace-bubble}, which become
\begin{equation}\label{eq:bubblecyl}
v_{\mathrm{bub}}^s(t,\theta)
:=\bigl(\cosh(t-s)+\cos\theta\bigr)^{-\frac{n-2}{2}},
\qquad s=-\log\mu\in\mathbb R
\end{equation}
in the Emden--Fowler variables. These profiles satisfy 
$$\mathcal H(v_{\mathrm{bub}}^s)=0.$$ 
The other is the unique positive stationary solution
$v_c=v_c(\theta)$, whose boundary trace is
\begin{equation}\label{eq:stationarytrace}
v_c(\pi/2)
=\left(\frac{4}{n-2}\right)^{\frac{n-2}{2}}
 \left(\frac{\Gamma(\frac n4)}{\Gamma(\frac{n-2}{4})}\right)^{n-2}.
\end{equation}
Its Hamiltonian is
\begin{equation}\label{eq:stationaryHam}
H_c:=\mathcal H(v_c)
=-\frac{n-2}{4(n-1)}v_c(\pi/2)^{\frac{2(n-1)}{n-2}}<0.
\end{equation}
Let $\kappa_c>0$ be the principal linearized frequency at $v_c$,
defined by \eqref{eq:prinfre}.
See Section~\ref{sec:stationary} for details.

It is natural to ask whether the range of the Hamiltonian $\mathcal H$ on the set
of positive bounded solutions of \eqref{eq:cylinder} is exactly
$[H_c,0]$, and whether $\mathcal H(v)$ determines $v$ uniquely
up to translation in $t$.
The following theorem answers both questions affirmatively.

\begin{theorem}\label{thm:fowler-family}
Every positive bounded solution $v$ of \eqref{eq:cylinder} satisfies
$\mathcal H(v)\in[H_c,0]$.
If $\mathcal H(v)=H_c$ then $v=v_c$.
If $\mathcal H(v)=0$, then $v=v_{\mathrm{bub}}^s$ for some
$s\in\mathbb R$.
For every $H\in(H_c,0)$ there is a nonstationary periodic
positive bounded solution $v_H$ with $\mathcal H(v_H)=H$ such that
every positive bounded solution $v$ with $\mathcal H(v)=H$ has the
form $v=v_H(\,\cdot+s,\cdot)$ for some $s\in\mathbb R$.
Let $w_H(t)=v_H(t,\pi/2)$ and choose the translate for which $w_H$
attains a maximum at $t=0$. Let $T(H)$ be its least period. Then
\[
v_H(-t,\theta)=v_H(t,\theta)
\qquad\mbox{and}\qquad
(v_H)_t(t,\theta)<0
\]
when $0<t<T(H)/2$ and $0\le\theta\le\pi/2$.
For every integer $j\ge0$, the map
\[
H\longmapsto
\bigl(T(H),v_H(T(H)\,\cdot,\cdot)\bigr)
\]
is locally real analytic into
$\mathbb R\times
H^j((\mathbb R/\mathbb Z)\times\mathbb S^{n-1}_+)$.
At the two ends of the energy interval,
\[
T(H)\to\frac{2\pi}{\kappa_c},
\qquad
v_H(T(H)\,\cdot,\cdot)\to v_c
\quad\text{smoothly as }H\downarrow H_c,
\]
while
\[
T(H)\to\infty,\quad
v_H\to(\cosh t+\cos\theta)^{-\frac{n-2}{2}}
\quad\text{smoothly on compact cylinders as }H\uparrow0.
\]
\end{theorem}

We call $v_c$, the translates of the periodic profiles $v_H$
with $H\in(H_c,0)$, and their corresponding half-space solutions
Fowler solutions (also called Delaunay-type solutions).

\begin{remark}
The nonstationary Fowler family can also be parametrized, modulo
translation in $t$, by its boundary minimum (necksize) or maximum.
Indeed, Proposition~\ref{prop:fowler-extremal-parameters} in the appendix
shows that
\[
H\longmapsto\min_{t\in\mathbb R}v_H(t,\pi/2),
\qquad
H\longmapsto\max_{t\in\mathbb R}v_H(t,\pi/2)
\]
are real-analytic diffeomorphisms from $(H_c,0)$ onto $(0,v_c(\pi/2))$ and
$(v_c(\pi/2),1)$, respectively, with the first strictly decreasing and the
second strictly increasing.
\end{remark}

The main difficulty in proving Theorem \ref{thm:fowler-family} is that, even after axial symmetry and the Emden--Fowler transformation, \eqref{eq:cylinder} remains a PDE in the two variables $(t,\theta)$, rather than an autonomous second-order ODE. Thus the classical phase-plane analysis for ODEs is unavailable.  This difficulty also arises in the fractional Yamabe problem. The boundary Yamabe equation \eqref{eq:halfspace} considered here corresponds, through harmonic extension, to fractional order $1/2$ on $\mathbb R^{n-1}$. There have been several important contributions toward the classification of singular solutions of the fractional Yamabe equation. DelaTorre--del Pino--Gonz\'alez--Wei \cite{DelaTorreEtAl} constructed periodic Fowler solutions. Ao--DelaTorre--Gonz\'alez--Wei \cite{AoEtAlGluing2020} constructed solutions with prescribed isolated singularities by a gluing method. Cruz-Bl\'azquez--DelaTorre--Ruiz \cite{CruzDelaTorreRuiz2023} proved that every global singular solution has infinite Morse index and that its cylindrical profile, unless stationary, takes values on both sides of the stationary value. More recently, Andrade--DelaTorre--do \'O--Ratzkin--Wei \cite{AndradeEtAl2025} classified global singular solutions for fractional orders sufficiently close to $1$. Uniqueness and nondegeneracy for the fractional Laplacian have also been studied in the ground-state setting. Frank--Lenzmann \cite{FrankLenzmann2013} proved uniqueness up to translation and nondegeneracy of ground states in one dimension, and Frank--Lenzmann--Silvestre \cite{FrankLenzmannSilvestre2016} 
extended these results to arbitrary dimensions for the stationary
nonlinear fractional Schr\"odinger equation with a
Sobolev-subcritical nonlinearity. Their proofs combine spectral analysis of the linearized operator with continuation in the fractional order.

Our Hamiltonian approach to Theorem~\ref{thm:fowler-family} proceeds
as follows.  We first establish a two-dimensionality theorem for an exponentially weighted class of entire axially symmetric solutions of linear equations with a variable Robin boundary coefficient. Spectral analysis shows
that each solution in this class is uniquely determined by two scalar projections of its Cauchy data. Conversely, every pair of principal Cauchy data is realized by a unique solution in this class. Furthermore, the conserved Wronskian is nondegenerate on this two-dimensional linear solution space.


We next prove local rigidity of the Hamiltonian level sets: near
a nonstationary profile, solutions with the same energy differ only
by translation, while the stationary profile $v_c$ is isolated in its energy level.
To prove this, we apply the linear estimates to differences of
nonlinear solutions, obtaining local bi-Lipschitz
embeddings of the solution set into $\mathbb R^2$ through two fixed
scalar projections of the Cauchy data. These embeddings allow us
to deduce local rigidity at nonstationary solutions from the first
variation of the Hamiltonian and the nondegeneracy of the Wronskian.
The isolation of $v_c$ in its energy level follows from the positive
definiteness of the second variation on the two-dimensional space
of admissible linearized solutions at $v_c$.

Finally, we establish two local descriptions of the nonlinear solution space. Near any nonstationary periodic profile, all bounded entire axial solutions are periodic and form a real-analytic two-dimensional manifold, on which the least period is real analytic. Near $v_c$, bifurcation produces a family of periodic solutions whose translates exhaust all nearby nonstationary solutions. The positive quadratic term in the Hamiltonian expansion gives exactly one periodic translation orbit at each energy sufficiently close to and above $H_c$. Compactness and local energy rigidity then force every solution with energy in $(H_c,0)$ to be periodic. Moreover, periodic continuation and compactness show that the number of translation orbits at each such energy is finite and locally constant. The local classification near $v_c$ therefore gives exactly one translation orbit at every interior energy.

Related results on the dimension of radial solution spaces were obtained by Ao--Chan--DelaTorre--Fontelos--Gonz\'alez--Wei \cite{AoEtAl2019}, who proved that certain such spaces for fractional Hardy operators are at most two-dimensional, under growth restrictions on the solutions and asymptotic assumptions on the potential. They developed a weighted linear theory for fractional operators arising in singular Yamabe problems, including bounds on radial kernels, injectivity results, and uniformly bounded right inverses for their gluing construction. See also \cite{AoEtAl2020}, where the same authors develop Hamiltonian identities and Wronskian-type quantities for radial nonlocal equations, with applications to the nondegeneracy of fast-decaying singular solutions of the fractional Lane--Emden
equation. In the present work, the coefficient bounds required by the linear theory follow directly from the nonlinear equation. Consequently, the two-dimensionality theorem applies without prior knowledge of periodicity or asymptotic behavior at either end of the cylinder, and can therefore be used to establish the nonlinear classification.

The linear theory here also classifies solutions of the linearized equation. At a nonstationary Fowler profile, translation and variation of the energy generate two distinguished solutions, analogous to the translation and necksize Jacobi fields studied by Mazzeo--Pacard \cite{MazzeoPacard1999} for the Yamabe equation. The next theorem shows that these two solutions exhaust this growth class, without imposing axial symmetry on the linearized solution.

\begin{theorem}\label{thm:linearized-classification}
Let $H\in[H_c,0)$, set $v_{H_c}=v_c$ and $w_H(t)=v_H(t,\pi/2)$,
and let
$h\in C^\infty(
\mathbb R\times\overline{\mathbb S^{n-1}_+};\mathbb R)$
satisfy
\begin{equation}\label{eq:full-linearized-fowler}
\begin{cases}
h_{tt}+\Delta_{\mathbb S^{n-1}_+}h-\frac{(n-2)^2}{4}h=0
&\text{in }\mathbb R\times\mathbb S^{n-1}_+,\\
\partial_\nu h=\frac n2w_H(t)^{\frac2{n-2}}h
&\text{on }\mathbb R\times\partial\mathbb S^{n-1}_+,
\end{cases}
\end{equation}
where $\nu$ is the outward unit normal to
$\partial\mathbb S^{n-1}_+$.
Assume that, for some $C>0$ and
$0<\varepsilon<\min\{1,n/4\}$,
\begin{equation}\label{eq:linearized-growth}
\sup_{\omega\in\overline{\mathbb S^{n-1}_+}}
|h(t,\omega)|
\le Ce^{\varepsilon|t|}
\quad\mbox{for all }t\in\mathbb R.
\end{equation}
Then $h$ is axially symmetric and, for unique $A,B\in\mathbb R$,
\begin{equation}\label{eq:linearized-classification}
h(t,\theta)=
\begin{cases}
A(v_H)_t(t,\theta)+B\,\partial_Hv_H(t,\theta),
&H_c<H<0,\\
\phi_c(\theta)
\bigl(A\cos(\kappa_ct)+B\sin(\kappa_ct)\bigr),
&H=H_c.
\end{cases}
\end{equation}
Here $\phi_c>0$ is the principal eigenfunction of
\eqref{eq:eigenvalue}, normalized by $\phi_c(\pi/2)=1$,
and $\partial_Hv_H$ is taken at fixed $t$.
Conversely, every function in
\eqref{eq:linearized-classification} solves
\eqref{eq:full-linearized-fowler} and has at most linear growth in $t$.
\end{theorem}

The Fowler solutions also describe the local behavior near nonremovable isolated boundary singularities, as in the Yamabe problem studied by Caffarelli--Gidas--Spruck \cite{CGS} and Korevaar--Mazzeo--Pacard--Schoen \cite{KMPS}, and also in the $\sigma_k$-Yamabe equation studied by Han--Li--Teixeira \cite{HLT2010}. For higher order conformally invariant equations, Jin--Xiong \cite{JinXiong2021} established sharp growth estimates and asymptotic radial symmetry near isolated singularities. The local estimates of Caffarelli--Jin--Sire--Xiong \cite{CJSX} give the critical growth rate near a nonremovable isolated boundary singularity and hence uniform bounds for the corresponding profile on a forward half-cylinder. The global classification identifies the possible limits of its translates. Using the classification of linearized solutions, we prove convergence to a single global singular profile, with a uniquely determined phase modulo the period when the limit is nonstationary. The convergence is exponential in cylindrical variables and gives an algebraically decaying relative error in the original variables. 

For $r>0$, write
\[
B_r^+=\{X=(x,y)\in\mathbb R^n:|X|<r,\ y>0\},
\qquad
B_r'=\{x\in\mathbb R^{n-1}:|x|<r\}.
\]

\begin{theorem}\label{thm:local-asymptotics}
Let $r_0>0$, and let  $U\in C^2(B_{r_0}^+)\cap C^1(B_{r_0}^+\cup B_{r_0}'\setminus\{0\})$ be a positive solution of
\[
\begin{cases}
\begin{aligned}
\Delta U&=0
&& \text{in }B_{r_0}^+,\\
-\partial_yU&=U^{\frac{n}{n-2}}
&& \text{on }B_{r_0}'\setminus\{0\}.
\end{aligned}
\end{cases}
\]
Assume that $0$ is a nonremovable singularity.
Then there is a unique global singular solution $U_\infty$ of
\eqref{eq:halfspace} and a number $\beta>0$ such that
\[
U(X)=U_\infty(X)\bigl(1+O(|X|^\beta)\bigr)
\qquad\mbox{as }X\to0.
\]
\end{theorem}

In cylindrical variables, the limiting profile is either $v_c$ or
$v_H(t+s,\theta)$ for some $H\in(H_c,0)$. In the latter case, the
energy $H$ is unique and the phase $s$ is unique modulo $T(H)$.

\begin{remark}
The Hamiltonian method developed here has been applied in our follow-up paper \cite{FractionalYamabeFowler} to extend the results of the present paper to the fractional Yamabe equation via the Caffarelli--Silvestre extension \cite{CS}. The computations in the general fractional case are substantially more involved. The boundary Yamabe equation studied here, which corresponds to fractional order $1/2$, is of independent geometric interest as it arises from the classical problem of prescribing boundary mean curvature. Moreover, this setting allows us to present the main ideas of the Hamiltonian method clearly and without computational complications.
\end{remark}

The paper is organized as follows. In
Section~\ref{sec:linear-theory}, we establish the sharp a priori
estimates and the two-dimensional weighted linear theory. In
Section~\ref{sec:hamiltonian-global}, we construct local Hamiltonian
coordinates, classify the global Fowler family, and complete the
classification of solutions of the linearized equation. In
Section~\ref{sec:local-singularities}, we prove convergence of local
singular solutions to their unique global models and establish the
quantitative rate. In the Appendix, we establish the parametrizations
of the Fowler family by its boundary minimum (necksize) and maximum.

\section{A priori estimates and two-dimensional linear theory}\label{sec:linear-theory}

\subsection{Stationary solutions and the principal eigenvalue}\label{sec:stationary}

Recall that $m(\theta)=\sin^{n-2}\theta$. We study positive bounded solutions of \eqref{eq:cylinder},
\[
\begin{cases}
\begin{aligned}
v_{tt}+m^{-1}(mv_\theta)_\theta-\frac{(n-2)^2}{4}v&=0
&& \text{in }\mathbb R\times(0,\pi/2),\\
v_\theta(t,\pi/2)&=\frac{n-2}{2}v(t,\pi/2)^{\frac{n}{n-2}}
&& \text{for }t\in\mathbb R,\\
v_\theta(t,0)&=0
&& \text{for }t\in\mathbb R.
\end{aligned}
\end{cases}
\]
Let $f_0$ be the unique regular solution of
\begin{equation}\label{eq:f0}
(mf_0')'=\frac{(n-2)^2}{4}mf_0
\quad\text{in }(0,\pi/2),\qquad
f_0'(0)=0,\qquad f_0(\pi/2)=1.
\end{equation}
The maximum principle gives $f_0>0$. Evaluating the symbol formula
\eqref{eq:symbol} at zero frequency will give
\[
b_0:=f_0'(\pi/2)
=\frac{2\Gamma(\frac n4)^2}{\Gamma(\frac{n-2}{4})^2}>0.
\]
Every positive stationary solution of \eqref{eq:cylinder} is a multiple of $f_0$.
For $w_c=v_c(\pi/2)$, the boundary condition
$b_0w_c=\frac{n-2}{2}w_c^{\frac{n}{n-2}}$
determines $w_c>0$ uniquely, yielding \eqref{eq:stationarytrace} and
\[
v_c=w_cf_0. 
\]
Testing the stationary equation with $v_c$
gives
\[
\int_0^{\pi/2}
m\left((v_c')^2+\frac{(n-2)^2}{4}v_c^2\right)\,\mathrm d\theta
=w_cv_c'(\pi/2)
=\frac{n-2}{2}w_c^{\frac{2(n-1)}{n-2}},
\]
and hence, by \eqref{eq:hamiltonian}, one obtains \eqref{eq:stationaryHam}. 

The linearization of \eqref{eq:cylinder} at $v_c$ is
\[
\begin{cases}
\begin{aligned}
h_{tt}+m^{-1}(mh_\theta)_\theta-\frac{(n-2)^2}{4}h&=0,
&& \text{in }\mathbb R\times(0,\pi/2),\\
h_\theta(t,\pi/2)&=b_ch(t,\pi/2),
&& \text{for }t\in\mathbb R,\\
h_\theta(t,0)&=0,
&& \text{for }t\in\mathbb R,
\end{aligned}
\end{cases}
\]
where $$b_c:=\frac n2w_c^{\frac{2}{n-2}}=\frac{n}{n-2}b_0>b_0.$$ 
The corresponding angular eigenvalue problem is
\begin{equation}\label{eq:eigenvalue}
\begin{cases}
\begin{aligned}
-m^{-1}(m\phi')'+\frac{(n-2)^2}{4}\phi&=\lambda\phi,
\quad \text{in }(0,\pi/2),\\
\phi'(\pi/2)&=b_c\phi(\pi/2),\\
\phi'(0)&=0.
\end{aligned}
\end{cases}
\end{equation}
We use the inner product and norm
\[
\langle f,g\rangle
=\int_0^{\pi/2}mfg\,\mathrm d\theta,
\qquad
\|f\|_2^2=\langle f,f\rangle
\]
on $L^2((0,\pi/2),m\,\mathrm d\theta)$. Let $\lambda_0(b_c)$ be the principal eigenvalue of \eqref{eq:eigenvalue} and let $\phi_c>0$
be the corresponding eigenfunction, normalized by $\phi_c(\pi/2)=1$.
Testing the Rayleigh quotient with $v_c$ gives
\[
\begin{aligned}
\lambda_0(b_c)
&\le
\frac{\displaystyle\int_0^{\pi/2}m\left((v_c')^2+\frac{(n-2)^2}{4}v_c^2\right)\,\mathrm d\theta-b_cw_c^2}{\displaystyle\int_0^{\pi/2} mv_c^2\,\mathrm d\theta}=
\frac{(b_0-b_c)w_c^2}{\displaystyle\int_0^{\pi/2}mv_c^2\,\mathrm d\theta}<0.
\end{aligned}
\]
Set 
\begin{equation}\label{eq:prinfre}
\kappa_c=\sqrt{-\lambda_0(b_c)}>0.
\end{equation}
Thus
$\phi_c(\theta)\cos(\kappa_ct)$ and
$\phi_c(\theta)\sin(\kappa_ct)$ solve the linearized equation
and have period $2\pi/\kappa_c$.

\subsection{Sharp trace estimates}
The following alternative follows from Caffarelli--Jin--Sire--Xiong \cite{CJSX}, the classification of regular solutions by Li--Zhu \cite{LiZhu1995} (see also Jin--Li--Xiong \cite{JLX}), and Kelvin inversion.
\begin{theorem}[\cites{LiZhu1995,CJSX}]\label{thm:conformal-end-alternative}
Let $U>0$ be a weak solution of \eqref{eq:halfspace} on
$\overline{\mathbb R^n_+}\setminus\{0\}$. Then exactly one of the following holds.
\begin{enumerate}
\item The origin is removable, and $U$ is a bubble of the form
\eqref{eq:halfspace-bubble}.
\item The origin is nonremovable, and \eqref{eq:radial} holds. 
\end{enumerate}
Consequently, if a positive bounded solution $v$ of
\eqref{eq:cylinder} satisfies
\[
 \lim_{t\to+\infty}v(t,\pi/2)=0
 \quad\text{or}\quad
 \lim_{t\to-\infty}v(t,\pi/2)=0,
\]
then for some $s\in\mathbb R$,
\[
 v(t,\theta)=\bigl(\cosh(t-s)+\cos\theta\bigr)^{-\frac{n-2}{2}}.
\]
\end{theorem}

We next establish a sharp estimate for the boundary trace and its derivative.

\begin{proposition}\label{prop:sharp}
Every positive bounded entire cylindrical solution $v$ of \eqref{eq:cylinder} satisfies
\begin{equation}\label{eq:sharp}
0<w\le1,\qquad |w'|\le \frac{n-2}{2}w\sqrt{1-w^{\frac{4}{n-2}}},\qquad
0<v(t,\theta)\le1,
\end{equation}
where $w(t)=v(t,\pi/2)$.
\end{proposition}

\begin{proof}
For a bubble, $w(t)=\cosh(t-\tau)^{-\frac{n-2}{2}}$, so
$|w'|=\frac{n-2}{2}w\sqrt{1-w^{\frac{4}{n-2}}}$ and
$v(t,\theta)=\bigl(\cosh(t-\tau)+\cos\theta\bigr)^{-\frac{n-2}{2}}\le1$.
For a nonbubble, use \eqref{eq:radial}.
Set $r=|x|$ and
\[
\Phi=U^{-\frac{2}{n-2}},\qquad
Q=\Phi_r^2+\Phi_y^2-\frac{2\Phi\Phi_r}{r}=|\nabla\Phi|^2-\frac{2\Phi\Phi_r}{r}.
\]
The quotient $\Phi_r/r$ extends smoothly to the symmetry axis $\{r=0,\ y>0\}$, with value $\Phi_{rr}(0,y)$ at $r=0$. And $Q$ is also even in \(r\). Equation \eqref{eq:halfspace} gives
\begin{equation}\label{eq:phi-equation}
\begin{cases}
\displaystyle
\Phi\left(\Phi_{rr}+\frac{n-2}{r}\Phi_r+\Phi_{yy}\right)
=\frac n2|\nabla\Phi|^2
& \text{for }r>0,\ y>0,\\[1mm]
\Phi_y(r,0)=\dfrac{2}{n-2}
& \text{for }r>0.
\end{cases}
\end{equation}
In the $(r,y)$ variables, set 
\[
M=D^2\Phi-\frac{\Phi_r}{r}I\quad \mbox{and}\quad F=(\Phi_r-\Phi/r,\Phi_y).
\]
Differentiation yields
\begin{equation}\label{eq:Q-first}
\operatorname{tr} M=\frac{nQ}{2\Phi},\qquad
\nabla Q=2MF,\qquad |F|^2=Q+\frac{\Phi^2}{r^2}.
\end{equation}
Set
\[
L=\partial_{rr}+\partial_{yy}+\frac nr\partial_r
-\frac n\Phi\nabla\Phi\cdot\nabla.
\]
From \eqref{eq:phi-equation},
\[
\begin{aligned}
L\Phi&=\frac{2\Phi_r}{r}-\frac{n}{2\Phi}|\nabla\Phi|^2,\\
L\left(\frac{\Phi_r}{r}\right)&=\frac{n}{\Phi r}\bigl(\Phi_{rr}\Phi_r+\Phi_{ry}\Phi_y\bigr)
-\frac{n\Phi_r|\nabla\Phi|^2}{2r\Phi^2}-\frac n\Phi\nabla\Phi\cdot\nabla\left(\frac{\Phi_r}{r}\right),\\
L\bigl(|\nabla\Phi|^2\bigr)&=2|D^2\Phi|^2-\frac n{\Phi^2}|\nabla\Phi|^4
-2(n-2)\nabla\Phi\cdot\nabla\left(\frac{\Phi_r}{r}\right)
+\frac{2n}{r}\bigl(\Phi_{rr}\Phi_r+\Phi_{ry}\Phi_y\bigr).
\end{aligned}
\]
Applying
the product rule to $Q=|\nabla\Phi|^2-\frac{2\Phi\Phi_r}{r}$ gives
\begin{equation}\label{eq:Q-equation}
\begin{aligned}
LQ
&=2|D^2\Phi|^2-\frac{4\Phi_r^2}{r^2}-\frac n{\Phi^2}|\nabla\Phi|^4
+\frac{2n\Phi_r}{r\Phi}|\nabla\Phi|^2=2|M|^2-\frac n{\Phi^2}Q^2.
\end{aligned}
\end{equation}

We claim that $Q\le0$. For $\lambda>0$, set
$\Phi_\lambda(X)=\lambda^{-1}\Phi(\lambda X)$ and $Q_\lambda(X)=Q(\lambda X)$.
By \eqref{eq:radial} and elliptic estimates,
$\Phi_\lambda$ has uniform positive lower bounds and uniform $C^j$ bounds
on fixed closed half-annuli, including the symmetry axis.
Thus $Q$ is bounded throughout \(\overline{\mathbb R^n_+}\setminus\{0\}\). If $A=\sup Q>0$, choose $X_j$ with $Q(X_j)\to A$
and take $\lambda_j=|X_j|$. A subsequence of $\Phi_{\lambda_j}$ converges
in $C^4$ on compact half-annuli, and $X_j/|X_j|$ converges to $X_0$ with
$|X_0|=1$. The limiting function satisfies $Q\le A$ and $Q(X_0)=A$.

If $r(X_0)>0$ and $y(X_0)>0$, then $MF=0$ and $F\ne0$ by
\eqref{eq:Q-first}. Hence, $0$ is an eigenvalue of $M$, and the other eigenvalue must be $\operatorname{tr} M$. Hence
\[
0\ge LQ
=2\left(\frac{nQ}{2\Phi}\right)^2-\frac n{\Phi^2}Q^2
=\frac{n(n-2)}{2\Phi^2}Q^2>0,
\]
a contradiction. On the symmetry axis,
$M=\operatorname{diag}(0,\Phi_{yy}-\Phi_{rr})$ and $Q_r(0,y)=0$.
Thus $r^{-1}Q_r\to Q_{rr}$ as $r\downarrow0$, and at a maximum
\eqref{eq:Q-equation} gives
\[
 0\ge(n+1)Q_{rr}+Q_{yy}=LQ
 =\frac{n(n-2)}{2\Phi^2}Q^2>0,
\]
again a contradiction.

If $y(X_0)=0$, then $r(X_0)>0$.  Write
$M_{rr}=\Phi_{rr}-\Phi_r/r$ and $F_r=\Phi_r-\Phi/r$. By
\eqref{eq:phi-equation}, $\Phi_{ry}=\Phi_{rry}=0$ on this boundary, and by \eqref{eq:Q-first}, 
$0=Q_r=2(\Phi_{rr}-\Phi_r/r)(\Phi_r-\Phi/r)$ at $X_0$. If $\Phi_r-\Phi/r\ne0$, then $\Phi_{rr}-\Phi_r/r=0$. If $\Phi_r-\Phi/r=0$, then
\[
0\ge Q_{rr}=2\left(\Phi_{rr}-\frac{\Phi_r}{r}\right)\left(\Phi_{rr}-\frac{\Phi_r}{r}+\frac\Phi{r^2}\right),
\]
so $\Phi_{rr}-\Phi_r/r\le0$. In both cases,
\[
Q_y=\frac{4}{n-2}\left[\operatorname{tr} M-\left(\Phi_{rr}-\frac{\Phi_r}{r}\right)\right]=\frac{4}{n-2}\left[\frac{nQ}{2\Phi}-\left(\Phi_{rr}-\frac{\Phi_r}{r}\right)\right]>0.
\]
This contradicts $Q_y(X_0)\le0$. Hence $Q\le0$.

On $y=0$, the identities $w=\left(\frac{2r}{(n-2)\Phi}\right)^{\frac{n-2}{2}}$ and $t=-\log r$ therefore give
\[
\begin{aligned}
\left(\frac{2w'}{(n-2)w}\right)^2+w^{\frac{4}{n-2}}
&=\left(\frac{r\Phi_r}{\Phi}-1\right)^2
+\frac{4r^2}{(n-2)^2\Phi^2}=1+\frac{r^2Q}{\Phi^2}\le1.
\end{aligned}
\]
This proves $0<w\le 1$. Since $v$ is the bounded Dirichlet extension
of $w$ for the shifted cylindrical operator, comparison gives
$0<v\le1$. For completeness, fix $\varepsilon>0$ and compare \(v\) with the function $1+\varepsilon\cosh\!\left(\frac{n-2}{2}t\right)$. This is a supersolution. It dominates \(v\) on the equator, and also for large $|t|$ because \(v\) is bounded. The maximum principle, followed by \(\varepsilon\downarrow0\), gives \(v\le1\).
\end{proof}

\begin{lemma}\label{lem:secant}
Let $v_1,v_2\ge0$ be bounded solutions. Set
$w_i(t)=v_i(t,\pi/2)$, $D=v_1-v_2$, and
\begin{equation}\label{eq:secant}
 b(t)=\frac n2\int_0^1\bigl((1-s)w_2(t)+sw_1(t)\bigr)^{\frac{2}{n-2}}\,\mathrm ds.
\end{equation}
Then $D$ satisfies
\begin{equation}\label{eq:secant-equation}
\begin{cases}
\begin{aligned}
D_{tt}+m^{-1}(mD_\theta)_\theta-\frac{(n-2)^2}{4}D&=0 \quad&& \text{in }\mathbb R\times(0,\pi/2),\\
D_\theta(t,\pi/2)&=b(t)D(t,\pi/2)\quad && \text{for }t\in\mathbb R,\\
D_\theta(t,0)&=0\quad && \text{for }t\in\mathbb R.
\end{aligned}
\end{cases}
\end{equation}
Moreover,
\[
0\le b\le \frac n2,\qquad
|b'|\le b\sqrt{1-\left(\frac{2b}{n}\right)^2}\le b.
\]
If $(v_1,v_2)\ne(0,0)$, then $b>0$ on $\mathbb R$. If $v_1=v_2=v$, then $b=\frac n2w^{\frac{2}{n-2}}$.
\end{lemma}

\begin{proof}
The function $x\mapsto x\sqrt{1-x^{\frac{4}{n-2}}}$ is strictly concave on $(0,1)$, since its second  derivative is negative on $(0,1)$. Let $w_s=(1-s)w_2+sw_1$ and $z_s=w_s^{\frac{2}{n-2}}$. By
Proposition~\ref{prop:sharp} and concavity,
\[
\begin{aligned}
|w_s'|
&\le \frac{n-2}{2}\bigl((1-s)w_2\sqrt{1-w_2^{\frac{4}{n-2}}}+sw_1\sqrt{1-w_1^{\frac{4}{n-2}}}\bigr)\le \frac{n-2}{2}w_s\sqrt{1-w_s^{\frac{4}{n-2}}},\\
|z_s'|&\le z_s\sqrt{1-z_s^2}.
\end{aligned}
\]
Since $z\mapsto z\sqrt{1-z^2}$ is concave on $[0,1]$, and $b(t)=\frac n2\int_0^1z_s(t)\,\mathrm ds$, we obtain
\[
\frac{2|b'|}{n}=\left|\int_0^1z_s'\,\mathrm ds\right|
\le\int_0^1 z_s\sqrt{1-z_s^2}\,\mathrm{d} s
\le \frac{2b}{n}\sqrt{1-\left(\frac{2b}{n}\right)^2},
\]
where we used Jensen's inequality in the last inequality. 
\end{proof}

\subsection{Robin spectral gap and the principal eigenspace}
We next establish spectral estimates for the linear equations arising
in Lemma~\ref{lem:secant}. The separation of the principal eigenvalue
from the remaining spectrum is the starting point for the
two-dimensionality of the weighted solution space.

Consider Robin coefficients satisfying
\begin{equation}\label{eq:coefficient-class}
b\in W^{1,\infty}(\mathbb{R}),\qquad
0\le b\le \frac n2,\qquad |b'|\le b\quad\text{almost everywhere}.
\end{equation}
Every $b(t)$ in Lemma~\ref{lem:secant} satisfies these conditions.
By Gronwall's inequality, either $b>0$ on $\mathbb R$ or $b\equiv0$.

Recall that $\langle \cdot,\cdot\rangle$ and $\|\cdot\|_2$ are the real inner product and norm in $L^2((0,\pi/2),m\,\mathrm{d}\theta)$, respectively. Constants denoted by $C_n$ depend only on $n$.

Let $H^1_{\mathrm{ax}}$ be the closure of smooth regular axially symmetric functions on $\mathbb{S}^{n-1}_+$ under the norm
\[
\|f\|_{H^1_{\mathrm{ax}}}^2
=\int_0^{\pi/2} m(|f'|^2+|f|^2)\,\mathrm{d}\theta.
\]
We identify $H^1_{\mathrm{ax}}$ with the axially symmetric subspace of
$H^1(\mathbb S^{n-1}_+)$, suppressing the factor $|\mathbb S^{n-2}|$.
The map \(f\mapsto f(\pi/2)\) is continuous from \(H^1_{\mathrm{ax}}\) to \(\mathbb R\). Write $H^{-1}_{\mathrm{ax}}$ for its dual.
The closed symmetric form
\[
Q_b[f,\phi]
=\int_0^{\pi/2} m(f'\phi'+\frac{(n-2)^2}{4}f\phi)\,\mathrm{d}\theta-bf(\pi/2)\phi(\pi/2),
\]
with form domain $H^1_{\mathrm{ax}}$, defines the self-adjoint axial Robin operator with compact resolvent
\[
A_bf=-m^{-1}(mf')'+\frac{(n-2)^2}{4}f,\qquad f'(\pi/2)=bf(\pi/2),
\]
with regularity at the pole interpreted on the hemisphere.  Using the trace inequality
\begin{equation}\label{eq:trace}
|f(\pi/2)|^2
\le \varepsilon\int_0^{\pi/2}m|f'|^2\,\mathrm d\theta
+C_{\varepsilon}\|f\|_2^2,
\qquad \varepsilon>0,
\end{equation}
we have
\[
\frac12\|f\|_{H^1_{\mathrm{ax}}}^2
\le Q_b[f,f]+C_n\|f\|_2^2
\le C_n'\|f\|_{H^1_{\mathrm{ax}}}^2.
\]

We call $D$ a weak solution of 
\begin{equation}\label{eq:homoeq}
D_{tt}-A_{b(t)}D=0
\end{equation}  
if
\[
D\in L^2_{\mathrm{loc}}(\mathbb{R};H^1_{\mathrm{ax}}),\quad
D_t\in L^2_{\mathrm{loc}}(\mathbb{R};L^2(m\,\mathrm{d}\theta)),\quad
D_{tt}=A_{b(t)}D\quad\text{in } L^2_{\mathrm{loc}}(\mathbb{R};H^{-1}_{\mathrm{ax}}),
\]
where the time derivatives are understood in the sense of distributions.

We first record the local regularity needed to define the Cauchy
data at each section. For fixed nested intervals $I_1\Subset I_2$, testing with a time cutoff and using
\[
 \|f\|_{L^2(\partial\mathbb{S}^{n-1}_+)}^2
 \le\varepsilon\|\nabla f\|_{L^2(\mathbb{S}^{n-1}_+)}^2
 +C_\varepsilon\|f\|_{L^2(\mathbb{S}^{n-1}_+)}^2
\]
gives, uniformly over \eqref{eq:coefficient-class},
\begin{equation}
 \|D\|_{H^1(I_1\times\mathbb S^{n-1}_+)}
 \le C_{n,I_1,I_2}\|D\|_{L^2(I_2\times\mathbb S^{n-1}_+)}.
\end{equation}
Choose $I_1\Subset J\Subset I_2$. On the cylinder, $D$ satisfies
\[
\begin{cases}
-D_{tt}-\Delta_{\mathbb S^{n-1}}D+\frac{(n-2)^2}{4}D=0
&\text{in }J\times\mathbb S^{n-1}_+,\\
\partial_\nu D=bD
&\text{on }J\times\partial\mathbb S^{n-1}_+,
\end{cases}
\]
where $\nu$ is the outward unit normal. Since
$\|b\|_{L^\infty}+\|b'\|_{L^\infty}\le n$, the product rule and
interpolation between $L^2$ and $H^1$, followed by the trace theorem, give
\[
\begin{aligned}
\|bD\|_{H^{1/2}(J\times\partial\mathbb S^{n-1}_+)}
&\le C_{n,J}\|D\|_{H^{1/2}(J\times\partial\mathbb S^{n-1}_+)}\le C_{n,J}\|D\|_{H^1(J\times\mathbb S^{n-1}_+)}.
\end{aligned}
\]
Applying the local Neumann estimate
\cite{McLean}*{Theorem~4.18(ii), $r=0$} in smooth boundary charts and
interior regularity, then using the preceding $H^1$ estimate, we obtain
\begin{equation}\label{eq:local-H2}
\begin{aligned}
\|D\|_{H^2(I_1\times\mathbb S^{n-1}_+)}
&\le C_{n,I_1,J}\Bigl(\|D\|_{H^1(J\times\mathbb S^{n-1}_+)}
+\|bD\|_{H^{1/2}(J\times\partial\mathbb S^{n-1}_+)}\Bigr)\\
&\le C_{n,I_1,J}\|D\|_{H^1(J\times\mathbb S^{n-1}_+)}\\
&\le C_{n,I_1,I_2}\|D\|_{L^2(I_2\times\mathbb S^{n-1}_+)}.
\end{aligned}
\end{equation}
Hence $D(t_0)\in H^{3/2}(\mathbb S^{n-1}_+)$ and
$D_t(t_0)\in H^{1/2}(\mathbb S^{n-1}_+)$ for every \(t_0\in\mathbb R\).

Let $e_0(b)>0$ be the principal eigenfunction of $A_b$,
normalized by $\|e_0(b)\|_2=1$, and set $e(t)=e_0(b(t))$.

\begin{lemma}\label{lem:robin-spectral-estimates}
For $0\le b\le n/2$, the eigenvalues of $A_b$ satisfy
\[
-\frac{(n-1)(n+3)}4
\le\lambda_0(b)\le\frac{(n-2)^2}{4}
<\frac{n^2}{4}<\lambda_1(b),
\]
with $\lambda_0(b)<(n-2)^2/4$ for $b>0$.
Moreover, there is $c_n>0$, depending only on $n$, such that
\[
(Q_b-\lambda_0(b))[\phi,\phi]\ge c_n\|\phi\|_{H^1_{\mathrm{ax}}}^2,
\qquad \phi\in H^1_{\mathrm{ax}},\quad \phi\perp e_0(b).
\]
\end{lemma}

\begin{proof}
For comparison with the Robin problem, consider the axial eigenvalue equation
\[
-m^{-1}(mf')'+\frac{(n-2)^2}{4}f=\lambda f,\qquad 0<\theta<\pi/2,
\]
with smooth rotational extension at $\theta=0$. The Neumann and Dirichlet
problems impose $f'(\pi/2)=0$ and $f(\pi/2)=0$, respectively, in place of
the Robin condition $f'(\pi/2)=bf(\pi/2)$. The function \(f(\theta)=\cos\theta\) is a positive Dirichlet eigenfunction with eigenvalue \(n^2/4\). Hence the first Dirichlet eigenvalue is \(n^2/4\). Testing the Rayleigh quotient with $1$ gives
\[
\lambda_0(b)\le\frac{Q_b[1,1]}{\|1\|_2^2}
=\frac{(n-2)^2}{4}
-\left(\displaystyle\int_0^{\pi/2}m\,\mathrm d\theta\right)^{-1} b
\le\frac{(n-2)^2}{4},
\]
with strict inequality in the last step for $b>0$.
Let $e_0(b),e_1(b)$ be orthonormal eigenfunctions for
$\lambda_0(b),\lambda_1(b)$. Since their span is two-dimensional, choose
$0\ne f\in\operatorname{span}\{e_0(b),e_1(b)\}$ with $f(\pi/2)=0$.
The Dirichlet variational characterization gives
\begin{equation}\label{eq:Dirifirst}
\frac{n^2}{4}
=\min_{\substack{0\ne g\in H^1_{\mathrm{ax}}\\g(\pi/2)=0}}
\frac{Q_b[g,g]}{\|g\|_2^2}
\le\frac{Q_b[f,f]}{\|f\|_2^2}
\le\lambda_1(b).
\end{equation}

If $\lambda_1(b)=n^2/4$, then
\[
0=\lambda_1(b)\|f\|_2^2-Q_b[f,f]
=\bigl(\lambda_1(b)-\lambda_0(b)\bigr)
|\langle f,e_0(b)\rangle|^2.
\]
Since $\lambda_0(b)\le(n-2)^2/4<n^2/4$, we obtain
$f\in\operatorname{span}\{e_1(b)\}$ and hence
$f'(\pi/2)=bf(\pi/2)=0$. Uniqueness for the eigenvalue ODE with
$f(\pi/2)=f'(\pi/2)=0$ gives $f\equiv0$, a contradiction. This proves \(\lambda_1(b)>n^2/4\).

Completing a square in the weighted form gives
\[
\begin{aligned}
Q_b[f,f]
&=\int_0^{\pi/2} m\,|f'-b\sin\theta\,f|^2\,\mathrm{d}\theta\\
&\quad+\int_0^{\pi/2} m\bigl[\frac{(n-2)^2}{4}-b^2\sin^2\theta-b(n-1)\cos\theta\bigr]f^2\,\mathrm{d}\theta.
\end{aligned}
\]
Here $(\sin^{n-1}\theta)'=(n-1)\sin^{n-2}\theta\cos\theta$ and
$\sin^{n-1}\theta\,f(\theta)^2\to0$ as $\theta\downarrow0$.
Thus
\[
\begin{aligned}
 \frac{(n-2)^2}{4}-b^2\sin^2\theta-b(n-1)\cos\theta
 &=\frac{(n-2)^2}{4}-b^2+\left(b\cos\theta-\frac{n-1}{2}\right)^2-\frac{(n-1)^2}{4}\\
 &\ge \frac{(n-2)^2}{4}-\frac{n^2}{4}-\frac{(n-1)^2}{4}.
\end{aligned}
\]
Hence
\begin{equation}\label{eq:lambda-lower}
 \lambda_0(b)\ge-\frac{(n-1)(n+3)}4\quad\mbox{for every }n\ge3.
\end{equation}

If $\phi\in H^1_{\mathrm{ax}}$ and $\phi\perp e_0(b)$, then
\[
(Q_b-\lambda_0(b))[\phi,\phi]
\ge(\lambda_1(b)-\lambda_0(b))\|\phi\|_2^2
\ge(n-1)\|\phi\|_2^2.
\]
The trace inequality \eqref{eq:trace} also gives, uniformly for $0\le b\le \frac n2$,
\[
(Q_b-\lambda_0(b))[\phi,\phi]
\ge\frac12\int_0^{\pi/2}m|\phi'|^2\,\mathrm d\theta-C_n\|\phi\|_2^2.
\]
These two inequalities imply
\[
(Q_b-\lambda_0(b))[\phi,\phi]\ge c_n\|\phi\|_{H^1_{\mathrm{ax}}}^2,
\qquad \phi\perp e_0(b).
\]
\end{proof}

Write $\lambda_0=\lambda_0(b)$. Choose a real orthonormal eigenbasis
$e_j(b)$, with $e_0(b)>0$, and set $e_j^\partial(b)=e_j(b,\pi/2)$.

\begin{lemma}\label{lem:principal-eigenfunction-variation}
The map $b\mapsto e_0(b)$ is differentiable in $H^1_{\mathrm{ax}}$ for
$0\le b\le n/2$, and
\begin{equation}\label{eq:e0H1}
\sup_{0\le b\le n/2}
\left(\|e_0(b)\|_{H^1_{\mathrm{ax}}}
+\|\partial_be_0(b)\|_{H^1_{\mathrm{ax}}}\right)\le C_n.
\end{equation}
Moreover,
\[
\begin{gathered}
\partial_b\lambda_0=-(e_0^\partial)^2,\ \ 
\langle\partial_be_0,e_0\rangle=0,\ \ 
\langle\partial_be_0,e_j\rangle
=\frac{e_0^\partial e_j^\partial}{\lambda_j-\lambda_0}
\ \ (j\ge1),\ \ 
\|\partial_be_0\|_2
\le\frac{(e_0^\partial)^2}{2(\frac{n^2}{4}-\lambda_0)}.
\end{gathered}
\]
If $b$ satisfies \eqref{eq:coefficient-class} and $e(t)=e_0(b(t))$, then for almost every $t\in\mathbb R$,
\begin{equation}\label{eq:rotation}
\|e'(t)\|_2
\le\frac{|b'(t)|[e_0^\partial(b(t))]^2}
{2[\frac{n^2}{4}-\lambda_0(b(t))]}.
\end{equation}
\end{lemma}

\begin{proof}
By $\|e_0\|_2=1$ and Lemma~\ref{lem:robin-spectral-estimates}, we have
\[
\frac12\|e_0\|_{H^1_{\mathrm{ax}}}^2
\le Q_b[e_0,e_0]+C_n\|e_0\|_2^2
=\lambda_0+C_n\le\frac{(n-2)^2}{4}+C_n.
\]
The forms $Q_b$ depend analytically on $b$ on the common domain
$H^1_{\mathrm{ax}}$. Since $\lambda_0(b)$ is simple, $b\mapsto e_0(b)$
is differentiable in $H^1_{\mathrm{ax}}$; see
Kato~\cite{Kato}*{Chapter~VII, \S4}. Differentiating $\|e_0\|_2^2=1$ and
$Q_b[e_0,\phi]=\lambda_0\langle e_0,\phi\rangle$ gives $\langle\partial_be_0,e_0\rangle=0$ and
\[
(Q_b-\lambda_0)[\partial_be_0,\phi]
=e_0^\partial\phi(\pi/2)
+(\partial_b\lambda_0)\langle e_0,\phi\rangle,
\qquad \phi\in H^1_{\mathrm{ax}}.
\]
Taking $\phi=e_j$, $j\ge0$, gives the stated differentiation identities.
For $\phi\perp e_0$, the above identity reduces to
\[
(Q_b-\lambda_0)[\partial_be_0,\phi]=e_0^\partial\phi(\pi/2).
\]
By Lemma~\ref{lem:robin-spectral-estimates} and trace continuity,
\[
\begin{aligned}
c_n\|\partial_be_0\|_{H^1_{\mathrm{ax}}}^2
\le(Q_b-\lambda_0)[\partial_be_0,\partial_be_0]
=e_0^\partial(\partial_be_0)(\pi/2)&\le C_n\|e_0\|_{H^1_{\mathrm{ax}}}
\|\partial_be_0\|_{H^1_{\mathrm{ax}}}\\
&\le C_n\|\partial_be_0\|_{H^1_{\mathrm{ax}}}.
\end{aligned}
\]
This proves \eqref{eq:e0H1}.

To obtain the explicit $L^2$ bound for $\partial_be_0$, we evaluate
the Robin-to-Dirichlet map at the first Dirichlet eigenvalue $n^2/4$. For $\zeta\notin\sigma(A_b)$, the scalar Robin-to-Dirichlet map is
\begin{equation}\label{eq:specpara}
R_b(\zeta)
=\sum_{j\ge0}\frac{(e_j^\partial)^2}{\lambda_j-\zeta}.
\end{equation}
Equivalently, $R_b(\zeta)=u(\pi/2)$, where $u$ is regular at $\theta=0$
and satisfies
\[
\begin{gathered}
-m^{-1}(mu')'+\left(\frac{(n-2)^2}{4}-\zeta\right)u=0,
\qquad u'(\pi/2)-bu(\pi/2)=1,\\
(Q_b-\zeta)[u,\phi]=\phi(\pi/2),\qquad \phi\in H^1_{\mathrm{ax}}.
\end{gathered}
\]
At $\zeta=n^2/4$, the regular solution is $u=-\cos\theta$, so
$u(\pi/2)=0$ and $u'(\pi/2)=1$. Thus $R_b(n^2/4)=0$, and
\eqref{eq:specpara} gives
\begin{equation}\label{eq:eigenres}
\sum_{j\ge1}\frac{(e_j^\partial)^2}{\lambda_j-\frac{n^2}{4}}
=\frac{(e_0^\partial)^2}{\frac{n^2}{4}-\lambda_0}.
\end{equation}
Since $\lambda_j>n^2/4>\lambda_0$ and
$(\lambda_j-\lambda_0)^2\ge
4(\lambda_j-\frac{n^2}{4})(\frac{n^2}{4}-\lambda_0)$ for $j\ge1$,
\[
\begin{aligned}
\|\partial_be_0\|_2^2
&=\sum_{j\ge1}
\frac{(e_0^\partial)^2(e_j^\partial)^2}{(\lambda_j-\lambda_0)^2}\le\frac{(e_0^\partial)^2}{4(\frac{n^2}{4}-\lambda_0)}
\sum_{j\ge1}\frac{(e_j^\partial)^2}{\lambda_j-\frac{n^2}{4}}
=\frac{(e_0^\partial)^4}{4(\frac{n^2}{4}-\lambda_0)^2}.
\end{aligned}
\]
Finally, $e'(t)=b'(t)\partial_be_0(b(t))$ almost everywhere gives
\eqref{eq:rotation}.
\end{proof}

To sharpen \eqref{eq:rotation}, we use $|b'|\le b$ and estimate
$b(e_0^\partial)^2$ through the Fourier symbol of the cylindrical
Dirichlet-to-Neumann operator.
For a smooth bounded boundary function $w$, let $v$ be the unique
bounded solution of
\[
\begin{cases}
v_{tt}+m^{-1}(mv_\theta)_\theta-\frac{(n-2)^2}{4}v=0
& \text{in }\mathbb R\times(0,\pi/2),\\
v(t,\pi/2)=w(t)
& \text{for }t\in\mathbb R,
\end{cases}
\]
regular at $\theta=0$, and define
$$\mathcal P w:=v_\theta(\cdot,\pi/2).$$
For $w(t)=e^{i\xi t}$, $\xi\in\mathbb R$, the extension has the form
$v(t,\theta)=e^{i\xi t}\psi_\xi(\theta)$, where the regular factor
$\psi_\xi$ satisfies
\begin{equation}\label{eq:modexi}
\begin{cases}
m^{-1}(m\psi_\xi')' -\left(\frac{(n-2)^2}{4}+\xi^2\right)\psi_\xi=0
& 0<\theta<\pi/2,\\
\psi_\xi'(0)=0,\qquad \psi_\xi(\pi/2)=1.
\end{cases}
\end{equation}
Consequently,
$\widehat{\mathcal P w}(\xi)=\Theta(\xi)\widehat w(\xi)$,
where $\Theta(\xi)=\psi_\xi'(\pi/2)$ and the Fourier transforms
are taken in $t$, in the sense of tempered distributions. The formula of
DelaTorre--Gonz\'alez~\cite{DelaTorreGonzalez2018}*{Theorem~1.1},
applied in dimension $n-1$ with fractional order $1/2$ and
constant angular mode, gives
\begin{equation}\label{eq:symbol}
\Theta(\xi)
=2\frac{
\Gamma\left(\frac{n+2i\xi}{4}\right)
\Gamma\left(\frac{n-2i\xi}{4}\right)}
{\Gamma\left(\frac{n-2+2i\xi}{4}\right)
\Gamma\left(\frac{n-2-2i\xi}{4}\right)}.
\end{equation}
Stirling's formula gives
$\Theta(\xi)=|\xi|+O(|\xi|^{-1})$ as $|\xi|\to\infty$. If $w$ is the trace of a bounded solution of \eqref{eq:cylinder},
the boundary condition gives
\[
\mathcal P w=\frac{n-2}{2}w^{\frac{n}{n-2}}.
\]

We now use \eqref{eq:rotation} and \eqref{eq:symbol} to prove the following
estimate for $\|e'(t)\|_2$.
\begin{lemma}\label{lem:quantitative-gap}
Let $b$ satisfy \eqref{eq:coefficient-class}, and set $\lambda(t)=\lambda_0(b(t))$, $e(t)=e_0(b(t))$ and
$$\eta(t)=\frac12\left(\frac n2-\frac{|\frac{n^2}{4}+\lambda(t)|}{n}\right).$$
There is a constant $\varepsilon_n<1$, depending only on $n$, such that for
almost every $t\in\mathbb R$,
\begin{equation}\label{eq:strict-gap}
 \eta\ge\frac13,\qquad \|e'\|_2\le\varepsilon_n \eta.
\end{equation}
\end{lemma}
\begin{proof}
Taking $b=0$ in \eqref{eq:specpara}, we obtain the scalar
Neumann-to-Dirichlet map $R_0(\zeta)$. Then $R_0(\lambda)=\Theta(\sqrt{-\lambda})^{-1}$ if  $\lambda\le 0$. Both $R_0(\lambda)$ and $\Theta(i\sqrt\lambda)^{-1}$ ($i=\sqrt{-1}$) are holomorphic
for $|\lambda|<\frac{(n-2)^2}{4}$
(the gamma quotient is even in $\sqrt\lambda$).
By analytic continuation,
$R_0(\lambda)=\Theta(i\sqrt\lambda)^{-1}$ also holds for
$0<\lambda<\frac{(n-2)^2}{4}$. Substituting \eqref{eq:symbol} and using the Weierstrass product formula for the gamma function $\Gamma$, we obtain
\[
\frac{R_0(\lambda)}{R_0(0)}
=\prod_{j\ge0}
 \frac{1-4\lambda/(n+4j)^2}
      {1-4\lambda/(n-2+4j)^2}
\quad \mbox{for every }\lambda<\frac{(n-2)^2}{4}.
\]
The infinite product is positive, and both the paired product and its differentiated
logarithmic series converge locally uniformly. Hence
\begin{equation}\label{eq:log-resolvent}
\begin{aligned}
\Lambda(\lambda):=\frac{\mathrm d}{\mathrm d \lambda}\log R_0(\lambda)
&=\sum_{j\ge0}\left(
 \frac4{(n-2+4j)^2-4\lambda}
-\frac4{(n+4j)^2-4\lambda}\right).
\end{aligned}
\end{equation}
For $b>0$, differentiating $R_0(\lambda_0(b))=1/b$ and using Lemma \ref{lem:principal-eigenfunction-variation}, one obtains
\begin{equation}\label{eq:boundary-mass}
 b[e_0^\partial(b)]^2=\frac1{\Lambda(\lambda_0(b))}.
\end{equation}

Set $\delta=\frac{(n-2)^2}{4}-\lambda_0(b)$. By \eqref{eq:lambda-lower},
$0\le\delta\le\delta_{\max}:=\frac{2n^2-2n+1}{4}$.
The $j=0$ term in \eqref{eq:log-resolvent} for $\Lambda(\lambda_0(b))$ is
$(n-1)/[\delta(\delta+n-1)]$.  For the remaining terms write
\[
 \frac4{(n-2+4j)^2-4\lambda_0(b)}-\frac4{(n+4j)^2-4\lambda_0(b)}
 =\int_{\frac{n-2}{2}+2j}^{\frac n2+2j}\frac{2s}{(s^2-\lambda_0(b))^2}\,\mathrm{d} s.
\]
For $s\ge \frac{n+2}{2}$, \eqref{eq:lambda-lower} gives
$3\left(\frac{n+2}{2}\right)^2+\lambda_0(b)\ge\frac{2n^2+10n+15}{4}>0$. Hence
\[
 \frac{\mathrm{d}}{\mathrm{d} s}\frac{2s}{(s^2-\lambda_0(b))^2}
 =-\frac{2(3s^2+\lambda_0(b))}{(s^2-\lambda_0(b))^3}<0.
\]
Since $[\frac{n-2}{2}+2j,\frac n2+2j]$ is the first half of $[\frac{n-2}{2}+2j,\frac{n+2}{2}+2j]$, summation gives
\begin{equation}\label{eq:log-resolvent-strong}
 \Lambda(\lambda_0(b))\ge
 \frac{n-1}{\delta(\delta+n-1)}+
 \frac1{2(\delta+2n)}\quad\mbox{for }\delta>0.
\end{equation}
Recall that either \(b\equiv0\) or \(b>0\) everywhere. For $b\equiv0$, $\lambda_0(0)=(n-2)^2/4$, $\eta=(n-1)/(2n)\ge 1/3$, and $\|e'\|_2=0$. For $b>0$, we know from Lemma \ref{lem:robin-spectral-estimates} that $\delta>0$.
Combining \eqref{eq:rotation}, $|b'|\le b$, and
\eqref{eq:boundary-mass} gives
\begin{equation*}
 \|e'\|_2\le\frac1{2(\delta+n-1)\Lambda(\lambda_0(b))}.
\end{equation*}

Set $\delta_0=\frac{n^2-2n+2}{2}$. If $0<\delta\le\delta_0$, then $\frac{n^2}{4}+\lambda_0(b)\ge0$ and
\[
 \eta=\frac{\delta+n-1}{2n}\ge\frac{n-1}{2n}\ge\frac13,\qquad
 \|e'\|_2\le\frac{\delta}{2(n-1)}
\]
by the first term of \eqref{eq:log-resolvent-strong}. Thus
\[
 \frac{\|e'\|_2}{\eta}\le\frac{n\delta}{(n-1)(\delta+n-1)}
 \le\frac{n\delta_0}{(n-1)(\delta_0+n-1)}<1,
\]
since $(n-1)(\delta_0+n-1)-n\delta_0=\frac{n(n-2)}2>0$.

If $\delta_0\le\delta\le\delta_{\max}$, then $\frac{n^2}{4}+\lambda_0(b)\le0$ and
\[
 \eta=\frac{n^2-n+1-\delta}{2n},\qquad
 \|e'\|_2\le
 \frac{\delta(\delta+2n)}
 {\delta^2+3(n-1)\delta+4n(n-1)}:=G_n(\delta).
\]
 Then $\eta\ge (n^2-n+1-\delta_{\max})/(2n)=(2n^2-2n+3)/(8n)>1/3$ for all $n\ge 3$, and $\eta > 1$ for all $n\ge 5$. We observe that $\delta^2+3(n-1)\delta+4n(n-1)-\delta(\delta+2n)=(n-3)\delta+4n(n-1)>0$, and thus,
\[
\|e'\|_2\le G_n(\delta)<1<\eta\quad\mbox{for } n\ge 5.
\]
Furthermore, a direct differentiation gives $G_n'(\delta)>0$. Consequently, for \(n=3\) and \(n=4\), respectively, we have $G_3(\delta)\le G_3(13/4)=\frac{481}{865}<\frac58\le\eta$ and $G_4(\delta)\le G_4(25/4)=\frac{1425}{2293}<\frac{27}{32}
\le\eta.$ This proves $\|e'\|_2<\eta$.

The upper bounds for $\|e'\|_2/\eta$ are continuous and
strictly less than $1$ on the two compact intervals $[0,\delta_0]$ and $[\delta_0,\delta_{\max}]$. Their maximum gives $\varepsilon_n<1$ in \eqref{eq:strict-gap}.
\end{proof}

\subsection{Two-dimensionality of the linear equation}
Let 
\[
\Pi_0(b)f=\langle e_0(b),f\rangle e_0(b).
\]
Since \(\|e_0(b)\|_2=1\), this is the orthogonal projection onto \(\operatorname{span}\{e_0(b)\}\). For fixed $b_{\mathrm{ref}}\in[0,\frac n2]$,
define Kato's parallel transport \cite{Kato}*{Chapter~II, \S4} by
\begin{equation}\label{eq:mathcalT}
\partial_b\mathcal T=[\partial_b\Pi_0,\Pi_0]\mathcal T,
\qquad \mathcal T(b_{\mathrm{ref}})=I,
\end{equation}
where
\[
[\partial_b\Pi_0,\Pi_0]f
=\langle e_0,f\rangle\,\partial_be_0
 -\langle\partial_be_0,f\rangle\,e_0.
\]
The generator is skew-adjoint in $L^2$ and, together with its adjoint,
uniformly bounded on $H^1_{\mathrm{ax}}$. Thus \eqref{eq:mathcalT} has a unique solution. Moreover,  
$\partial_b\langle\mathcal T f,\mathcal T g\rangle=0$. Since $\mathcal T(b_{\mathrm{ref}})=I,$
$\mathcal T(b)$ is orthogonal
in $L^2$, satisfies $\mathcal T(b)e_0(b_{\mathrm{ref}})=e_0(b)$, and is a
uniform isomorphism of $H^1_{\mathrm{ax}}$ for $b\in[0,\frac n2]$.
Write 
\[
\lambda(t)=\lambda_0(b(t)),\quad e(t)=e_0(b(t)),\quad \mathcal T(t)=\mathcal T(b(t)).
\]

With respect to the fixed decomposition
$\operatorname{span}\{e_0(b_{\mathrm{ref}})\}\oplus e_0(b_{\mathrm{ref}})^\perp$, write
\[
\mathcal T(t)(q,z):=\mathcal T(t)\bigl(qe_0(b_{\mathrm{ref}})+z\bigr),
\qquad q\in\mathbb R,\quad z\perp e_0(b_{\mathrm{ref}}).
\]
For a weak solution $D$ of \eqref{eq:homoeq}, define the scalar functions
\[
\begin{aligned}
q(t)&=\langle e(t),D(t)\rangle
=\int_0^{\pi/2}m(\theta)e(t,\theta)D(t,\theta)\,\mathrm d\theta,\\
p(t)&=\langle e(t),D_t(t)\rangle
=\int_0^{\pi/2}m(\theta)e(t,\theta)D_t(t,\theta)\,\mathrm d\theta.
\end{aligned}
\]
Since $\mathcal T(t)$ is orthogonal and
$\mathcal T(t)e_0(b_{\mathrm{ref}})=e(t)$, the decompositions
\[
\begin{aligned}
\mathcal T(t)^{-1}D(t)&=q(t)e_0(b_{\mathrm{ref}})+z(t),\\
\mathcal T(t)^{-1}D_t(t)&=p(t)e_0(b_{\mathrm{ref}})+y(t)
\end{aligned}
\]
uniquely define $z(t),y(t)\in e_0(b_{\mathrm{ref}})^\perp$.
Thus $D=\mathcal T(q,z)=qe+\mathcal Tz$ and
$D_t=\mathcal T(p,y)=pe+\mathcal Ty$, with
$\mathcal Tz,\mathcal Ty\perp e(t)$.

Set $c_\perp=\mathcal T^{-1}e'(t)$. Then $\|c_\perp\|_2=\|\mathcal T^{-1}e'\|_2=\|e'\|_2$. Since $\|e(t)\|_2=1$,
\[
\begin{aligned}
\langle c_\perp,e_0(b_{\mathrm{ref}})\rangle
&=\langle e',\mathcal T e_0(b_{\mathrm{ref}})\rangle
=\langle e',e\rangle
=\frac12\frac{\mathrm d}{\mathrm dt}\|e\|_2^2=0.
\end{aligned}
\]
Hence $c_\perp\in e_0(b_{\mathrm{ref}})^\perp$. In general $p\ne q'$, since
\[
q'=\langle e,D_t\rangle+\langle e',D\rangle=p + \langle \mathcal T^{-1}e',\mathcal T^{-1}D\rangle
=p+\langle c_\perp,qe_0(b_{\mathrm{ref}})+z\rangle
=p+\langle c_\perp,z\rangle.
\]
Differentiating the decompositions and using \eqref{eq:mathcalT} and
$D_{tt}=A_{b(t)}D$ gives, almost everywhere in $t$,
\begin{equation}\label{eq:frame-system}
\begin{aligned}
q'&=p+\langle c_\perp,z\rangle,&
z'&=y-c_\perp q,\\
p'&=\lambda q+\langle c_\perp,y\rangle,&
y'&=A_\perp(t)z-c_\perp p,
\end{aligned}
\end{equation}
where $A_\perp(t)$ is the  self-adjoint operator in  \(e_0(b_{\mathrm{ref}})^\perp\subset L^2\) associated with the form
\[
Q_\perp(t)[z,\zeta]
:=Q_{b(t)}[\mathcal T(t)(0,z),\mathcal T(t)(0,\zeta)],
\]
with form domain $H^1_\perp=H^1_{\mathrm{ax}}\cap e_0(b_{\mathrm{ref}})^\perp$. Write
$H^{-1}_\perp$ for the dual of $H^1_\perp$. The spectral gap gives
$Q_\perp(t)[z,z]\ge \frac{n^2}{4}\|z\|_2^2$. Together with
$Q_b[f,f]\ge\frac12\int_0^{\pi/2} m|f'|^2-C_n\|f\|_2^2$ and trace continuity, this
shows that $Q_\perp(t)[z,z]$ is uniformly equivalent to $\|z\|_{H^1_\perp}^2$.

Set 
\[
E_0^2=\frac{n^2}{4}q^2+p^2\quad\mbox{and}\quad E_\perp^2=\frac{n^2}{4}\|z\|_2^2+\|y\|_2^2.
\]
By \eqref{eq:frame-system}, 
\[
\begin{aligned}
\frac12(E_0^2)'
&=\frac{n^2}{4}qq'+pp'=\left(\frac{n^2}{4}+\lambda\right)qp
+\frac{n^2}{4}q\langle c_\perp,z\rangle
+p\langle c_\perp,y\rangle.
\end{aligned}
\]
Using $n|qp|\le E_0^2$, the Cauchy--Schwarz inequality, and
$\|c_\perp\|_2=\|e'\|_2$, we obtain
\[
\begin{aligned}
\left|\frac12(E_0^2)'\right|
&\le\frac{\left|\frac{n^2}{4}+\lambda\right|}{n}E_0^2
+\|c_\perp\|_2
\left(\frac{n^2}{4}|q|\|z\|_2+|p|\|y\|_2\right)\le\left(\frac n2-2\eta\right)E_0^2
+\|e'\|_2E_0E_\perp,
\end{aligned}
\]
where $\eta$ is the one in Lemma \ref{lem:quantitative-gap}.
On $\{E_0>0\}$, division by $E_0$ gives
\begin{equation}\label{eq:low-differential}
|E_0'|\le\left(\frac n2-2\eta\right)E_0+\|e'\|_2 E_\perp.
\end{equation}
 To include the zero set, set
$E_{0,\varepsilon}=(E_0^2+\varepsilon^2)^{1/2}$. The same estimate yields
\[
\begin{aligned}
|E_{0,\varepsilon}'|
&\le\left(\frac n2-2\eta\right)\frac{E_0^2}{E_{0,\varepsilon}}
+\|e'\|_2E_\perp\frac{E_0}{E_{0,\varepsilon}}\le\left(\frac n2-2\eta\right)E_0+\|e'\|_2E_\perp.
\end{aligned}
\]
Integrating on any finite interval and letting $\varepsilon\downarrow0$
proves \eqref{eq:low-differential} almost everywhere.

Fix $t_0\in\mathbb R$ and define the two-sided weight
\[
\omega_{t_0}(t)
:=\exp\left(-2\left|\int_{t_0}^t\Big(\frac n2-\eta(\tau)\Big)\,\mathrm d\tau\right|\right),
\]
where $\eta$ is the one in Lemma \ref{lem:quantitative-gap}.
Write $\omega=\omega_{t_0}$. Since $0<\eta\le n/4$,
\[
0<\omega(t)\le e^{-\frac n2|t-t_0|},
\qquad |(\sqrt\omega)'|\le\frac n2\sqrt\omega.
\]

\begin{lemma}\label{lem:weighted-energy}
Let $b$ satisfy \eqref{eq:coefficient-class}, and let $D$ be a weak
solution of \eqref{eq:homoeq}. If
$\int_{\mathbb R}\omega\|D\|_2^2\,\mathrm dt<\infty$, then
\begin{equation}\label{eq:weighted-form}
\int_{\mathbb R}\omega
\bigl(\|D_t\|_2^2+\|D\|_{H^1_{\mathrm{ax}}}^2\bigr)\,\mathrm dt
\le C_n\int_{\mathbb R}\omega\|D\|_2^2\,\mathrm dt,
\end{equation}
where $C_n$ depends only on $n$.
\end{lemma}

\begin{proof}
For a compactly supported Lipschitz function $\chi$, testing the weak
equation with $\chi^2D$, justified by density, gives
\[
\int_{\mathbb R}\|(\chi D)_t\|_2^2\,\mathrm dt
+\int_{\mathbb R}\chi^2Q_b[D,D]\,\mathrm dt
=\int_{\mathbb R}|\chi'|^2\|D\|_2^2\,\mathrm dt.
\]
Using
$Q_b[f,f]+C_n\|f\|_2^2\ge\frac12\|f\|_{H^1_{\mathrm{ax}}}^2$
and $\chi D_t=(\chi D)_t-\chi'D$, we obtain
\[
\int_{\mathbb R}\chi^2
\bigl(\|D_t\|_2^2+\|D\|_{H^1_{\mathrm{ax}}}^2\bigr)\,\mathrm dt
\le C_n\int_{\mathbb R}(|\chi'|^2+\chi^2)\|D\|_2^2\,\mathrm dt.
\]
Choose $0\le\chi_R\le1$ in $C_c^\infty(\mathbb R)$ with
$\chi_R=1$ on $[-R,R]$, $\operatorname{supp}\chi_R\subset[-2R,2R]$,
and $|\chi_R'|\le C/R$.
Taking $\chi=\sqrt\omega\,\chi_R$ and using
$|(\sqrt\omega)'|\le\frac n2\sqrt\omega$ gives, for $R\ge1$,
\[
\begin{aligned}
&\int_{\mathbb R}\omega\chi_R^2
\bigl(\|D_t\|_2^2+\|D\|_{H^1_{\mathrm{ax}}}^2\bigr)\,\mathrm dt\le C_n\int_{\mathbb R}
\omega(\chi_R^2+|\chi_R'|^2)\|D\|_2^2\,\mathrm dt
\le C_n\int_{\mathbb R}\omega\|D\|_2^2\,\mathrm dt.
\end{aligned}
\]
Fatou's lemma yields \eqref{eq:weighted-form} as $R\to\infty$.
Finally, if $D\in L^\infty(\mathbb R;L^2)$, then
\[
\int_{\mathbb R}\omega\|D\|_2^2\,\mathrm dt
\le\frac4n\|D\|_{L^\infty(\mathbb R;L^2)}^2<\infty.
\]
\end{proof}

\begin{lemma}\label{lem:principal-cauchy-control}
Let $b$ satisfy \eqref{eq:coefficient-class}, fix $t_0\in\mathbb R$,
and write $\omega=\omega_{t_0}$. If $D$ is a weak solution of
\eqref{eq:homoeq} with
$\int_{\mathbb R}\omega\|D\|_2^2\,\mathrm dt<\infty$, then
\begin{equation}\label{eq:observation}
(1-\varepsilon_n)\int_\mathbb{R}\omega \eta
 \bigl(\frac{n^2}{4}\|D\|_2^2+\|D_t\|_2^2\bigr)\,\mathrm{d} t
\le \frac{n^2}{4}q(t_0)^2+p(t_0)^2.
\end{equation}
where $\varepsilon_n\in(0,1)$ is the one in Lemma \ref{lem:quantitative-gap}. In particular, $q(t_0)=p(t_0)=0$ implies $D\equiv0$.
\end{lemma}

\begin{proof}

Lemma~\ref{lem:weighted-energy} gives
\[
\int_{\mathbb R}\omega(E_0^2+E_\perp^2)\,\mathrm dt
=\int_{\mathbb R}\omega
\left(\frac{n^2}{4}\|D\|_2^2+\|D_t\|_2^2\right)\,\mathrm dt<\infty.
\]
On $t>t_0$, differentiation of $\omega E_0^2$ using
\eqref{eq:low-differential} gives
\[
(\omega E_0^2)'\le-2\omega \eta E_0^2+2\omega \|e'\|_2 E_0 E_\perp.
\]
Integrating on both half-lines, with time reversed for $t<t_0$, gives
\begin{equation}\label{eq:low-integrated}
\int_\mathbb{R}\omega \eta E_0^2\,\mathrm{d} t
\le E_0(t_0)^2+\int_\mathbb{R}\omega \|e'\|_2 E_0 E_\perp\,\mathrm{d} t.
\end{equation}

Recall  the Hilbert-space Sobolev product rule 
\begin{equation}\label{eq:weak-product}
\frac{\mathrm d}{\mathrm dt}\langle z,y\rangle
=\langle z',y\rangle+\langle y',z\rangle_{H^{-1}_\perp,H^1_\perp},
\end{equation}
where $z\in H^1_{\rm loc}(L^2)\cap L^2_{\rm loc}(H^1_\perp)$,
$y\in L^2_{\rm loc}(L^2)$, and $y'\in L^2_{\rm loc}(H^{-1}_\perp)$.
Then it follows from \eqref{eq:local-H2} and \eqref{eq:weak-product} that
\[
\frac{\mathrm d}{\mathrm d t}\langle z,y\rangle
=\|y\|_2^2+Q_\perp[z,z]
-q\langle c_\perp,y\rangle-p\langle c_\perp,z\rangle,
\qquad
n|\langle z,y\rangle|\le E_\perp^2.
\]
For the coupling terms,
\[
\begin{aligned}
\left|q\langle c_\perp,y\rangle+p\langle c_\perp,z\rangle\right|
&\le\|e'\|_2\bigl(|q|\|y\|_2+|p|\|z\|_2\bigr)\\
&=\frac2n\|e'\|_2
\left(\frac n2|q|\|y\|_2+|p|\frac n2\|z\|_2\right)\\
&\le\frac2n\|e'\|_2E_0E_\perp.
\end{aligned}
\]
Since $Q_\perp[z,z]\ge\frac{n^2}{4}\|z\|_2^2$, it follows that
\[
E_\perp^2
\le\frac{\mathrm d}{\mathrm dt}\langle z,y\rangle
+\frac2n\|e'\|_2E_0E_\perp.
\]
Integrating by parts and using $n|\langle z,y\rangle|\le E_\perp^2$
and $|\omega'|\le2(\frac n2-\eta)\omega$, we obtain
\[
\begin{aligned}
\int_{\mathbb R}\omega E_\perp^2\,\mathrm dt
&\le-\int_{\mathbb R}\omega'\langle z,y\rangle\,\mathrm dt
+\frac2n\int_{\mathbb R}\omega\|e'\|_2E_0E_\perp\,\mathrm dt\\
&\le\frac1n\int_{\mathbb R}|\omega'|E_\perp^2\,\mathrm dt
+\frac2n\int_{\mathbb R}\omega\|e'\|_2E_0E_\perp\,\mathrm dt\\
&\le\int_{\mathbb R}\omega\left(1-\frac{2\eta}{n}\right)E_\perp^2\,\mathrm dt
+\frac2n\int_{\mathbb R}\omega\|e'\|_2E_0E_\perp\,\mathrm dt,
\end{aligned}
\]
where the integration by parts in the first inequality can be justified by inserting the cutoffs $\chi_R$
and passing to the limit using \eqref{eq:weighted-form}.
Rearranging gives
\begin{equation}\label{eq:high-integrated}
\int_\mathbb{R}\omega \eta E_\perp^2\,\mathrm{d} t
\le\int_\mathbb{R}\omega \|e'\|_2 E_0 E_\perp\,\mathrm{d} t.
\end{equation}
Adding \eqref{eq:low-integrated} and \eqref{eq:high-integrated}, and
using $\|e'\|_2\le\varepsilon_n\eta$ in \eqref{eq:strict-gap} and
$2E_0 E_\perp\le E_0^2+E_\perp^2$, proves \eqref{eq:observation}.
If $q(t_0)=p(t_0)=0$, then $\omega\eta>0$ and
$\varepsilon_n<1$ imply $D\equiv0$.
\end{proof}

\begin{definition}\label{def:weighted-solution-space}
Let $b$ satisfy \eqref{eq:coefficient-class} and set $\omega=\omega_0$.
We denote by $\mathcal C_b$ the space of weak solutions of
\eqref{eq:homoeq} satisfying
\[
\int_{\mathbb R}\omega
\bigl(\|D_t\|_2^2+\|D\|_{H^1_{\mathrm{ax}}}^2\bigr)\,\mathrm dt<\infty.
\]
\end{definition}

The space $\mathcal C_b$ is independent of the choice of $t_0$, since $e^{-n|t_0|}\omega_0(t)\le\omega_{t_0}(t)\le e^{n|t_0|}\omega_0(t)$.
By \eqref{eq:observation}, $q(0)=p(0)=0$ implies $D\equiv0$.
Thus the two principal Cauchy data determine every
$D\in\mathcal C_b$ uniquely.

\begin{proposition}\label{prop:linear-two-dimensionality}
For every coefficient satisfying \eqref{eq:coefficient-class}, the principal
Cauchy-data map
\[
D\longmapsto
\bigl(\langle e(0),D(0)\rangle,\langle e(0),D_t(0)\rangle\bigr)
\]
is a bijective linear map from $\mathcal C_b$ onto $\mathbb R^2$.
In particular, $\dim\mathcal C_b=2$.
\end{proposition}

\begin{proof}
It is clear that the principal Cauchy-data map is linear. Injectivity follows from \eqref{eq:observation}. For surjectivity, first solve
\begin{equation}\label{eq:high-forced}
z'=y+F_1,\qquad y'=A_\perp(t)z+F_2.
\end{equation}
Since $\eta\ge1/3$, our weight satisfies
$|\omega'|/(2\omega)\le s_*:=n/2-1/3<n/2$.
The following variational construction applies to any positive
locally Lipschitz weight satisfying this bound.
Complete the compactly supported smooth $H^1_\perp$-valued functions under
\[
\|z\|_{\omega}^2
=\int_\mathbb{R}\omega\bigl(\|z'\|_2^2+Q_\perp(t)[z,z]\bigr)\,\mathrm{d} t.
\]
Since $Q_\perp(t)[z,z]\asymp\|z\|_{H^1_\perp}^2$, this is the weighted space
with $z'\in L^2(\omega\,\mathrm dt;L^2)$ and
$z\in L^2(\omega\,\mathrm dt;H^1_\perp)$. Density follows by truncation and
time convolution, using local positive lower and upper bounds for $\omega$.

Eliminating $y$ from \eqref{eq:high-forced} gives
$-z''+A_\perp z=-F_1'-F_2$. Its weighted variational formulation is
$B(z,\phi)=\ell(\phi)$ for every $\phi$ in the weighted energy space, where
\begin{align}
B(z,\phi)
&=\int_\mathbb{R}\omega\left(
 \langle z',\phi'\rangle+\frac{\omega'}{\omega}\langle z',\phi\rangle+Q_\perp(t)[z,\phi]\right)\,\mathrm{d} t,\label{eq:bilinearform}\\
 \ell(\phi)&=\int_\mathbb{R}\omega\left(
 \langle F_1,\phi'\rangle+\frac{\omega'}{\omega}\langle F_1,\phi\rangle-\langle F_2,\phi\rangle\right)\,\mathrm{d} t. \nonumber
\end{align}
The form $B$ is bounded, and
\[
\left|\frac{\omega'}{\omega}\langle z',z\rangle\right|
\le\frac{2s_*}{n}\bigl(\|z'\|_2^2+\frac{n^2}{4}\|z\|_2^2\bigr),
\qquad
B(z,z)\ge(1-2s_*/n)\|z\|_{\omega}^2.
\]
Thus Lax--Milgram gives a unique $z$ for weighted $L^2$
data $F_1,F_2$. With $y=z'-F_1$, this is the unique weak solution of
\eqref{eq:high-forced} in the weighted energy class.

Set $F_\perp^2=\frac{n^2}{4}\|F_1\|_2^2+\|F_2\|_2^2$. Recall $E_0^2=\frac{n^2}{4}q^2+p^2$ and  $E_\perp^2=\frac{n^2}{4}\|z\|_2^2+\|y\|_2^2$. Then
\[
\frac{\mathrm d}{\mathrm dt}\langle z,y\rangle
=\|y\|_2^2+Q_\perp[z,z]+\langle F_1,y\rangle+\langle z,F_2\rangle,\qquad
|\langle F_1,y\rangle+\langle z,F_2\rangle|\le \frac2n F_\perp E_\perp.
\]
Since \(Q_\perp[z,z]\ge\frac{n^2}{4}\|z\|_2^2\), we have \[E_\perp^2\le(\langle z,y\rangle)'+\frac2nF_\perp E_\perp.\]
Multiplication by $\omega$, using
$n|\langle z,y\rangle|\le E_\perp^2$ and $|\omega'|\le2(\frac n2-\eta)\omega$, followed by the cutoff argument above, gives, as in \eqref{eq:high-integrated}
\[
\int_\mathbb{R}\omega \eta E_\perp^2\,\mathrm{d} t
\le\int_\mathbb{R}\omega F_\perp E_\perp\,\mathrm{d} t.
\]
Consequently,
\begin{equation}\label{eq:EFinverse}
\|E_\perp\|_{L^2(\omega \eta\,\mathrm{d} t)}
\le\|F_\perp\|_{L^2(\omega/\eta\,\mathrm{d} t)}.
\end{equation}

For the system governing the principal components $q$ and $p$, set
$\mathbf u=(\frac n2q,p)$ and
\[
M_0(t)=\begin{pmatrix}0&\frac n2\\ \frac{2\lambda(t)}{n}&0\end{pmatrix}.
\]
Its evolution operator $S(t,\tau)$, solving $\partial_tS(t,\tau)=M_0(t)S(t,\tau)$ and $S(\tau,\tau)=I$, satisfies
\[
\|S(t,\tau)\|
\le\exp\left(\left|\int_\tau^t\left(\frac n2-2\eta(s)\right)\,\mathrm{d}s\right|\right).
\]

Let $G\in L^2(\omega/\eta\,\mathrm dt;\mathbb R^2)$ and let
$\mathbf u$ solve $\mathbf u'=M_0\mathbf u+G$ with
$\mathbf u(0)=\mathbf u_0$.
Since $|\mathbf u\cdot M_0\mathbf u|
\le(\frac n2-2\eta)|\mathbf u|^2$, for almost every $t>0$,
\[
\begin{aligned}
(\omega|\mathbf u|^2)'
&=-2\Bigl(\frac n2-\eta\Bigr)\omega|\mathbf u|^2
  +2\omega\mathbf u\cdot(M_0\mathbf u+G)\le-2\omega\eta|\mathbf u|^2
  +2\omega|G||\mathbf u|.
\end{aligned}
\]
Integrating on $[0,R]$ and reversing time on $[-R,0]$ gives
\[
\begin{aligned}
\int_{-R}^R\omega\eta|\mathbf u|^2\,\mathrm dt
&+\frac12\omega(R)|\mathbf u(R)|^2
 +\frac12\omega(-R)|\mathbf u(-R)|^2\le|\mathbf u_0|^2
 +\int_{-R}^R\omega|G||\mathbf u|\,\mathrm dt.
\end{aligned}
\]
For $\mathbf u_0=0$, the Cauchy--Schwarz inequality yields
\[
\begin{aligned}
\int_{-R}^R\omega\eta|\mathbf u|^2\,\mathrm dt
&\le
\left(\int_{-R}^R\frac{\omega}{\eta}|G|^2\,\mathrm dt\right)^{1/2}
\left(\int_{-R}^R\omega\eta|\mathbf u|^2\,\mathrm dt\right)^{1/2}.
\end{aligned}
\]
Letting $R\to\infty$, we obtain
\begin{equation}\label{eq:low-inverse}
\|\mathbf u\|_{L^2(\omega\eta\,\mathrm dt)}
\le\|G\|_{L^2(\omega/\eta\,\mathrm dt)}.
\end{equation}
For $G=0$, the same integrated inequality gives
\[
\int_{\mathbb R}\omega\eta|\mathbf u|^2\,\mathrm dt
\le|\mathbf u_0|^2.
\]

In \eqref{eq:frame-system},
\[
\begin{aligned}
G&=\left(\frac n2\langle c_\perp,z\rangle,\langle c_\perp,y\rangle\right),
& |G|&\le\|e'\|_2E_\perp,\\
F_1&=-c_\perp q,\qquad F_2=-c_\perp p,
& F_\perp&\le\|e'\|_2E_0.
\end{aligned}
\]
Consider the Hilbert space of triples of functions $(\mathbf u,z,y)$, where
$\mathbf u=(\frac n2q,p)$ and $z(t),y(t)\perp e_0(b_{\mathrm{ref}})$, with finite norm
\[
\|(\mathbf u,z,y)\|^2
=\int_{\mathbb R}\omega\eta
\left(|\mathbf u|^2+\frac{n^2}{4}\|z\|_2^2+\|y\|_2^2\right)\,\mathrm dt
=\int_{\mathbb R}\omega\eta(E_0^2+E_\perp^2)\,\mathrm dt.
\]
By \eqref{eq:strict-gap},
\[
\begin{aligned}
\int_\mathbb{R}\frac{\omega}{\eta}\|e'\|_2^2E_\perp^2\,\mathrm{d} t
&\le\varepsilon_n^2\int_\mathbb{R}\omega \eta E_\perp^2\,\mathrm{d} t,\\
\int_\mathbb{R}\frac{\omega}{\eta}\|e'\|_2^2E_0^2\,\mathrm{d} t
&\le\varepsilon_n^2\int_\mathbb{R}\omega \eta E_0^2\,\mathrm{d} t.
\end{aligned}
\]

Fix $\mathbf u_0\in\mathbb R^2$. For each triple, take $G,F_1,F_2$ as above and
define $\mathcal K(\mathbf u,z,y)=(\widetilde{\mathbf u},\widetilde z,\widetilde y)$ by
\[
\begin{aligned}
\widetilde{\mathbf u}(t)
&=S(t,0)\mathbf u_0+\int_0^t S(t,\tau)G(\tau)\,\mathrm d\tau,\\
\widetilde z'&=\widetilde y+F_1,\qquad
\widetilde y'=A_\perp(t)\widetilde z+F_2,
\end{aligned}
\]
where $(\widetilde z,\widetilde y)$ is the weighted variational solution.
The preceding estimates show that $\mathcal K$ maps this space into itself.
For two inputs and their images, write $\delta$ for the corresponding
differences. Since $\delta\widetilde{\mathbf u}(0)=0$, \eqref{eq:EFinverse}, \eqref{eq:low-inverse} and \eqref{eq:strict-gap} give
\[
\begin{aligned}
\int_{\mathbb R}\omega\eta|\delta\widetilde{\mathbf u}|^2\,\mathrm dt
&\le\int_{\mathbb R}\frac{\omega}{\eta}|\delta G|^2\,\mathrm dt\le\varepsilon_n^2\int_{\mathbb R}\omega\eta
\left(\frac{n^2}{4}\|\delta z\|_2^2+\|\delta y\|_2^2\right)\,\mathrm dt,\\
\int_{\mathbb R}\omega\eta
\left(\frac{n^2}{4}\|\delta\widetilde z\|_2^2
+\|\delta\widetilde y\|_2^2\right)\,\mathrm dt
&\le\int_{\mathbb R}\frac{\omega}{\eta}
\left(\frac{n^2}{4}\|\delta F_1\|_2^2+\|\delta F_2\|_2^2\right)\,\mathrm dt\le\varepsilon_n^2\int_{\mathbb R}\omega\eta|\delta\mathbf u|^2\,\mathrm dt.
\end{aligned}
\]
Adding yields
\[
\|(\delta\widetilde{\mathbf u},\delta\widetilde z,\delta\widetilde y)\|
\le\varepsilon_n\|(\delta\mathbf u,\delta z,\delta y)\|,
\qquad \varepsilon_n<1.
\]
By Banach's fixed-point theorem, there is a unique triple satisfying
\[
(\mathbf u,z,y)=\mathcal K(\mathbf u,z,y),\qquad
\mathbf u(0)=\mathbf u_0.
\]
At this fixed point, \eqref{eq:frame-system} holds, and the variational
construction gives
\[
z\in L^2(\omega\,\mathrm dt;H^1_\perp),\qquad
z'\in L^2(\omega\,\mathrm dt;L^2).
\]

Set $D=\mathcal T(q,z)$. By \eqref{eq:mathcalT} and \eqref{eq:frame-system},
\[
\begin{aligned}
D_t
&=\mathcal T\bigl(q'-\langle c_\perp,z\rangle,\,z'+c_\perp q\bigr)
=\mathcal T(p,y),\\
D_{tt}
&=\mathcal T\bigl(p'-\langle c_\perp,y\rangle,\,y'+c_\perp p\bigr)=\mathcal T(\lambda q,A_\perp z)=A_{b(t)}D,
\end{aligned}
\]
where the second identity holds in $H^{-1}_{\mathrm{ax}}$.
Since $\eta\ge1/3$, the fixed-point norm and the preceding regularity,
together with the uniform $H^1_{\mathrm{ax}}$ bound for $\mathcal T$, give
\[
\begin{aligned}
&\int_{\mathbb R}\omega
\left(\|D\|_{H^1_{\mathrm{ax}}}^2+\|D_t\|_2^2\right)\,\mathrm dt\le C_n\int_{\mathbb R}\omega
\left(q^2+p^2+\|z\|_{H^1_\perp}^2+\|y\|_2^2\right)\,\mathrm dt<\infty.
\end{aligned}
\]
Thus $D\in\mathcal C_b$. The Cauchy traces are well defined by
\eqref{eq:local-H2}, and the orthogonality of $\mathcal T$ gives
\[
\langle e(0),D(0)\rangle=q(0),\qquad
\langle e(0),D_t(0)\rangle=p(0).
\]
Since $(\frac n2q(0),p(0))=\mathbf u_0$ and $\mathbf u_0\in\mathbb R^2$
is arbitrary, the principal Cauchy-data map is surjective.
\end{proof}

\begin{proposition}\label{prop:symplectic}
The bilinear form
\[
 \Omega(D,E)=\langle D,E_t\rangle-\langle D_t,E\rangle,\quad D,E\in\mathcal C_b,
\]
is independent of time and nondegenerate on $\mathcal C_b$.
\end{proposition}

The conservation of $\Omega$ is analogous to the Wronskian identity studied by Ao--Chan--DelaTorre--Fontelos--Gonz\'alez--Wei  \cite{AoEtAl2020}.

\begin{proof}
Let $D,E\in\mathcal C_b$. Applying \eqref{eq:weak-product} on the full
axial spaces with $(z,y)=(D,E_t)$ and $(z,y)=(E,D_t)$ gives
\[
\begin{aligned}
\frac{\mathrm d}{\mathrm dt}\langle D,E_t\rangle
&=\langle D_t,E_t\rangle
+\langle E_{tt},D\rangle_{H^{-1}_{\mathrm{ax}},H^1_{\mathrm{ax}}},\\
\frac{\mathrm d}{\mathrm dt}\langle E,D_t\rangle
&=\langle E_t,D_t\rangle
+\langle D_{tt},E\rangle_{H^{-1}_{\mathrm{ax}},H^1_{\mathrm{ax}}}.
\end{aligned}
\]
Since $\langle D_t,E\rangle=\langle E,D_t\rangle$,
$D_{tt}=A_{b(t)}D$, and $E_{tt}=A_{b(t)}E$, subtraction yields
\[
\begin{aligned}
\frac{\mathrm d}{\mathrm dt}\Omega(D,E)
&=\langle E_{tt},D\rangle_{H^{-1}_{\mathrm{ax}},H^1_{\mathrm{ax}}}
-\langle D_{tt},E\rangle_{H^{-1}_{\mathrm{ax}},H^1_{\mathrm{ax}}}=Q_{b(t)}[E,D]-Q_{b(t)}[D,E]=0.
\end{aligned}
\]
So $\Omega(D,E)$ is independent of $t$.

To prove nondegeneracy, define
\[
\omega_+(t)=\exp\left(-2\int_0^t\left(\frac n2-\eta(s)\right)\,\mathrm ds\right),\qquad
\omega_-(t)=\exp\left(2\int_0^t\left(\frac n2-\eta(s)\right)\,\mathrm ds\right).
\]
Since $\omega_\pm'=\mp2(\frac n2-\eta)\omega_\pm$, taking
$s_*=\frac n2-\frac13$ gives
\[
\begin{aligned}
\frac{|\omega_\pm'|}{2\omega_\pm}
&=\frac n2-\eta\le s_*<\frac n2,\qquad
B(z,z)\ge\left(1-\frac{2s_*}{n}\right)\|z\|_{\omega_\pm}^2
=\frac{2}{3n}\|z\|_{\omega_\pm}^2.
\end{aligned}
\]
Here $B$ is the variational form \eqref{eq:bilinearform} with weight $\omega_\pm$.
Thus Lax--Milgram gives a unique weighted variational solution of
\eqref{eq:high-forced}. The proof of \eqref{eq:EFinverse} uses only
$|\omega'|\le2(\frac n2-\eta)\omega$, and hence gives
\[
\begin{aligned}
&\int_{\mathbb R}\omega_\pm\eta
\left(\frac{n^2}{4}\|z\|_2^2+\|y\|_2^2\right)\,\mathrm dt\le\int_{\mathbb R}\frac{\omega_\pm}{\eta}
\left(\frac{n^2}{4}\|F_1\|_2^2+\|F_2\|_2^2\right)\,\mathrm dt.
\end{aligned}
\]

For the equation $\mathbf u'=M_0\mathbf u+G$, let
$G\in L^2(\omega_+/\eta\,\mathrm dt;\mathbb R^2)$ and set
\[
\mathbf u(t)=\int_{-\infty}^t S(t,\tau)G(\tau)\,\mathrm d\tau.
\]
For $t\ge\tau$, the evolution estimate gives
\[
\begin{aligned}
&\left\|\omega_+(t)^{\frac12}S(t,\tau)\omega_+(\tau)^{-\frac12}\right\|\le\exp\left(-\int_\tau^t\left(\frac n2-\eta\right)\,\mathrm ds
+\int_\tau^t\left(\frac n2-2\eta\right)\,\mathrm ds\right)=e^{-\int_\tau^t\eta(s)\,\mathrm ds}.
\end{aligned}
\]
Since $\eta\ge1/3$,
\[
\int_{-\infty}^t\eta(\tau)e^{-\int_\tau^t\eta(s)\,\mathrm ds}\,\mathrm d\tau=1,
\qquad
\int_\tau^\infty\eta(t)e^{-\int_\tau^t\eta(s)\,\mathrm ds}\,\mathrm dt=1.
\]
The Cauchy--Schwarz inequality therefore gives
\[
\begin{aligned}
\omega_+(t)|\mathbf u(t)|^2
&\le\left(\int_{-\infty}^te^{-\int_\tau^t\eta(s)\,\mathrm ds}\sqrt{\omega_+(\tau)}\,|G(\tau)|\,\mathrm d\tau\right)^2\\
&\le\left(\int_{-\infty}^t
e^{-\int_\tau^t\eta(s)\,\mathrm ds}\eta(\tau)\,\mathrm d\tau\right)\left(\int_{-\infty}^t e^{-\int_\tau^t\eta(s)\,\mathrm ds}
\frac{\omega_+(\tau)}{\eta(\tau)}|G(\tau)|^2\,\mathrm d\tau\right)\\
&=\int_{-\infty}^t e^{-\int_\tau^t\eta(s)\,\mathrm ds}
\frac{\omega_+(\tau)}{\eta(\tau)}|G(\tau)|^2\,\mathrm d\tau.
\end{aligned}
\]
The same bound applied to $\int_{-\infty}^t\|S(t,\tau)G(\tau)\|\,\mathrm d\tau$
shows that \(\mathbf u(t)\) is well defined.  Multiplication by $\eta(t)$ and Tonelli's theorem yield
\[
\begin{aligned}
\int_{\mathbb R}\omega_+\eta|\mathbf u|^2\,\mathrm dt
&\le\int_{\mathbb R}\frac{\omega_+(\tau)}{\eta(\tau)}|G(\tau)|^2
\left(\int_\tau^\infty\eta(t)e^{-\int_\tau^t\eta(s)\,\mathrm ds}\,\mathrm dt\right)\,\mathrm d\tau=\int_{\mathbb R}\frac{\omega_+}{\eta}|G|^2\,\mathrm dt.
\end{aligned}
\]
Moreover, $\partial_tS(t,\tau)=M_0(t)S(t,\tau)$ and $S(t,t)=I$ give
\[
\begin{aligned}
\mathbf u'(t)
&=G(t)+\int_{-\infty}^t M_0(t)S(t,\tau)G(\tau)\,\mathrm d\tau
=M_0(t)\mathbf u(t)+G(t),\\
\mathbf u(0)&=\int_{-\infty}^0S(0,\tau)G(\tau)\,\mathrm d\tau.
\end{aligned}
\]
Thus \eqref{eq:low-inverse} holds with $\omega_+$, with the initial value
determined by $G$. 

If $\mathbf u'=M_0\mathbf u$ and
$\mathbf u\in L^2(\omega_+\eta\,\mathrm dt;\mathbb R^2)$, choose
$s_j\to-\infty$ with $\omega_+(s_j)|\mathbf u(s_j)|^2\to0$. For fixed $t$
and large $j$,
\[
\begin{aligned}
\omega_+(t)^{1/2}|\mathbf u(t)|&\le \left\|\omega_+(t)^{1/2}S(t,s_j)\omega_+(s_j)^{-1/2}\right\| \omega_+(s_j)^{1/2}|\mathbf u(s_j)|\\
&\le e^{-(t-s_j)/3}\omega_+(s_j)^{1/2}|\mathbf u(s_j)|\\
&\to 0.
\end{aligned}
\]
Hence the weighted solution is unique.

For $\omega_-$, set
\[
\begin{aligned}
\mathbf u(t)&=-\int_t^\infty S(t,\tau)G(\tau)\,\mathrm d\tau,\\
\mathbf u'(t)&=G(t)-\int_t^\infty M_0(t)S(t,\tau)G(\tau)\,\mathrm d\tau
=M_0(t)\mathbf u(t)+G(t).
\end{aligned}
\]
For $t\le\tau$,
\[
\left\|\omega_-(t)^{1/2}S(t,\tau)\omega_-(\tau)^{-1/2}\right\|
\le e^{-\int_t^\tau\eta(s)\,\mathrm ds}\le e^{-(\tau-t)/3}.
\]
Time reversal gives the same estimate and uniqueness. In particular,
\[
\|\mathbf u\|_{L^2(\omega_\pm\eta\,\mathrm dt)}
\le\|G\|_{L^2(\omega_\pm/\eta\,\mathrm dt)}.
\]

Let $F\in L^2(\mathbb R;L^2)$ have compact support. The arguments in the proof of Proposition \ref{prop:linear-two-dimensionality} show that
\begin{equation}\label{eq:one-sided}
(-\partial_t^2+A_{b(t)})Y_\pm=F
\end{equation}
has a unique solution in each weighted space $\{Y\in L^2(\omega_\pm\,\mathrm dt;H^1_{\mathrm{ax}}),\ \ Y_t\in L^2(\omega_\pm\,\mathrm dt;L^2)\}$.

Fix $D\in\mathcal C_b\setminus\{0\}$. If $Y_\pm$ solves
\eqref{eq:one-sided}, then
\[
\frac{\mathrm{d}}{\mathrm{d} t}\Omega(D,Y_\pm)=-\langle D,F\rangle.
\]
Choose $0\le\chi\in C_c^\infty(\mathbb R)$
with $\int_\mathbb R\chi(t)\|D(t)\|_2^2\,\mathrm dt>0$. Set $F=\chi D$.
Let $Y_\pm$ be the two solutions just constructed. Since
$\omega\le\omega_\pm$, both belong to the desired weighted space, and
$E=Y_--Y_+\in\mathcal C_b$.

Choose $R>0$ with $\operatorname{supp}F\subset[-R,R]$. Since $t_0=0$,
\[
 \omega(t)\omega_+(t)=1\quad\mbox{for }t<0,\qquad
 \omega(t)\omega_-(t)=1\quad\mbox{for }t>0.
\]
By Cauchy--Schwarz,
\[
\begin{aligned}
|\Omega(D,Y_+)|
&\le\|D\|_2\|(Y_+)_t\|_2+\|D_t\|_2\|Y_+\|_2\le\left(\|D\|_2^2+\|D_t\|_2^2\right)^{\frac12}
\left(\|Y_+\|_2^2+\|(Y_+)_t\|_2^2\right)^{\frac12},
\end{aligned}
\]
and hence, by Cauchy--Schwarz again,
\[
\begin{aligned}
&\int_{-\infty}^{-R}|\Omega(D,Y_+)|\,\mathrm dt<\infty.
\end{aligned}
\]
Similarly, $\int_R^\infty|\Omega(D,Y_-)|\,\mathrm dt<\infty$.
Since $F=0$ on $(-\infty,-R)\cup(R,+\infty)$, $\Omega(D,Y_+)$ is constant on $(-\infty,-R)$ and $\Omega(D,Y_-)$ is constant on $(R,+\infty)$. Their integrability implies
\[
\Omega(D,Y_+)=0\quad\mbox{for }t<-R,\qquad
\Omega(D,Y_-)=0\quad\mbox{for }t>R.
\]
Integrating $\frac{\mathrm{d}}{\mathrm{d} t}\Omega(D,Y_\pm)=-\langle D,F\rangle$ gives
\[
\begin{aligned}
\Omega(D,Y_+)(t)
&=-\int_{-\infty}^t\langle D(s),F(s)\rangle\,\mathrm ds,\\
\Omega(D,Y_-)(t)
&=\int_t^\infty\langle D(s),F(s)\rangle\,\mathrm ds.
\end{aligned}
\]
Since $E=Y_--Y_+$ and $F=\chi D$,
\[
\begin{aligned}
\Omega(D,E)
&=\Omega(D,Y_-)(t)-\Omega(D,Y_+)(t)\\
&=\int_t^\infty\langle D(s),F(s)\rangle\,\mathrm ds
+\int_{-\infty}^t\langle D(s),F(s)\rangle\,\mathrm ds\\
&=\int_{\mathbb R}\langle D,F\rangle\,\mathrm dt\\
&=\int_{\mathbb R}\chi(t)\|D(t)\|_2^2\,\mathrm dt\\
&>0.
\end{aligned}
\]
Thus $\Omega$ is nondegenerate on $\mathcal C_b$.
\end{proof}

\subsection{Solutions of the linearized equation}\label{sec:full-linearized}

We next show that, at a positive periodic solution of \eqref{eq:cylinder},
every solution of the linearized equation with sufficiently slow
exponential growth is axially symmetric and belongs to $\mathcal C_b$.

\begin{lemma}\label{lem:linearized-axial-symmetry}
Let $T>0$ and let $v$ be a positive smooth $T$-periodic solution of
\eqref{eq:cylinder}, including the stationary case. Set
$w(t)=v(t,\pi/2)$ and $b(t)=\frac n2w(t)^{\frac2{n-2}}$.
Suppose $h\in C^\infty(\mathbb R\times\overline{\mathbb S^{n-1}_+};\mathbb R)$
satisfies
\[
\begin{cases}
h_{tt}+\Delta_{\mathbb S^{n-1}_+}h-\frac{(n-2)^2}{4}h=0
&\text{in }\mathbb R\times\mathbb S^{n-1}_+,\\
\partial_\nu h=b(t)h
&\text{on }\mathbb R\times\partial\mathbb S^{n-1}_+,
\end{cases}
\]
and \eqref{eq:linearized-growth} for some $C>0$ and
$0<\varepsilon<\min\{1,n/4\}$. Then $h$ is axially symmetric and
$h\in\mathcal C_b$.
\end{lemma}

\begin{proof}
Write $\omega=(\sin\theta\,\zeta,\cos\theta)$, with
$\zeta\in\mathbb S^{n-2}$. For a real spherical harmonic $Y$ satisfying
$-\Delta_{\mathbb S^{n-2}}Y=\ell(\ell+n-3)Y$ and
$\int_{\mathbb S^{n-2}}Y^2\,\mathrm d\zeta=1$, $\ell\ge1$, set
$$u(t,\theta)=\int_{\mathbb S^{n-2}}h(t,\theta,\zeta)Y(\zeta)\,\mathrm d\zeta.$$
Then
\[
\begin{cases}
\mathcal L_\ell u:=-u_{tt}-m^{-1}(mu_\theta)_\theta
+\left(\frac{(n-2)^2}{4}
+\frac{\ell(\ell+n-3)}{\sin^2\theta}\right)u=0,\\
u_\theta(t,\pi/2)=b(t)u(t,\pi/2),\\
u(t,\theta)=O(\sin^\ell\theta)\quad\mbox{as }\theta\to0.
\end{cases}
\]
The last assertion follows from smoothness at the pole.

Define
\[
q_\pm=\frac{n-2}{2}v\pm v_t-\cot\theta\,v_\theta.
\]
Equation~\eqref{eq:cylinder} gives
\[
\begin{aligned}
-(q_\pm)_{tt}\mp2(q_\pm)_t
-\sin^{-n}\theta\bigl(\sin^n\theta(q_\pm)_\theta\bigr)_\theta
+\left(\frac{n^2}{4}-1\right)q_\pm&=0,\\
(q_\pm)_\theta(t,\pi/2)&=b(t)q_\pm(t,\pi/2).
\end{aligned}
\]
The functions $q_\pm$ extend smoothly and evenly across $\theta=0$.
Moreover, \eqref{eq:sharp} gives
\[
\begin{aligned}
q_\pm(t,\pi/2)
&=\frac{n-2}{2}w(t)\pm w'(t)\ge\frac{n-2}{2}w(t)
\left(1-\sqrt{1-w(t)^{\frac4{n-2}}}\right)>0.
\end{aligned}
\]
By viewing $\sin^{-n}\theta(\sin^n\theta\,\partial_\theta)_\theta$
as the axial Laplacian on $\mathbb S^{n+1}_+$, the maximum principle on
$(\mathbb R/T\mathbb Z)\times\mathbb S^{n+1}_+$ and
$\frac{n^2}{4}-1>0$ imply
$0<C_v^{-1}\le q_\pm\le C_v$.
Consequently, for
\[
Z_\pm(t,\theta):=e^{\pm t}\sin\theta\,q_\pm(t,\theta)>0
\quad\mbox{with }0<\theta\le\pi/2,
\]
these functions satisfy
\begin{equation}\label{eq:positive-angular-modes}
\begin{cases}
\mathcal L_1Z_\pm=0,\\
(Z_\pm)_\theta(t,\pi/2)=b(t)Z_\pm(t,\pi/2).
\end{cases}
\end{equation}

Local boundary estimates, periodicity of $b$, and
\eqref{eq:linearized-growth} yield
$\|h(t,\cdot)\|_{C^1(\overline{\mathbb S^{n-1}_+})}\le Ce^{\varepsilon|t|}$.
Since $u(t,0)=0$, it follows that
$|u(t,\theta)|\le C_Ye^{\varepsilon|t|}\sin\theta$.
For $R>1$, choose
\[
\Psi_R=C_Ye^{-(1-\varepsilon)R}(Z_++Z_-),
\]
with $C_Y$ large enough, independent of $R$. Then $|u(\pm R,\theta)|\le \Psi_R(\pm R,\theta)$, 
$(\Psi_R)_\theta(t,\pi/2)=b(t)\Psi_R(t,\pi/2)$, and
\[
\mathcal L_\ell \Psi_R
=\frac{\ell(\ell+n-3)-(n-2)}{\sin^2\theta}\Psi_R\ge0.
\]
Set $f=(\pm u-\Psi_R)_+$ on $(-R,R)\times(0,\pi/2)$, so that
$f(\pm R,\cdot)=0$. Testing with $f$ and using
\eqref{eq:positive-angular-modes}, integration by parts gives
\[
\begin{aligned}
0&\ge\int_{-R}^R\int_0^{\pi/2}
 m\,\mathcal L_\ell(\pm u-\Psi_R)f\,\mathrm d\theta\,\mathrm dt\\
&=\int_{-R}^R\int_0^{\pi/2}mZ_+^2
 \left(\left|\partial_t\frac f{Z_+}\right|^2
 +\left|\partial_\theta\frac f{Z_+}\right|^2\right)
 \,\mathrm d\theta\,\mathrm dt\\
&\quad+\bigl(\ell(\ell+n-3)-(n-2)\bigr)
 \int_{-R}^R\int_0^{\pi/2}\frac m{\sin^2\theta}f^2
 \,\mathrm d\theta\,\mathrm dt\ge0.
\end{aligned}
\]
Thus $f/Z_+$ is constant. Since $f$ is zero at $t=\pm R$, we have $f\equiv 0$.
Hence $|u|\le \Psi_R$. Letting $R\to\infty$ and using $\varepsilon<1$
gives $u=0$. Since $Y$ was arbitrary with $\ell\ge1$, completeness of
the spherical harmonics implies that $h$ is axially symmetric.

By Lemma~\ref{lem:secant}, $b$ satisfies \eqref{eq:coefficient-class}.
For the weight $\omega_0$ defining $\mathcal C_b$, we have
$\omega_0(t)\le e^{-\frac n2|t|}$. Since $\varepsilon<n/4$,
\[
\int_{\mathbb R}\omega_0(t)\|h(t)\|_2^2\,\mathrm dt
\le C\int_{\mathbb R}e^{-(\frac n2-2\varepsilon)|t|}\,\mathrm dt<\infty.
\]
Estimate~\eqref{eq:weighted-form} gives $h\in\mathcal C_b$.
\end{proof}

\section{Hamiltonian structure and the global Fowler family}\label{sec:hamiltonian-global}

\subsection{Local coordinates and Hamiltonian level sets}
Let $\mathcal S:=\{0\}\cup\{v:\  v\text{ is a positive bounded entire axially symmetric solution of \eqref{eq:cylinder}}\}$.
Endow $\mathcal S$ with the topology induced by
\[
 v\longmapsto (v,v_t)\in
 C_{\mathrm{loc}}\bigl(\mathbb{R};H^1_{\mathrm{ax}}\times L^2(m\,\mathrm{d}\theta)\bigr).
\]
More explicitly, \(v_j\to v\) in \(\mathcal S\) means that, for every \(R>0\),
\[
\sup_{|t|\le R}
\left(
\|v_j(t)-v(t)\|_{H^1_{\mathrm{ax}}}
+\|(v_j)_t(t)-v_t(t)\|_2
\right)\rightarrow0.
\]

We describe nearby solutions by two fixed scalar projections of their Cauchy data at \(t=0\).

\begin{proposition}\label{prop:coordinates}
For every $v_0\in\mathcal S$ there are a neighborhood $\mathcal U$ of $v_0$ in $\mathcal S$, a number $\bar b\in[0,\frac n2]$, and the normalized eigenfunction $e_*=e_0(\bar b)$ such that
\[
 \mathbf c(v)=\bigl(\frac n2\langle e_*,v(0)\rangle,\langle e_*,v_t(0)\rangle\bigr),
 \qquad v\in\mathcal U,
\]
is injective.  For all $v_1,v_2\in\mathcal U$,
\begin{equation}
 C^{-1}|\mathbf c(v_1)-\mathbf c(v_2)|
 \le
 \|(v_1-v_2,(v_1-v_2)_t)(0)\|_{H^1_{\mathrm{ax}}\times L^2}
 \le C|\mathbf c(v_1)-\mathbf c(v_2)|.
\end{equation}
Moreover, for every $R>0$ there is $C_R>0$ such that
\begin{equation}
 \|v_1-v_2\|_{H^2([-R,R]\times\mathbb S^{n-1}_+)}
 +\|v_1(\,\cdot,\pi/2)-v_2(\,\cdot,\pi/2)\|_{C^0([-R,R])}
 \le C_R|\mathbf c(v_1)-\mathbf c(v_2)|.
\end{equation}
Thus, with respect to the Cauchy-data distance at $t=0$,
$\mathbf c$ is a bi-Lipschitz embedding of $\mathcal U$ into
$\mathbb R^2$, and it is a homeomorphism onto its image for
the topology of $\mathcal S$.
\end{proposition}
\begin{proof}
Let $v_1,v_2\in\mathcal S$, set $D=v_1-v_2$, and let $b$
be the corresponding coefficient defined in \eqref{eq:secant}.
Since $\eta\ge1/3$ and $\omega_{t_0}\ge e^{-nR}$ on
$[t_0-R,t_0+R]$, \eqref{eq:observation} gives
\begin{equation}\label{eq:local-L2}
\int_{t_0-R}^{t_0+R}
 \bigl(\frac{n^2}{4}\|D\|_2^2+\|D_t\|_2^2\bigr)\,\mathrm{d} t
\le\frac{3e^{nR}}{1-\varepsilon_n}E_0(t_0)^2.
\end{equation}

By \eqref{eq:local-H2} and the trace theorem, $D(t_0)\in H^{3/2}$ and
$D_t(t_0)\in H^{1/2}$, with bounds controlled by $E_0(t_0)$. In particular,
\begin{equation}\label{eq:section}
\|D(t_0)\|_{H^1_{\mathrm{ax}}}+\|D_t(t_0)\|_2
\le C_nE_0(t_0).
\end{equation}
On $\pi/2-\delta<\theta\le\pi/2$, the equation is uniformly elliptic in
$(t,\theta)$. The local $H^2$ estimate and $H^2\hookrightarrow C^0$ in two
dimensions give
\begin{equation}\label{eq:top-trace}
|D(t_0,\pi/2)|\le C_nE_0(t_0).
\end{equation}

Fix $\bar b\in[0,\frac n2]$ and set $e_*=e_0(\bar b)$, chosen close to $e_0(b(t_0))$.
Set
\[
\begin{aligned}
\mathbf c_b&=(\frac n2\langle e_0(b(t_0)),D(t_0)\rangle,
         \langle e_0(b(t_0)),D_t(t_0)\rangle),\\
\mathbf c_*&=(\frac n2\langle e_*,D(t_0)\rangle,
         \langle e_*,D_t(t_0)\rangle).
\end{aligned}
\]
If $\|e_0(b(t_0))-e_*\|_2\le\delta$, then
\[
|\mathbf c_b-\mathbf c_*|
\le\delta\bigl(\frac{n^2}{4}\|D(t_0)\|_2^2+\|D_t(t_0)\|_2^2\bigr)^{1/2}
\le\delta C_n E_0(t_0) =\delta C_n|\mathbf c_b|.
\]
Choosing $\delta C_n<1$ yields
\begin{equation}\label{eq:fixed-projection}
|\mathbf c_b|\le\frac{|\mathbf c_*|}{1-\delta C_n},\qquad
\bigl(\frac{n^2}{4}\|D(t_0)\|_2^2+\|D_t(t_0)\|_2^2\bigr)^{1/2}
\le\frac{C_n|\mathbf c_*|}{1-\delta C_n}.
\end{equation}
The corresponding bounds for the $H^1_{\mathrm{ax}}\times L^2$ section norm,
equatorial trace, and local $H^2$ norms follow from
\eqref{eq:section}, \eqref{eq:top-trace}, and \eqref{eq:local-L2}.
If $\mathbf c_*=0$, then $D\equiv0$ by \eqref{eq:observation}.

All these solutions satisfy $0\le v\le1$ by Proposition~\ref{prop:sharp}.
With $w_i(t)=v_i(t,\pi/2)$, \eqref{eq:secant} gives
\[
\frac n2\min_{i=1,2}w_i(t_0)^{\frac{2}{n-2}}
\le b(t_0)\le
\frac n2\max_{i=1,2}w_i(t_0)^{\frac{2}{n-2}}.
\]
By \eqref{eq:e0H1}, the map $b\mapsto e_0(b)$ is uniformly continuous on $[0,\frac n2]$.
Choose finitely many relatively open intervals $I_j$ covering $[0,1]$,
with $s_j\in I_j$, narrow enough that
\[
\|e_0(b)-e_0(\frac n2s_j)\|_2<\frac1{2C_n}
\]
whenever $2b/n\in I_j$.
Let $b_j=\frac n2s_j$ and set
\[
\mathcal U_j
:=\left\{v\in\mathcal S:
v(t_0,\pi/2)^{\frac{2}{n-2}}\in I_j\right\}.
\]
By trace continuity, each $\mathcal U_j$ is relatively open in $\mathcal S$,
and $\bigcup_j\mathcal U_j=\mathcal S$. Here $I_j$ is open relative to
$[0,1]$, also when it contains $0$ or $1$. On $\mathcal U_j$, define
\begin{equation*}
\mathbf c_j(v)=
 \bigl(\frac n2\langle e_0(b_j),v(t_0)\rangle,
       \langle e_0(b_j),v_t(t_0)\rangle\bigr).
\end{equation*}
For $v_1,v_2\in\mathcal U_j$, we have
$2b(t_0)/n\in I_j$. Take $e_*=e_0(b_j)$ in
\eqref{eq:fixed-projection}. Since
$\mathbf c_*=\mathbf c_j(v_1)-\mathbf c_j(v_2)$ and
$\|e_0(b(t_0))-e_0(b_j)\|_2<1/(2C_n)$, we obtain
\[
E_0(t_0)=|\mathbf c_b|
\le2|\mathbf c_j(v_1)-\mathbf c_j(v_2)|.
\]
Thus $\mathbf c_j(v_1)=\mathbf c_j(v_2)$ implies $E_0(t_0)=0$,
and hence $D\equiv0$ by \eqref{eq:observation}. This proves injectivity.
By \eqref{eq:section}, \eqref{eq:local-H2}, \eqref{eq:local-L2},
and \eqref{eq:top-trace},
\[
\begin{aligned}
\|(D,D_t)(t_0)\|_{H^1_{\mathrm{ax}}\times L^2}
&\le C_n|\mathbf c_j(v_1)-\mathbf c_j(v_2)|,\\
\|D\|_{H^2([t_0-R,t_0+R]\times\mathbb S^{n-1}_+)}
+\|D(\,\cdot,\pi/2)\|_{C^0([t_0-R,t_0+R])}&\le C_R|\mathbf c_j(v_1)-\mathbf c_j(v_2)|.
\end{aligned}
\]
Since $\|e_0(b_k)\|_2=1$, the Cauchy--Schwarz inequality and the first bound give,
for $v_1,v_2\in\mathcal U_j\cap\mathcal U_k$,
\[
\begin{aligned}
|\mathbf c_k(v_1)-\mathbf c_k(v_2)|
&\le\left(\frac{n^2}{4}\|D(t_0)\|_2^2+\|D_t(t_0)\|_2^2\right)^{1/2}\\
&\le C_n\|(D,D_t)(t_0)\|_{H^1_{\mathrm{ax}}\times L^2}\\
&\le C_n|\mathbf c_j(v_1)-\mathbf c_j(v_2)|.
\end{aligned}
\]
Hence $\mathbf c_k\circ\mathbf c_j^{-1}$ is Lipschitz on
$\mathbf c_j(\mathcal U_j\cap\mathcal U_k)$; taking $k=j$ also gives
the forward bound. Finally, take $t_0=0$. For each $v_0\in\mathcal S$,
choose $j$ with $v_0\in\mathcal U_j$ and set
$\mathcal U=\mathcal U_j$, $\bar b=b_j$, and $e_*=e_0(b_j)$.

\end{proof}

The local estimates also give compactness.
\begin{lemma}\label{lem:solution-compactness}
The space $\mathcal S$ is compact, and $\tau_s v:=v(\,\cdot+s,\cdot)$, $s\in\mathbb R$,
defines a continuous flow on $\mathcal S$. Moreover,
\[
 v_j\to v\ \text{in }\mathcal S
 \quad\Longrightarrow\quad
 v_j(\,\cdot,\pi/2)\to v(\,\cdot,\pi/2)
 \ \text{in }C^0_{\mathrm{loc}}(\mathbb R).
\]
\end{lemma}
\begin{proof}
Write $v_\theta(t,\pi/2)=\frac{n-2}{2}w(t)^{\frac{2}{n-2}}v(t,\pi/2)$. By \eqref{eq:sharp},
$0\le \frac{n-2}{2}w^{\frac{2}{n-2}}\le \frac{n-2}{2}$ and $\left|\left(\frac{n-2}{2}w^{\frac{2}{n-2}}\right)'\right|\le \frac{n-2}{2}w^{\frac{2}{n-2}}$.
Thus \eqref{eq:local-H2} gives uniform $H^2$ bounds on fixed cylinders.

At each section of a fixed $t$, $H^{3/2}\Subset H^1$ and $H^{1/2}\Subset L^2$ give
compactness of $(v,v_t)$ in $H^1_{\mathrm{ax}}\times L^2$. Equicontinuity follows from
\[
\begin{aligned}
\|v(t)-v(s)\|_{H^1_{\mathrm{ax}}}
&\le |t-s|^{1/2}\|v_t\|_{L^2((s,t);H^1_{\mathrm{ax}})},\\
\|v_t(t)-v_t(s)\|_{L^2}
&\le |t-s|^{1/2}\|v_{tt}\|_{L^2((s,t);L^2)}.
\end{aligned}
\]
On a collar of the equator, ellipticity in $(t,\theta)$ and the compact
embedding $H^2\Subset C^0$ give convergence of the equatorial traces on
compact time intervals. A diagonal subsequence converges in the stated
topology and satisfies the equation and boundary condition in the limit.
By the strong maximum principle and the Hopf lemma combined with the Robin condition,
the limit is either $0$ or positive. Thus $\mathcal S$ is compact and
the translation flow is continuous. At a positive limit, elliptic
regularity gives smooth local convergence.
\end{proof}

Recall $\mathcal H(v)$ defined in \eqref{eq:hamiltonian}. For $H\in\mathbb R$, set $$\mathcal S_H:=\{v\in\mathcal S:\mathcal H(v)=H\}. $$

Recall from the introduction and Section~\ref{sec:stationary} that
\[
\begin{aligned}
&b_0=\Theta(0)=2\frac{\Gamma(\frac n4)^2}{\Gamma(\frac{n-2}{4})^2}, &&\qquad
w_c=\left(\frac{2b_0}{n-2}\right)^{\frac{n-2}{2}},\\
&b_c=\frac n2w_c^{\frac{2}{n-2}}=\frac{n}{n-2}b_0,
&&\qquad
H_c=-\frac{n-2}{4(n-1)}w_c^{\frac{2(n-1)}{n-2}}<0,
\end{aligned}
\]
where $\Theta$ is given in \eqref{eq:symbol}, and $v_c$ is the regular stationary solution with trace $w_c$.
The frequency in Theorem~\ref{thm:fowler-family} satisfies
$\Theta(\kappa_c)=b_c$, and $\phi_c=e_0(b_c)/e_0^\partial(b_c)$.

\begin{proposition}\label{prop:energy-rigidity}
Let $v\in\mathcal S$ be nonstationary and write $H=\mathcal H(v)$. There exist
$\varepsilon>0$ and a neighborhood $\mathcal U$ of $v$ in $\mathcal S$ such that
\[
 \mathcal S_H\cap\mathcal U
 =\{v(\,\cdot+s,\cdot): |s|<\varepsilon\}\cap\mathcal U.
\]
Moreover, there is a neighborhood $\mathcal U_c$ of $v_c$ in $\mathcal S$ for which
\[
 \mathcal S_{H_c}\cap\mathcal U_c=\{v_c\}.
\]
\end{proposition}
\begin{proof}
Fix a positive solution $v\in\mathcal S$ and set
$b(t)=\frac n2w(t)^{\frac{2}{n-2}}$.
For $h\in\mathcal C_b$, integration by parts gives 
\begin{equation}\label{eq:H-flux}
\mathrm d\mathcal H_v[h]
:=
\left.\frac{\mathrm d}{\mathrm ds}\right|_{s=0}
\mathcal H(v+s h)
 =\langle v_t,h_t\rangle-\langle v_{tt},h\rangle
 =\Omega(v_t,h).
\end{equation}
For nonstationary $v$, $v_t\in\mathcal C_b\setminus\{0\}$.
Since $\dim\mathcal C_b=2$ and $\Omega(v_t,v_t)=0$, Proposition~\ref{prop:symplectic} gives
\begin{equation}\label{eq:dHkernel}
 \operatorname{ker}\bigl(\mathrm{d}\mathcal H_v|_{\mathcal C_b}\bigr)=\operatorname{span}\{v_t\}.
\end{equation}

Fix nonstationary $v$ and a chart $\mathbf c$ from
Proposition~\ref{prop:coordinates}. Since $v_t\in\mathcal C_b\setminus\{0\}$,
uniqueness in Proposition~\ref{prop:linear-two-dimensionality} gives
$(v_t(0),v_{tt}(0))\ne(0,0)$. On the section space
$H^1_{\mathrm{ax}}\times L^2$, define
\[
\ell(f,g):=\langle v_t(0),f\rangle+\langle v_{tt}(0),g\rangle.
\]
By Cauchy--Schwarz,
$|\ell(f,g)|\le C_v(\|f\|_{H^1_{\mathrm{ax}}}+\|g\|_2)$, so $\ell$ is continuous.
For a function $h$, write $\ell(h):=\ell(h(0),h_t(0))$. Then
\[
\ell(v_t)=\ell(v_t(0),v_{tt}(0))
=\|v_t(0)\|_2^2+\|v_{tt}(0)\|_2^2>0.
\]

If the assertion fails,
there are $v_j\to v$, $\mathcal H(v_j)=\mathcal H(v)$, which are not small
translates of $v$. Set
\[
F_j(s):=\ell\bigl(v_j(s)-v(0),(v_j)_t(s)-v_t(0)\bigr).
\]
Since $v>0$, Lemma~\ref{lem:solution-compactness} gives smooth local
convergence. Hence $F_j(0)\to0$ and
\[
F_j'(s)=\ell\bigl((v_j)_t(s),(v_j)_{tt}(s)\bigr)
\longrightarrow\ell\bigl(v_t(s),v_{tt}(s)\bigr)
\]
locally uniformly in $s$. Since $\ell(v_t)>0$, choose $\delta>0$ so that,
for all large $j$, $F_j'(s)\ge\frac12\ell(v_t)$ on $[-\delta,\delta]$ and
$|F_j(0)|<\frac{\delta}{2}\ell(v_t)$. Then
\[
\begin{aligned}
F_j(-\delta)&\le F_j(0)-\frac{\delta}{2}\ell(v_t)<0,\\
F_j(\delta)&\ge F_j(0)+\frac{\delta}{2}\ell(v_t)>0.
\end{aligned}
\]
Thus $F_j$ has a unique zero $s_j\in(-\delta,\delta)$, and the mean value
theorem gives
\[
|s_j|\le\frac{2|F_j(0)|}{\ell(v_t)}\longrightarrow0.
\]
By continuity of the translation flow, $\tau_{s_j}v_j\to v$, and
$\mathcal H(\tau_{s_j}v_j)=\mathcal H(v_j)=\mathcal H(v)$.
For large $j$, $\tau_{s_j}v_j\ne v$, since otherwise
$v_j=\tau_{-s_j}v$ would be a small translate of $v$.
Relabel $\tau_{s_j}v_j$ as $v_j$ and set $D_j=v_j-v$. Then
\begin{equation}\label{eq:translation-normalization}
\ell(D_j)=F_j(s_j)=0.
\end{equation}
Set
\[
 \varepsilon_j=|\mathbf c(v_j)-\mathbf c(v)|>0,
 \qquad U_j=D_j/\varepsilon_j.
\]
Then
\[
\frac{n^2}{4}\bigl|\langle e_*,U_j(0)\rangle\bigr|^2
+\bigl|\langle e_*,(U_j)_t(0)\rangle\bigr|^2
=1.
\]
By Lemma~\ref{lem:secant}, $D_j$ satisfies \eqref{eq:secant-equation} with
coefficient $b_j$ from \eqref{eq:secant}. Estimates
\eqref{eq:fixed-projection}, \eqref{eq:observation}, and \eqref{eq:weighted-form}
give uniform weighted energy bounds for $U_j$. Passing to a subsequence, we have
\[
 U_j\rightharpoonup U\quad\text{in }H^2_{\mathrm{loc}},
 \qquad
 (U_j,(U_j)_t)\longrightarrow(U,U_t)
 \quad\text{in }C_{\mathrm{loc}}
 \bigl(\mathbb R;H^1_{\mathrm{ax}}\times L^2\bigr).
\]
Passing to the limit in the weak formulation using the local convergence
of the coefficients shows that $U$ satisfies the limiting equation.
On each bounded time interval, local convergence of the weights and weak
lower semicontinuity give a weighted energy bound for $U$, with a constant
independent of the interval. Letting these intervals exhaust the real
line, we obtain $U\in\mathcal C_b$ and
\[
\frac{n^2}{4}|\langle e_*,U(0)\rangle|^2
+|\langle e_*,U_t(0)\rangle|^2=1.
\]
Hence, $U\ne0$.

The Hamiltonian is $C^2$ on a positive value interval. At $t=0$,
\[
 0=\frac{\mathcal H(v_j)-\mathcal H(v)}{\varepsilon_j}
   \longrightarrow \mathrm{d}\mathcal H_v[U].
\]
By \eqref{eq:dHkernel}, $U=\alpha v_t$. Passing
\eqref{eq:translation-normalization} to the limit gives $\alpha=0$, contradicting
$U\ne0$. This proves the first assertion.

At $v_c$, $\mathrm d\mathcal H_{v_c}=0$ and $A_{b_c}\phi_c=-\kappa_c^2\phi_c$.
Since $0<w_c\le1$, we have $0<b_c=\frac n2w_c^{\frac{2}{n-2}}\le\frac n2$.
The Robin spectral gap gives
\[
\lambda_0(b_c)=-\kappa_c^2<0,\qquad
\lambda_j(b_c)\ge\lambda_1(b_c)>\frac{n^2}{4}
\quad\mbox{for every }j\ge1.
\]
For $U\in\mathcal C_{b_c}$, projecting $U_{tt}=A_{b_c}U$ onto the
time-independent eigenfunctions $e_j(b_c)$ gives
\[
\begin{gathered}
\frac{\mathrm d^2}{\mathrm dt^2}\langle U(t),e_j(b_c)\rangle
=\lambda_j(b_c)\langle U(t),e_j(b_c)\rangle,\qquad j\ge0,\\
\langle U(t),e_j(b_c)\rangle=
\begin{cases}
A_0\cos(\kappa_ct)+B_0\sin(\kappa_ct),&j=0,\\
A_j e^{\sqrt{\lambda_j(b_c)}\,t}
+B_j e^{-\sqrt{\lambda_j(b_c)}\,t},&j\ge1.
\end{cases}
\end{gathered}
\]
Moreover,
\[
\int_{\mathbb R}\omega(t)|\langle U(t),e_j(b_c)\rangle|^2\,\mathrm dt
\le\int_{\mathbb R}\omega(t)\|U(t)\|_2^2\,\mathrm dt<\infty,
\]
where $\omega(t)=e^{-2(\frac n2-\eta)|t|}$, and $\eta$ is constant such that $\frac13\le\eta\le\frac n4$. 
For $j\ge1$, $\sqrt{\lambda_j(b_c)}>\frac n2>\frac n2-\eta$.
If $A_j\ne0$, then
\[
\omega(t)|\langle U(t),e_j(b_c)\rangle|^2
\sim |A_j|^2e^{2(\sqrt{\lambda_j(b_c)}-\frac n2+\eta)t}
\quad\mbox{as }t\to+\infty,
\]
contradicting integrability. Thus $A_j=0$, and the same argument as
$t\to-\infty$ gives $B_j=0$.
By completeness of $\{e_j(b_c)\}_{j\ge0}$, 
\[
\mathcal C_{b_c}
=\operatorname{span}\{\phi_c(\theta)\cos(\kappa_c t), \phi_c(\theta)\sin(\kappa_c t)\}.
\]
For $U=\phi_c(A\cos\kappa_c t+B\sin\kappa_c t)\ne0$,
\begin{equation}\label{eq:cylinder-quadratic}
 \frac12\mathrm{d}^2\mathcal H_{v_c}[U,U]= \frac12\left.\frac{\mathrm d^2}{\mathrm ds^2}\right|_{s=0}
\mathcal H(v_c+sU)
 =\frac{\kappa_c^2}{2}\|\phi_c\|_2^2(A^2+B^2)>0.
\end{equation}
If $v_j\ne v_c$, $\mathcal H(v_j)=H_c$, and $v_j\to v_c$, then set
$\varepsilon_j=|\mathbf c(v_j)-\mathbf c(v_c)|$ and
$U_j=(v_j-v_c)/\varepsilon_j$. By \eqref{eq:observation},
\eqref{eq:fixed-projection}, and local $H^2$ estimates, a subsequence
converges to a nonzero $U$ as above. Since $\mathrm d\mathcal H_{v_c}=0$, Taylor's formula gives
\[
0=\frac{\mathcal H(v_j)-H_c}{\varepsilon_j^2}
  \longrightarrow \frac12\mathrm{d}^2\mathcal H_{v_c}[U,U],
\]
contradicting \eqref{eq:cylinder-quadratic}. Hence $v_c$ is isolated in its
energy level.
\end{proof}

\subsection{Periodic continuation and bifurcation at the cylinder}
Recall that $v\in\mathcal S$ is periodic if there exists $T>0$ such that
\[
v(t+T,\theta)=v(t,\theta)
\qquad\text{for all }(t,\theta)\in\mathbb R\times[0,\pi/2].
\]
If $v_t\not\equiv0$, its least period is the smallest such $T>0$.
For a $T$-periodic boundary trace $w$, set $W(s)=w(Ts)$, so that $W$ is
$1$-periodic. Denote by $\mathcal P_T$ the Dirichlet-to-Neumann operator
in this rescaled variable. Then
\begin{equation}\label{eq:periodic-symbol}
 \widehat{\mathcal P_T W}(j)
 =\Theta(2\pi j/T)\widehat W(j),\qquad j\in\mathbb Z.
\end{equation}

\begin{proposition}\label{prop:periodic-neighborhood}
Let $v_0\in\mathcal S$ be a nonstationary periodic solution with least period $T_0$.
There is a neighborhood $\mathcal U$ of $v_0$ in $\mathcal S$ which is a
real-analytic two-dimensional manifold in the local Cauchy coordinates of
Proposition~\ref{prop:coordinates}. Every $v\in\mathcal U$ is periodic, and its
least period $T:\mathcal U\to(0,\infty)$ is real analytic, with $T(v_0)=T_0$.
\end{proposition}

\begin{proof}
Fix a nonstationary periodic solution of least period $T_0$, and put
$W_0(s)=w(T_0s)$ for $s\in\mathbb{R}/\mathbb Z$.  For \(T>0\) and a positive \(1\)-periodic function \(W\), define
\[
\mathcal F(T,W)=\mathcal P_T W-\frac{n-2}{2}W^{\frac{n}{n-2}},
\]
where $\mathcal P_T$ is given by \eqref{eq:periodic-symbol}. Then \(\mathcal F(T_0,W_0)=0\).
For $T\in[T_-,T_+]\Subset(0,\infty)$, consider
\[
 \mathcal L_T=-T^{-2}\partial_s^2-\Delta_{\mathbb S^{n-1}_+}+\frac{(n-2)^2}{4}
\]
on $(\mathbb R/\mathbb Z)\times\mathbb S^{n-1}_+$ with the $T$-independent
Dirichlet domain $H^2\cap H^1_0$, where $H^1_0$ imposes zero trace on
$(\mathbb R/\mathbb Z)\times\partial\mathbb S^{n-1}_+$.
All norms in this paragraph are on this product. For $u\in H^2\cap H^1_0$,
integration by parts gives
\[
\begin{aligned}
(\mathcal L_Tu,u)_{L^2}
&=T^{-2}\|u_s\|_{L^2}^2
+\|\nabla_{\mathbb S^{n-1}}u\|_{L^2}^2
+\frac{(n-2)^2}{4}\|u\|_{L^2}^2\ge C(T_+,n)\|u\|_{H^1}^2.
\end{aligned}
\]
By Lax--Milgram and the $H^2$ estimates,
each $f\in L^2$ has a unique solution $u\in H^2\cap H^1_0$ of
$\mathcal L_Tu=f$, with
\[
\begin{aligned}
\|u\|_{H^2}
&\le C_{n,T_-,T_+}\bigl(\|f\|_{L^2}+\|u\|_{L^2}\bigr)\le C_{n,T_-,T_+}\|f\|_{L^2}.
\end{aligned}
\]
Thus $\mathcal L_T:H^2\cap H^1_0\to L^2$ is uniformly invertible.
Since
\[
\mathcal L_T-\mathcal L_{T_0}
=-\bigl(T^{-2}-T_0^{-2}\bigr)\partial_s^2,
\]
the map $T\mapsto\mathcal L_T$ is real analytic in operator norm.
For $T$ sufficiently close to $T_0$,
\[
\begin{aligned}
&\left\|(\mathcal L_T-\mathcal L_{T_0})\mathcal L_{T_0}^{-1}
\right\|_{L^2\to L^2}<1,\\
&\mathcal L_T^{-1}
=\mathcal L_{T_0}^{-1}\sum_{j=0}^{\infty}
\left[-(\mathcal L_T-\mathcal L_{T_0})\mathcal L_{T_0}^{-1}\right]^j.
\end{aligned}
\]
The series converges in operator norm from $L^2$ to $H^2\cap H^1_0$
and gives the real-analytic dependence of $\mathcal L_T^{-1}$ on $T$.

For $\ell>5/2$ and $W\in H^\ell(\mathbb R/\mathbb Z)$ (the choice of $\ell>5/2$ provides convenient \(C^2\) regularity), let $v_{T,W}$
denote the $1$-periodic solution, regular at $\theta=0$, with boundary trace
$W$, so that
\[
\begin{cases}
T^{-2}(v_{T,W})_{ss}+m^{-1}\bigl(m(v_{T,W})_\theta\bigr)_\theta-\frac{(n-2)^2}{4}v_{T,W}=0
& \text{in }(\mathbb R/\mathbb Z)\times(0,\pi/2),\\
v_{T,W}(s,\pi/2)=W(s)
& \text{for }s\in\mathbb R/\mathbb Z.
\end{cases}
\]
Elliptic boundary regularity gives, uniformly for $T\in[T_-,T_+]$,
\[
\begin{aligned}
\|v_{T,W}\|_{H^{\ell+1/2}}
&\le C_{n,\ell,T_-,T_+}\|W\|_{H^\ell},\\
\|\mathcal L_T^{-1}f\|_{H^{\ell+1/2}}
&\le C_{n,\ell,T_-,T_+}\|f\|_{H^{\ell-3/2}}.
\end{aligned}
\]
Interior norms are on the product, and boundary norms on $\mathbb R/\mathbb Z$.
Since $\partial_s^2\mathcal L_{T_0}^{-1}:H^{\ell-3/2}\to H^{\ell-3/2}$
is bounded, the same Neumann series gives real-analytic dependence of
$\mathcal L_T^{-1}$ in
$\mathcal L(H^{\ell-3/2},H^{\ell+1/2}\cap H^1_0)$, where \(\mathcal L(E,Y)\) denotes the space of bounded linear operators from \(E\) to \(Y\).
Moreover,
\[
\begin{aligned}
\mathcal L_T(v_{T,W}-v_{T_0,W})
&=(T^{-2}-T_0^{-2})\partial_s^2v_{T_0,W},\\
(v_{T,W}-v_{T_0,W})(s,\pi/2)&=0.
\end{aligned}
\]
Hence
\[
v_{T,W}=v_{T_0,W}
+(T^{-2}-T_0^{-2})\mathcal L_T^{-1}\partial_s^2v_{T_0,W}.
\]
Thus the Poisson operator $W\mapsto v_{T,W}$ depends real analytically on
$T$ in $\mathcal L(H^\ell,H^{\ell+1/2})$. On axial functions, the normal-trace map is bounded
from $H^{\ell+1/2}$ to $H^{\ell-1}$ and independent of $T$, so
\[
\begin{aligned}
\mathcal P_TW&=(v_{T,W})_\theta(\cdot,\pi/2),\\
\|\mathcal P_TW\|_{H^{\ell-1}}
&\le C\|v_{T,W}\|_{H^{\ell+1/2}}\le C\|W\|_{H^\ell}.
\end{aligned}
\]
Consequently, $T\mapsto\mathcal P_T$ is real analytic in
$\mathcal L(H^\ell(\mathbb R/\mathbb Z),H^{\ell-1}(\mathbb R/\mathbb Z))$.
Thus
$\mathcal F:\mathbb R_+\times H^\ell\to H^{\ell-1}$ is real analytic on the
positive cone \(\left\{W\in H^\ell(\mathbb R/\mathbb Z): \min_{\mathbb R/\mathbb Z}W>0\right\}\) for $\ell>5/2$, with
\[
 D_W\mathcal F(T,W)=\mathcal P_T-\frac n2W^{\frac{2}{n-2}}.
\]
By $\Theta(\xi)=|\xi|+O(|\xi|^{-1})$, this is a self-adjoint
first-order elliptic operator and hence Fredholm of index zero.

Let $\Pi_\perp$ be the fixed $L^2$ projection onto
$\operatorname{span}\{(W_0)_s\}^\perp$. Translation invariance gives $\mathcal F(T_0,W_0(\cdot+\sigma))=0$. Hence, $D_W\mathcal F(T_0,W_0)(W_0)_s=0$. By self-adjointness,
\[
\int_0^1
\bigl(D_W\mathcal F(T_0,W_0)h\bigr)(W_0)_s\,\mathrm ds
=
\int_0^1h\,D_W\mathcal F(T_0,W_0)(W_0)_s\,\mathrm ds
=0.
\]
Hence the range of this linearized operator lies in \(\operatorname{span}\{(W_0)_s\}^{\perp}\). 

Choose a neighborhood \(\mathcal U\) of \(v_0\) and the associated normalized eigenfunction \(e_{\mathrm{ref}}\) as in Proposition~\ref{prop:coordinates}.
Define the Cauchy-coordinate map
\[
\mathbf c_{\mathrm{ref}}(T,W)=
\left(
 \frac{n}{2}\langle e_{\mathrm{ref}},v_{T,W}(0)\rangle,
 T^{-1}\langle e_{\mathrm{ref}},\partial_sv_{T,W}(0)\rangle
\right).
\]
The analytic dependence of the Poisson operator and the trace theorem imply that \((T,W)\mapsto\mathbf c_{\mathrm{ref}}(T,W)\) is real analytic. 
The derivative at $(T_0,W_0)$ of
\begin{equation}\label{eq:periodic-map}
(T,W)\longmapsto\bigl(\Pi_\perp\mathcal F(T,W),\mathbf c_{\mathrm{ref}}(T,W)\bigr)
\end{equation}
has index $0+1+1-2=0$: one period variable, one codomain projection,
and two Cauchy coordinates. Using the $1$-periodicity of $W$, and the translation invariance and the self-adjointness of $\mathcal P_T$, we have
\begin{equation}\label{eq:translation-identity}
\int_0^1\mathcal F(T,W)W_s\,\mathrm ds=0.
\end{equation}
If $(\delta T,\delta W)$ belongs to the kernel of the derivative of the map defined in 
\eqref{eq:periodic-map}, then
\[
 D\mathcal F_{(T_0,W_0)}(\delta T,\delta W)=\lambda(W_0)_s
\]
for some $\lambda\in\mathbb R$. Differentiating
\eqref{eq:translation-identity} at $(T_0,W_0)$ and using \(\mathcal F(T_0,W_0)=0\), we obtain
\[
 0=\int_0^1
 D\mathcal F_{(T_0,W_0)}(\delta T,\delta W)(W_0)_s\,\mathrm ds
 =\lambda\|(W_0)_s\|_{L^2}^2,
\]
so $\lambda=0$.

Set $T_\varepsilon=T_0+\varepsilon\delta T$ and
$W_\varepsilon=W_0+\varepsilon\delta W$, and define
\[
\delta v(s,\theta)
:=\left.\frac{\mathrm d}{\mathrm d\varepsilon}\right|_{\varepsilon=0}
v_{T_\varepsilon,W_\varepsilon}(s,\theta).
\]
Then $\delta v(s+1,\theta)=\delta v(s,\theta)$ and
$\delta v(s,\pi/2)=\delta W(s)$. Since
$v_{T_0,W_0}(s,\theta)=v(T_0s,\theta)$, the chain rule gives
\[
\begin{aligned}
\phi(t,\theta)
&:=\left.\frac{\mathrm d}{\mathrm d\varepsilon}\right|_{\varepsilon=0}
v_{T_\varepsilon,W_\varepsilon}(t/T_\varepsilon,\theta)\\
&=\delta v(t/T_0,\theta)
-\frac{\delta T}{T_0^2}\,t\,\partial_s v_{T_0,W_0}(t/T_0,\theta)\\
&=\delta v(t/T_0,\theta)
-\frac{\delta T}{T_0}\,t\,v_t(t,\theta).
\end{aligned}
\]
It satisfies
\[
\begin{cases}
\phi_{tt}+m^{-1}(m\phi_\theta)_\theta-\frac{(n-2)^2}{4}\phi=0
& \text{in }\mathbb R\times(0,\pi/2),\\
\phi_\theta(t,\pi/2)=b(t)\phi(t,\pi/2)
& \text{for }t\in\mathbb R,
\end{cases}
\]
with $b=\frac n2w^{\frac{2}{n-2}}$ and regularity at $\theta=0$.
Since $\phi$ grows at most linearly, $\phi\in\mathcal C_b$.
At $t=0$,
\[
\phi(0)=\delta v(0),\qquad
\phi_t(0)=T_0^{-1}(\delta v)_s(0)-\frac{\delta T}{T_0}v_t(0).
\]
Since $e_{\mathrm{ref}}$ is fixed and
$\partial_s v_{T_0,W_0}(0)=T_0v_t(0)$, the kernel assumption in
\eqref{eq:periodic-map} gives
\[
\begin{aligned}
0&=D\mathbf c_{\mathrm{ref},(T_0,W_0)}(\delta T,\delta W)\\
&=\begin{pmatrix}
\frac n2\langle e_{\mathrm{ref}},\delta v(0)\rangle\\[2pt]
T_0^{-1}\langle e_{\mathrm{ref}},(\delta v)_s(0)\rangle
-\frac{\delta T}{T_0}\langle e_{\mathrm{ref}},v_t(0)\rangle
\end{pmatrix}\\
&=\begin{pmatrix}
\frac n2\langle e_{\mathrm{ref}},\phi(0)\rangle\\[2pt]
\langle e_{\mathrm{ref}},\phi_t(0)\rangle
\end{pmatrix}.
\end{aligned}
\]
Hence $\langle e_{\mathrm{ref}},\phi(0)\rangle=
\langle e_{\mathrm{ref}},\phi_t(0)\rangle=0$.
Thus \eqref{eq:fixed-projection} and
Proposition~\ref{prop:linear-two-dimensionality} imply $\phi=0$.
The periodicity of $\delta v$ gives
\[
\phi(t+T_0)-\phi(t)=-\delta T\,v_t(t).
\]
Since $v_t\not\equiv0$, $\delta T=0$ and then $\delta v=0$. Hence $\delta W=0$.
Thus the derivative in \eqref{eq:periodic-map} is injective. Since this derivative is Fredholm of index zero, it is also surjective. Since it is a bounded linear operator, the bounded inverse theorem gives a bounded inverse. Hence, it is an isomorphism.

The analytic implicit function theorem gives, for every
$\mathbf c\in\mathbb R^2$ sufficiently close to
$\mathbf c_{\mathrm{ref}}(T_0,W_0)$, a unique pair $(T,W)$ near
$(T_0,W_0)$ satisfying
\[
\Pi_\perp\mathcal F(T,W)=0,\qquad
\mathbf c_{\mathrm{ref}}(T,W)=\mathbf c.
\]
Moreover, $(T,W)$ depends real analytically on $\mathbf c$.
Since $\ker\Pi_\perp=\operatorname{span}\{(W_0)_s\}$, there is
$\lambda(T,W)\in\mathbb R$ such that
$\mathcal F(T,W)=\lambda(T,W)(W_0)_s$.
By \eqref{eq:translation-identity},
\[
0=\int_0^1\mathcal F(T,W)W_s\,\mathrm ds
=\lambda(T,W)\int_0^1(W_0)_sW_s\,\mathrm ds.
\]
For $W$ sufficiently close to $W_0$ in $H^\ell$,
\[
\int_0^1(W_0)_sW_s\,\mathrm ds
\ge\frac12\int_0^1|(W_0)_s|^2\,\mathrm ds>0.
\]
Hence $\lambda(T,W)=0$, and therefore
$\mathcal F(T,W)=0$.

After shrinking the neighborhood, $H^\ell\hookrightarrow C^0$ gives
$\min W\ge\frac12\min W_0>0$, and the maximum principle yields
$v_{T,W}>0$.
Hence $v_{T,W}(t/T,\theta)$ is a positive entire $T$-periodic solution of
\eqref{eq:cylinder}, with smooth rotational extension at $\theta=0$
and bounded trace. It therefore belongs to $\mathcal S$.

By Proposition~\ref{prop:coordinates}, these exhaust all nearby solutions.
After shrinking the neighborhood, $T(v)$ is the least period. Indeed, suppose 
$W_j\to W_0$ and $W_j$ has proper periods $1/m_j$, with the integers $m_j\ge2$. If $\{m_j\}$ is a bounded sequence, then locally uniform convergence of $W_j$  would give a smaller period of $W_0$. If $m_j\to\infty$, the uniform $C^1$ bound gives $\operatorname{osc}W_j\le C/m_j\to0$, which is also a contradiction.
\end{proof}

We next apply simple-eigenvalue bifurcation at $v_c$.

\begin{proposition}\label{prop:cylinder-neighborhood}
The function $\phi_c$ satisfies
$A_{b_c}\phi_c=-\kappa_c^2\phi_c$, and the following hold.
\begin{enumerate}
\item There are $\varepsilon_0>0$ and a real-analytic family of even periodic
solutions $\{v_\varepsilon:|\varepsilon|<\varepsilon_0\}$ with
$v_0=v_c$.  For $\varepsilon\ne0$, let $T(\varepsilon)$ be the least period
and set
\[
 \widetilde v_\varepsilon(s,\theta)
 :=v_\varepsilon(T(\varepsilon)s,\theta),
 \qquad s\in\mathbb R/\mathbb Z.
\]
Then $T(\varepsilon)\to2\pi/\kappa_c$ as $\varepsilon\to0$, and for every integer $j\ge0$,
\begin{equation}\label{eq:cylinder-branch-statement}
 \bigl\|\widetilde v_\varepsilon-v_c
 -\varepsilon \phi_c(\theta)\cos(2\pi s)\bigr\|_
 {C^j((\mathbb R/\mathbb Z)\times\overline{\mathbb S^{n-1}_+})}
 \le C_j\varepsilon^2.
\end{equation}
\item Every solution sufficiently close to $v_c$ in $\mathcal S$ is either
$v_c$ or a time translate of one of the $v_\varepsilon$.
\item Along the branch,
\[
 \mathcal H(v_\varepsilon)
 =H_c+\frac{\kappa_c^2\|\phi_c\|_2^2}{2}\varepsilon^2+O(\varepsilon^3).
\]
Consequently, there are a neighborhood $\mathcal U$ of $v_c$ in
$\mathcal S$ and $\delta>0$ such that, for every $H\in(H_c,H_c+\delta)$,
there is exactly one nonstationary periodic translation orbit in $\mathcal S_H\cap\mathcal U$.
\end{enumerate}
\end{proposition}

\begin{proof}
For $x\ge0$, let $\xi=\sqrt{x}$ and $\psi_{\xi}$ be the solution of \eqref{eq:modexi}. Denote $f_x=\psi_{\xi}$. 
Multiplying \eqref{eq:modexi} by $mf_x$ and integrating by parts gives
\[
\Theta(\sqrt x)=f_x'(\pi/2)
=\int_0^{\pi/2}m\left((f_x')^2+
\left(\frac{(n-2)^2}{4}+x\right)f_x^2\right)\,\mathrm d\theta.
\]
Differentiating it and using \eqref{eq:modexi}, we therefore obtain
\[
\frac{\mathrm d}{\mathrm dx}\Theta(\sqrt x)
=\int_0^{\pi/2}mf_x^2\,\mathrm d\theta>0.
\]

Restrict $\mathcal F$ to even functions on $\mathbb{R}/\mathbb Z$.
Set $T_c:=2\pi/\kappa_c$ and note that $\mathcal F(T,w_c)=0$ for all $T>0$.
For $h\in H^\ell(\mathbb R/\mathbb Z)$ even, write
$h(s)=a_0+\sum_{j\ge1}a_j\cos(2\pi js)$. Then
\[
\begin{aligned}
D_W\mathcal F(T_c,w_c)h
&=\mathcal P_{T_c}h-b_ch=(b_0-b_c)a_0+
\sum_{j\ge1}\bigl(\Theta(j\kappa_c)-b_c\bigr)a_j\cos(2\pi js).
\end{aligned}
\]
Since $b_0-b_c=-2b_0/(n-2)<0$, $\Theta(\kappa_c)=b_c$, and
$\Theta(j\kappa_c)>b_c$ for $j\ge2$, we obtain
\[
\operatorname{ker}D_W\mathcal F(T_c,w_c)
=\operatorname{span}\{\cos(2\pi s)\}.
\]
Its range in the even subspace of $H^{\ell-1}(\mathbb R/\mathbb Z)$
has codimension one and is $L^2$-orthogonal to $\cos(2\pi s)$.
Moreover, the derivative identity above gives $\Theta'(\kappa_c)>0$, so
\[
\left.\frac{\mathrm{d}}{\mathrm{d} T}\Theta(2\pi/T)\right|_{T=T_c}
=-\frac{2\pi}{T_c^2}\Theta'(\kappa_c)
=-\frac{\kappa_c}{T_c}\Theta'(\kappa_c)<0.
\]
Hence the transversality condition holds:
\[
\begin{aligned}
\partial_TD_W\mathcal F(T_c,w_c)[\cos(2\pi s)]
&=-\frac{\kappa_c}{T_c}\Theta'(\kappa_c)\cos(2\pi s)\notin\operatorname{Range}D_W\mathcal F(T_c,w_c),
\end{aligned}
\]
since it is not $L^2$-orthogonal to $\cos(2\pi s)$.

By the simple-eigenvalue bifurcation theorem of
Crandall--Rabinowitz~\cite{CR}*{Theorem~1.7}, there is a $C^1$ curve
$\varepsilon\mapsto(T(\varepsilon),W_\varepsilon)$ of even solutions
of $\mathcal F(T,W)=0$, with
\[
T(0)=T_c,\qquad
W_\varepsilon=w_c+\varepsilon\cos(2\pi s)+o_{H^\ell}(\varepsilon).
\]
Every nearby even solution with $W\ne w_c$ lies on this curve.
We next obtain a real-analytic parametrization by the first
cosine coefficient. Write
\[
W=w_c+\varepsilon\cos(2\pi s)+z,\qquad
z(-s)=z(s),\qquad
\int_0^1z(s)\cos(2\pi s)\,\mathrm ds=0.
\]
On the even subspaces orthogonal to $\cos(2\pi s)$,
$D_W\mathcal F(T_c,w_c):H^\ell\to H^{\ell-1}$ is an isomorphism.
Applying the analytic implicit function theorem to the component of
$\mathcal F$ orthogonal to this mode gives a unique real-analytic
$z=z(T,\varepsilon)$ such that
\[
\mathcal F\bigl(T,w_c+\varepsilon\cos(2\pi s)+z(T,\varepsilon)\bigr)
\in\operatorname{span}\{\cos(2\pi s)\}.
\]
Since $\mathcal F(T,w_c)=0$, uniqueness gives $z(T,0)=0$.
Moreover,
$D_W\mathcal F(T,w_c)[\cos(2\pi s)]
=(\Theta(2\pi/T)-b_c)\cos(2\pi s)$, so differentiating the
complementary equation at $\varepsilon=0$ gives
$\partial_\varepsilon z(T,0)=0$.
Thus there is a real-analytic scalar function $g$ such that
\[
\begin{aligned}
\mathcal F\bigl(T,w_c+\varepsilon\cos(2\pi s)+z(T,\varepsilon)\bigr)
&=\varepsilon g(T,\varepsilon)\cos(2\pi s),\\
g(T,0)&=\Theta(2\pi/T)-b_c,
\end{aligned}
\]
where the first equation follows from $\mathcal F(T,w_c)=0$ and the second equation follows from differentiating the first one at $\varepsilon=0$. 
In particular,
\[
g(T_c,0)=0,\qquad
\partial_Tg(T_c,0)=-\frac{\kappa_c}{T_c}\Theta'(\kappa_c)\ne0.
\]
The analytic implicit function theorem gives a unique real-analytic
$T=T(\varepsilon)$ near $T_c$ with $T(0)=T_c$ and
$g(T(\varepsilon),\varepsilon)=0$. Set
\[
W_\varepsilon=w_c+\varepsilon\cos(2\pi s)
+z(T(\varepsilon),\varepsilon).
\]
Then $\mathcal F(T(\varepsilon),W_\varepsilon)=0$.
By uniqueness in the two implicit function theorems, every even solution
sufficiently close to $(T_c,w_c)$ in $\mathbb R\times H^\ell$ with
$W\ne w_c$ has this form. The orthogonality of $z$ gives
\begin{equation}
 2\int_0^1\bigl(W_\varepsilon(s)-w_c\bigr)\cos(2\pi s)\,\mathrm ds
 =\varepsilon.
\end{equation}
Since $z(T,0)=\partial_\varepsilon z(T,0)=0$,
\begin{equation}\label{eq:cylinder-branch}
T(\varepsilon)=T_c+O(\varepsilon),\qquad
W_\varepsilon=w_c+\varepsilon\cos(2\pi s)+O(\varepsilon^2),
\end{equation}
initially in $H^\ell$. Since $\ell>5/2$ is arbitrary, local uniqueness
identifies the branches obtained at different Sobolev orders.
Sobolev embedding gives the remainder estimate in every fixed smooth norm.
If $W_\varepsilon(s+\tau)=W_\varepsilon(s)$ for some $\tau>0$, then
periodicity and evenness give
\[
\begin{aligned}
\varepsilon
&=2\int_0^1\bigl(W_\varepsilon(s+\tau)-w_c\bigr)
\cos(2\pi s)\,\mathrm ds=\varepsilon\cos(2\pi\tau).
\end{aligned}
\]
For $\varepsilon\ne0$, this implies $\tau\in\mathbb N$; hence the least
period in $s$ is one. By recalling 
\[
\begin{aligned}
v_\varepsilon(t,\theta)
&:=v_{T(\varepsilon),W_\varepsilon}
\!\left(\frac{t}{T(\varepsilon)},\theta\right),\qquad
v_\varepsilon(t,\pi/2)
=W_\varepsilon\!\left(\frac{t}{T(\varepsilon)}\right),
\end{aligned}
\]
we prove part (i).

Use $e_0(b_c)$ as the coordinate reference and, for $\varphi\in\mathbb R/2\pi\mathbb Z$, set
\[
 v_{\varepsilon,\varphi}(t,\theta)
 :=v_\varepsilon\!\left(t+\frac{T(\varepsilon)}{2\pi}\varphi,\theta\right).
\]
Then
\[
\mathbf c(v_{\varepsilon,\varphi})=\mathbf c(v_c)
+\varepsilon\langle e_0(b_c),\phi_c\rangle
\left(\frac n2\cos\varphi,-\kappa_c\sin\varphi\right)+O(\varepsilon^2).
\]
The half-period translate remains even, since
\[
\widetilde v_\varepsilon(-s+1/2,\theta)
=\widetilde v_\varepsilon(s-1/2,\theta)
=\widetilde v_\varepsilon(s+1/2,\theta).
\]
Moreover, by $1$-periodicity,
\[
\begin{aligned}
2\int_0^1\bigl(W_\varepsilon(s+1/2)-w_c\bigr)\cos(2\pi s)\,\mathrm ds
&=2\int_{1/2}^{3/2}\bigl(W_\varepsilon(s)-w_c\bigr)
\cos\bigl(2\pi(s-1/2)\bigr)\,\mathrm ds\\
&=-2\int_0^1\bigl(W_\varepsilon(s)-w_c\bigr)\cos(2\pi s)\,\mathrm ds\\
&=-\varepsilon.
\end{aligned}
\]
Thus $(T(\varepsilon),W_\varepsilon(\cdot+1/2))$ is a nearby even
solution with first cosine coefficient $-\varepsilon$. Local uniqueness
of the branch and uniqueness of the Dirichlet extension give
\[
 T(-\varepsilon)=T(\varepsilon),\qquad
 \widetilde v_{-\varepsilon}(s,\theta)=\widetilde v_\varepsilon(s+1/2,\theta).
\]
Returning to the original time variable,
\[
\begin{aligned}
v_{-\varepsilon}(t,\theta)
&=v_\varepsilon\!\left(t+\frac{T(\varepsilon)}2,\theta\right),\\
v_{-\varepsilon,\varphi+\pi}(t,\theta)
&=v_\varepsilon\!\left(t+\frac{T(\varepsilon)}{2\pi}\varphi
+T(\varepsilon),\theta\right)
=v_{\varepsilon,\varphi}(t,\theta).
\end{aligned}
\]
Consequently, it suffices to take $0\le\varepsilon<\varepsilon_0$.
Set $x_1=\varepsilon\cos\varphi$ and $x_2=\varepsilon\sin\varphi$.
These parameters form a disk, and the coordinate expansion becomes
\[
\mathbf c(v_{\varepsilon,\varphi})-\mathbf c(v_c)
=\langle e_0(b_c),\phi_c\rangle
\begin{pmatrix}\frac n2&0\\0&-\kappa_c\end{pmatrix}
\binom{x_1}{x_2}+R(\varepsilon,\varphi).
\]
Analyticity of the branch and smoothness in $\varphi$ give, uniformly in
$\varphi$,
\[
|R|+\varepsilon|\partial_\varepsilon R|+|\partial_\varphi R|
\le C\varepsilon^2.
\]
For $\varepsilon>0$,
\[
\partial_{x_1}=\cos\varphi\,\partial_\varepsilon
-\frac{\sin\varphi}{\varepsilon}\partial_\varphi,
\qquad
\partial_{x_2}=\sin\varphi\,\partial_\varepsilon
+\frac{\cos\varphi}{\varepsilon}\partial_\varphi,
\]
so $|D_{(x_1,x_2)}R|\le C\varepsilon$. Hence
$(x_1,x_2)\mapsto\mathbf c(v_{\varepsilon,\varphi})$, with value
$\mathbf c(v_c)$ at $(0,0)$, is $C^1$ and has an invertible derivative
there, since $\langle e_0(b_c),\phi_c\rangle>0$ and $\kappa_c>0$.
The inverse function theorem shows that its image contains a neighborhood
of $\mathbf c(v_c)$ in $\mathbb R^2$.
For every $v\in\mathcal S$ sufficiently close to $v_c$, we can therefore
choose $\varepsilon,\varphi$ with
$\mathbf c(v)=\mathbf c(v_{\varepsilon,\varphi})$.
Injectivity in Proposition~\ref{prop:coordinates} gives
$v=v_{\varepsilon,\varphi}$, proving part~(ii).

Since $v_\varepsilon(t,\theta)=\widetilde v_\varepsilon(t/T(\varepsilon),\theta)$,
$T(0)=T_c$, and $(\widetilde v_0)_s=0$, the chain rule and
\eqref{eq:cylinder-branch-statement} give
\[
\begin{aligned}
\left.\partial_\varepsilon v_\varepsilon(t,\theta)\right|_{\varepsilon=0}
&=\phi_c(\theta)\cos\!\left(\frac{2\pi t}{T_c}\right)
-\frac{tT'(0)}{T_c^2}(\widetilde v_0)_s\!\left(\frac{t}{T_c},\theta\right)=\phi_c(\theta)\cos(\kappa_c t)=:U(t,\theta).
\end{aligned}
\]
Since $\mathrm d\mathcal H_{v_c}=0$, the identity
\eqref{eq:cylinder-quadratic} with $A=1$ and $B=0$ gives
\[
\begin{aligned}
\left.\frac{\mathrm d}{\mathrm d\varepsilon}
\mathcal H(v_\varepsilon)\right|_{\varepsilon=0}
&=\mathrm d\mathcal H_{v_c}[U]=0,\\
\left.\frac{\mathrm d^2}{\mathrm d\varepsilon^2}
\mathcal H(v_\varepsilon)\right|_{\varepsilon=0}
&=\mathrm d^2\mathcal H_{v_c}[U,U]
+\mathrm d\mathcal H_{v_c}
\left[\left.\partial_\varepsilon^2v_\varepsilon\right|_{\varepsilon=0}\right]=\kappa_c^2\|\phi_c\|_2^2.
\end{aligned}
\]
Analyticity gives the Taylor expansions
\[
\begin{aligned}
\mathcal H(v_\varepsilon)
&=H_c+\frac{\kappa_c^2\|\phi_c\|_2^2}{2}\varepsilon^2+O(\varepsilon^3),\\
\frac{\mathrm d}{\mathrm d\varepsilon}\mathcal H(v_\varepsilon)
&=\kappa_c^2\|\phi_c\|_2^2\varepsilon+O(\varepsilon^2)>0
\qquad\mbox{for }0<\varepsilon\ll1.
\end{aligned}
\]
Shrink $\varepsilon_0$ so that the derivative is positive on
$(0,\varepsilon_0)$, choose $0<\varepsilon_1<\varepsilon_0$, and set
$\delta=\mathcal H(v_{\varepsilon_1})-H_c>0$. Then for every $H\in(H_c,H_c+\delta)$, there exists a unique $\varepsilon(H)\in(0,\varepsilon_1)$ such that $\mathcal H(v_{\varepsilon(H)})=H$. By part~(ii), every solution in $\mathcal S_H\cap\mathcal U$
is a translate of some $v_\eta$ with $\eta>0$ small. Since $\mathcal H(v_\eta)=H=\mathcal H(v_{\varepsilon(H)})$, the strict monotonicity of
$\varepsilon\mapsto\mathcal H(v_\varepsilon)$ implies
$\eta=\varepsilon(H)$.
Thus there is exactly one periodic translation orbit
in $\mathcal S_H\cap\mathcal U$.
\end{proof}

\subsection{The global Fowler family}
We now use compactness to classify all energy levels.

\begin{proposition}\label{prop:global-energy-classification}
Recall  $v_{\mathrm{bub}}^s$ in \eqref{eq:bubblecyl}. 
We have
\begin{equation}\label{eq:energy-range}
 \mathcal H(\mathcal S)=[H_c,0],\qquad
 \mathcal S_{H_c}=\{v_c\},\qquad
 \mathcal S_0=\{0\}\cup\{v_{\mathrm{bub}}^s:s\in\mathbb R\}.
\end{equation}
For every $H\in(H_c,0)$ there is a nonstationary periodic profile $v_H$, unique
up to time translation, such that
\begin{equation}\label{eq:single-periodic-level}
 \mathcal S_H
 =\{v_H(\,\cdot+s,\cdot):s\in\mathbb R\}.
\end{equation}
\end{proposition}

\begin{proof}
The stationary trace satisfies $b_0w=\frac{n-2}{2}w^{\frac{n}{n-2}}$, so $w\in\{0,w_c\}$.
Uniqueness of the regular extension gives  that the only stationary solutions are $0$ and
$v_c$, with energies $0$ and $H_c$. Since $v_{\mathrm{bub}}^s\to0$ as $t\to\pm\infty$, conservation and continuity give $\mathcal H(v_{\mathrm{bub}}^s)=0$.

By compactness of $\mathcal S$ (Lemma \ref{lem:solution-compactness}) and continuity of $\mathcal H$, $\mathcal H$ attains its minimum $H_-$
and maximum $H_+$. Fix $H\in\{H_-,H_+\}$ and suppose that
$H\notin\{H_c,0\}$. Choose $v\in\mathcal S_H$. By compactness, there
are $t_j\to+\infty$ such that $\tau_{t_j}v\to u\in\mathcal S_H$.
Since $u$ is nonstationary, Proposition~\ref{prop:energy-rigidity}
gives $\varepsilon>0$ and, for all large $j$,
\[
\tau_{t_j}v=\tau_{s_j}u,\qquad |s_j|<\varepsilon.
\]
For two such indices $j,k$,
\[
\tau_{t_j-s_j}v=u=\tau_{t_k-s_k}v
\quad\Longrightarrow\quad
\tau_{t_j-t_k-s_j+s_k}v=v.
\]
Fixing $k$ and taking $j$ sufficiently large gives
$t_j-t_k-s_j+s_k>0$, so $v$ is periodic.

By Proposition~\ref{prop:periodic-neighborhood}, nearby solutions
depend real analytically on their two Cauchy coordinates:
\[
\mathbf c=(c_1,c_2)\longmapsto v_{\mathbf c},\qquad
\mathbf c(v_{\mathbf c})=\mathbf c,\qquad
v_{\mathbf c(v)}=v,
\]
with $\mathbf c$ ranging over an open neighborhood of $\mathbf c(v)$
in $\mathbb R^2$. Set $b(t)=\frac n2v(t,\pi/2)^{\frac{2}{n-2}}$ and
\[
h_i:=\left.\partial_{c_i}v_{\mathbf c}\right|_{\mathbf c=\mathbf c(v)},
\qquad i=1,2.
\]
As in the proof of Proposition~\ref{prop:periodic-neighborhood},
$h_1,h_2$ grow at most linearly in $t$, hence
belong to $\mathcal C_b$. Differentiating the Cauchy-coordinate identity
and using $\dim\mathcal C_b=2$ gives
\[
\mathbf c(h_1)=(1,0),\qquad \mathbf c(h_2)=(0,1),\qquad
\mathcal C_b=\operatorname{span}\{h_1,h_2\}.
\]
Since $H$ is an extremal value, \eqref{eq:H-flux} gives
\[
0=\left.\partial_{c_i}\mathcal H(v_{\mathbf c})
\right|_{\mathbf c=\mathbf c(v)}
=\mathrm d\mathcal H_v[h_i]=\Omega(v_t,h_i),
\qquad i=1,2.
\]
Thus $\Omega(v_t,h)=0$ for every $h\in\mathcal C_b$.
Since $v_t\in\mathcal C_b$, Proposition~\ref{prop:symplectic} implies
$v_t=0$, a contradiction. Therefore
\[
H_-,H_+\in\{H_c,0\},
\qquad
H_-\le H_c<0\le H_+.
\]
Consequently,
\[
 H_-=H_c,
 \qquad
 H_+=0,
 \qquad
 \mathcal H(\mathcal S)\subset[H_c,0].
\]

For $H\in(H_c,0)$, $\mathcal S_H$ is compact and contains no stationary
solution. The preceding compactness and translation argument shows that
every $v\in\mathcal S_H$ is periodic. By
Proposition~\ref{prop:energy-rigidity}, each $v$ has an open neighborhood
in $\mathcal S_H$ containing only its translates. Taking a finite
subcover, we conclude that $\mathcal S_H$ is a finite union of periodic
translation orbits.

For a nonstationary periodic solution $v$, $\Omega(v_t,h_i)$ cannot vanish for both $i=1,2$. Suppose, without loss of generality, that \(\Omega(v_t,h_1)>0\). The other possibilities are handled by choosing the second coordinate or reversing the direction. We can therefore define
\[
v_\sigma:=v_{(c_1(v)+\sigma,\;c_2(v))}
\]
for both positive and negative sufficiently small \(\sigma\). Then
\[
v_0=v,\qquad
\left.\frac{\mathrm d}{\mathrm d\sigma}v_\sigma\right|_{\sigma=0}=h_1.
\]
Taylor’s formula gives
\[
\mathcal H(v_\sigma)
=\mathcal H(v)+\sigma\,\Omega(v_t,h_1)+o(\sigma).
\]
Therefore,
\[
\mathcal H(v_{-\sigma})<\mathcal H(v)<\mathcal H(v_\sigma)
\qquad\mbox{for }0<\sigma\ll1.
\]
Thus neither $\mathcal S_{H_c}$ nor $\mathcal S_0$ contains a
nonstationary periodic solution.

Fix a nonstationary nonperiodic $v\in\mathcal S_H$ with $H\in\{H_c,0\}$.
For any sequence $|t_j|\to\infty$, compactness gives, after passing to
a subsequence, $\tau_{t_j}v\to u\in\mathcal S_H$.
If $u$ were nonstationary, we proved at the beginning of this proof that $v$ must be periodic. 
Hence every such limit is stationary. Since the only stationary
solutions are $v_c$ and $0$, with energies $H_c$ and $0$, respectively,
compactness gives
\[
\begin{cases}
\tau_s v\longrightarrow v_c,& H=H_c,\\
\tau_s v\longrightarrow0,& H=0,
\end{cases}
\qquad\text{in }\mathcal S\text{ as }s\to\pm\infty.
\]
Indeed, otherwise a sequence of translates remaining outside a
neighborhood of the indicated limit would have a subsequence converging
to a different limit.

By Proposition~\ref{prop:cylinder-neighborhood}, every solution
sufficiently close to $v_c$ is either $v_c$ or a translate of a small
branch solution, and
\[
\mathcal H(v_\varepsilon)-H_c
=\frac{\kappa_c^2\|\phi_c\|_2^2}{2}\varepsilon^2+O(\varepsilon^3)>0
\qquad\mbox{for }0<|\varepsilon|\ll1.
\]
Thus $\mathcal S_{H_c}\cap\mathcal U_c=\{v_c\}$ for some neighborhood
$\mathcal U_c$ of $v_c$. If $v\in\mathcal S_{H_c}$ were nonstationary,
then, for sufficiently large $s$,
\[
\tau_s v\in\mathcal S_{H_c}\cap\mathcal U_c
\quad\Longrightarrow\quad
\tau_s v=v_c
\quad\Longrightarrow\quad
v=\tau_{-s}v_c=v_c,
\]
a contradiction. Hence, we obtain 
\[
\mathcal S_{H_c}=\{v_c\}.
\]

For $v\in\mathcal S_0\setminus\{0\}$, since $\tau_s v\longrightarrow0$ as $s\to\pm\infty$, Lemma~\ref{lem:solution-compactness} gives
\[
v(s,\pi/2)=(\tau_s v)(0,\pi/2)\rightarrow0
\quad\mbox{as }s\to\pm\infty.
\]
Theorem~\ref{thm:conformal-end-alternative} yields
$v=v_{\mathrm{bub}}^s$ for some $s\in\mathbb R$.
Since $0$ and all $v_{\mathrm{bub}}^s$ have energy $0$, we obtain
\[
\mathcal S_0=\{0\}\cup\{v_{\mathrm{bub}}^s:s\in\mathbb R\}.
\]

Let $\mathcal S_{\mathrm{per}}$ be the nonstationary periodic solutions and
$I=(H_c,0)$. The image of the map $\mathcal H:\mathcal S_{\mathrm{per}}\to I$ is open by $\mathrm d\mathcal H\ne0$, closed by compactness and continuity, and
nonempty by Proposition~\ref{prop:cylinder-neighborhood}. Thus
\[
\mathcal H(\mathcal S_{\mathrm{per}})=I,
\qquad \mathcal H(\mathcal S)=[H_c,0].
\]

For $H\in I$, set
\[
N(H):=\#\{\text{periodic translation orbits in }\mathcal S_H\}.
\]
The preceding argument gives $1\le N(H)<\infty$. Fix $H_0\in I$ and choose
$v_1,\ldots,v_{N(H_0)}$ representing the distinct orbits in
$\mathcal S_{H_0}$. By Proposition~\ref{prop:periodic-neighborhood} and
$\mathrm d\mathcal H_{v_i}\ne0$, the implicit function theorem in the
two Cauchy coordinates gives a real-analytic curve $H\mapsto v_i(H)$,
defined for $|H-H_0|<\delta_i$, with
\[
v_i(H_0)=v_i,\qquad \mathcal H(v_i(H))=H.
\]
Since
\[
\mathrm d\mathcal H_{v_i}[\partial_Hv_i(H_0)]=1,\qquad
\mathrm d\mathcal H_{v_i}[(v_i)_t]=0,\qquad (v_i)_t\ne0,
\]
the two derivatives of $(H,s)\mapsto\tau_s v_i(H)$ at $(H_0,0)$ are
linearly independent. The inverse function theorem therefore gives a
neighborhood $\mathcal U_i$ of $v_i$ such that, after decreasing $\delta_i$,
\[
\mathcal S_H\cap\mathcal U_i
\subset\{\tau_s v_i(H):s\in\mathbb R\}
\quad\mbox{for }|H-H_0|<\delta_i.
\]
Since $N(H_0)<\infty$, we may take
$\delta=\min_{1\le i\le N(H_0)}\delta_i>0$.

After decreasing $\delta$, these continued orbits are pairwise distinct.
Otherwise, there are $H_j\to H_0$ and $i_j\ne k_j$ such that $v_{i_j}(H_j)$
and $v_{k_j}(H_j)$ are translates. Since $1\le i_j,k_j\le N(H_0)$, we can pass to a
subsequence with $i_j=i$ and $k_j=k$. By periodicity, we may write
\[
v_i(H_j)=\tau_{s_j}v_k(H_j),\qquad
0\le s_j<T(v_k(H_j)).
\]
Since $v_k(H_j)\rightarrow v_k(H_0)=v_k$ in $\mathcal S$, Proposition~\ref{prop:periodic-neighborhood} gives
$T(v_k(H_j))\to T(v_k)$. Thus, after passing to a subsequence,
$s_j\to s$, and continuity of translation yields $v_i=\tau_s v_k$,
a contradiction. Hence $N(H)\ge N(H_0)$ for $|H-H_0|<\delta$.

To prove the reverse inequality, suppose that there are $H_j\to H_0$ and
\[
u_j\in\mathcal S_{H_j}\setminus
\bigcup_{i=1}^{N(H_0)}\{\tau_s v_i(H_j):s\in\mathbb R\}.
\]
By compactness and continuity of $\mathcal H$, a subsequence satisfies
$u_j\to u\in\mathcal S_{H_0}$. Write $u=\tau_s v_i$. Then
$\tau_{-s}u_j\to v_i$, so $\tau_{-s}u_j\in\mathcal U_i$ for all large $j$.
The local representation above gives
\[
\tau_{-s}u_j=\tau_{\sigma_j}v_i(H_j)
\quad\Longrightarrow\quad
u_j=\tau_{s+\sigma_j}v_i(H_j),
\]
a contradiction. Therefore $N(H)\le N(H_0)$ for $H$ sufficiently close
to $H_0$, and $N$ is locally constant on $I$.

Fix an open neighborhood $\mathcal U$ of $v_c$ on which the local
classification in Proposition~\ref{prop:cylinder-neighborhood} applies.
We claim that $\mathcal S_H\subset\mathcal U$ for all sufficiently small
$H-H_c>0$. Otherwise, choose $H_j\downarrow H_c$ and
$u_j\in\mathcal S_{H_j}\setminus\mathcal U$. By compactness, a subsequence
satisfies $u_j\to u\in\mathcal S$. Continuity of $\mathcal H$ gives
\[
\mathcal H(u)=\lim_{j\to\infty}\mathcal H(u_j)
=\lim_{j\to\infty}H_j=H_c,
\quad \mbox{and hence},\quad u\in\mathcal S_{H_c}=\{v_c\}.
\]
Thus $u_j\to v_c\in\mathcal U$, so $u_j\in\mathcal U$ for all large $j$,
a contradiction.

By Proposition~\ref{prop:cylinder-neighborhood}, the small periodic
branch gives exactly one translation orbit at each energy $H>H_c$
sufficiently close to $H_c$, and every solution in $\mathcal U$ with
that energy belongs to this orbit. Since the entire level
$\mathcal S_H$ lies in $\mathcal U$, we obtain $N(H)=1$ near $H_c$.
Since $I$ is an interval and $N$ is locally constant, $N(H)=1$ for every $H\in I$.
\end{proof}

To prove strict monotonicity of the periodic profiles, we use the kernel of $\mathcal P^{-1}$.
\begin{lemma}\label{lem:positivekernel}
The inverse $\mathcal P^{-1}$ is convolution with an even function
$g\in L^1(\mathbb R)$ satisfying
\[
 g(t)>0\quad\mbox{for }t\ne0,\quad\mbox{and}\quad
 g'(t)<0\quad\mbox{for }t>0.
\]
For $T>0$, $g_T(t):=\sum_{\ell\in\mathbb Z}g(t+\ell T)$ satisfies
\[
 g_T'(t)<0\qquad\mbox{for }0<t<T/2.
\]
\end{lemma}
\begin{proof}
At $b=0$, set $\kappa_j=\sqrt{\lambda_j(0)}$.
Since $Q_0[f,f]\ge\frac{(n-2)^2}{4}\|f\|_2^2$, we have
$\kappa_j\ge(n-2)/2>0$. Also, $(e_j^\partial(0))^2>0$ by uniqueness for the
eigenvalue ODE, since $\partial_\theta e_j(0,\pi/2)=0$.
The spectral expansion \eqref{eq:specpara} gives
\[
 \frac1{\Theta(\xi)}=R_0(-\xi^2)
 =\sum_{j\ge0}\frac{(e_j^\partial(0))^2}{\kappa_j^2+\xi^2},
 \qquad \xi\in\mathbb R.
\]
Define
\[
 g(t):=\sum_{j\ge0}\frac{(e_j^\partial(0))^2}{2\kappa_j}e^{-\kappa_j|t|},
 \qquad t\ne0.
\]
Tonelli's theorem gives
\[
 \int_{\mathbb R}g(t)\,\mathrm dt
 =\sum_{j\ge0}\frac{(e_j^\partial(0))^2}{\kappa_j^2}
 =R_0(0)=\frac1{\Theta(0)}<\infty.
\]
This summability implies that the series defining $g$ and all
its derivatives converge locally uniformly away from $0$. Consequently,
$g(-t)=g(t)>0$ for $t\ne0$, and
\[
 g'(t)=-\frac12\sum_{j\ge0}(e_j^\partial(0))^2e^{-\kappa_jt}<0
 \quad \mbox{for }t>0.
\]
Since the Fourier transform of $e^{-\kappa|t|}/(2\kappa)$ is
$(\kappa^2+\xi^2)^{-1}$, convergence in $L^1(\mathbb R)$ permits
termwise Fourier transformation and gives $\widehat g(\xi)=\Theta(\xi)^{-1}$.
For $f\in L^2(\mathbb R)$, Young's inequality gives
$\|g*f\|_{L^2}\le\|g\|_{L^1}\|f\|_{L^2}$, where
$(g*f)(t)=\int_{\mathbb R}g(t-s)f(s)\,\mathrm ds$.
The convolution theorem gives
\[
 \widehat{g*f}(\xi)=\widehat g(\xi)\widehat f(\xi)
 =\frac{\widehat f(\xi)}{\Theta(\xi)}.
\]
Since $\Theta>0$ and $\Theta(\xi)\asymp1+|\xi|$, we have
$g*f\in H^1(\mathbb R)$ and
\[
 \widehat{\mathcal P(g*f)}(\xi)
 =\Theta(\xi)\widehat{g*f}(\xi)=\widehat f(\xi),
 \qquad \mathcal P(g*f)=f.
\]
Moreover, $\mathcal P w=0$ implies $\Theta\widehat w=0$, hence $w=0$.
Thus $\mathcal P^{-1}f=g*f$.

For $T>0$, Tonelli's theorem gives
$\int_0^Tg_T(t)\,\mathrm dt=\int_{\mathbb R}g(t)\,\mathrm dt<\infty$.
For $\mu>0$ and $0\le t\le T$, summing the two geometric series yields
\[
 \sum_{\ell\in\mathbb Z}e^{-\mu|t+\ell T|}
 =\frac{e^{-\mu t}+e^{-\mu(T-t)}}{1-e^{-\mu T}}
 =\frac{\cosh(\mu(T/2-t))}{\sinh(\mu T/2)},
\]
and therefore
\[
 g_T'(t)=-\frac12\sum_{j\ge0}(e_j^\partial(0))^2
 \frac{\sinh\left(\kappa_j(T/2-t)\right)}
 {\sinh\left(\kappa_jT/2\right)}
 <0\quad \mbox{for }0<t<T/2.
\]
\end{proof}

\begin{proposition}\label{prop:fowler-family-properties}
For $H\in(H_c,0)$, let $v_H$ be a representative in
\eqref{eq:single-periodic-level} whose boundary trace $w_H(t)=v_H(t,\pi/2)$ attains its maximum at \(t=0\), and let
$T(H)$ be its least period.  Then
$w_H''(0)<0<w_H''(T(H)/2)$, and
\[
\begin{aligned}
 v_H(-t,\theta)&=v_H(t,\theta)=v_H(T(H)-t,\theta),\\
 (v_H)_t(t,\theta)&<0
 \quad\mbox{for all }0<t<T(H)/2,\ 0\le\theta\le\pi/2.
\end{aligned}
\]
For every integer $j\ge0$ the map
\[
 H\longmapsto\bigl(T(H),v_H(T(H)\,\cdot,\cdot)\bigr)
\]
is locally real analytic into
$\mathbb R\times H^j((\mathbb R/\mathbb Z)\times\mathbb S^{n-1}_+)$.  Moreover,
for every $A>0$ and integer $j\ge0$,
\[
 T(H)\rightarrow\frac{2\pi}{\kappa_c},\qquad
 \sup_{\tau\in\mathbb R}
 \|v_H(\,\cdot+\tau,\cdot)-v_c\|_{C^j([-A,A]\times\overline{\mathbb S^{n-1}_+})}\rightarrow0
 \quad\mbox{as }H\downarrow H_c,
\]
while
\[
 T(H)\rightarrow\infty,
 \qquad
 \|v_H-(\cosh t+\cos\theta)^{-\frac{n-2}{2}}\|_
 {C^j([-A,A]\times\overline{\mathbb S^{n-1}_+})}\rightarrow0
 \quad\mbox{as }H\uparrow0.
\]
\end{proposition}

\begin{proof}
Set $(Rv)(t,\theta):=v(-t,\theta)$. Then $R$ preserves the equation and
$\mathcal H$, and $R\tau_s=\tau_{-s}R$.
By Proposition~\ref{prop:global-energy-classification}, $R$ maps the unique
periodic orbit at a fixed energy onto itself.  If $Rv=\tau_\sigma v$ on an
orbit of least period $T$, then $R\tau_s v=\tau_{-s} Rv=\tau_{\sigma-s}v$. A number \(s\) is a reflection center of \(v\) precisely when
\(v(s+t,\theta)=v(s-t,\theta)\), or equivalently \(R\tau_s v=\tau_s v\). Hence $\tau_s v=\tau_{\sigma-s}v$, that is, \(s\in\frac{\sigma}{2}+\frac T2\mathbb Z.\) Therefore, we choose a
representative even about $0$ and $T/2$.

Let $I_{\mathrm{mon}}\subset(H_c,0)$ consist of the energies admitting an even representative $v_H$ such that
\[
-w_H'(t)>0\quad\text{on }(0,T(H)/2),\qquad
w_H''(0)\ne0,\qquad w_H''(T(H)/2)\ne0.
\] The branch expansion
\eqref{eq:cylinder-branch-statement} shows that $I_{\mathrm{mon}}$ is nonempty.
Fix $H_0\in I_{\mathrm{mon}}$. By
Proposition~\ref{prop:periodic-neighborhood} and $\mathrm d\mathcal H\ne0$,
the implicit function theorem gives a real-analytic family
$H\mapsto\widetilde v_H$ near $H_0$ with
$\mathcal H(\widetilde v_H)=H$ and $\widetilde v_{H_0}=v_{H_0}$.
Write $\widetilde w_H(t)=\widetilde v_H(t,\pi/2)$. Since
\[
 \widetilde w_{H_0}'(0)=0,\qquad
 \left.\partial_s\widetilde w_H'(s)\right|_{(H,s)=(H_0,0)}
 =w_{H_0}''(0)<0,
\]
the scalar implicit function theorem gives a real-analytic $s=s(H)$
such that $s(H_0)=0$ and $\widetilde w_H'(s(H))=0$.
For $v_H:=\tau_{s(H)}\widetilde v_H$, continuity gives
\[
 w_H'(0)=0,\qquad
 w_H''(0)=\widetilde w_H''(s(H))<0
 \quad\mbox{for }|H-H_0|\ll1.
\]
Similarly, $w_{H_0}''(T(H_0)/2)>0$ allows us to continue the
nondegenerate minimum. The profiles converge in $C^2$ after rescaling
to period one, hence the strict sign of the first derivative away from the
two extrema and the signs of the second derivatives near them persist.
So these are the only critical points in each period. They are therefore
the two reflection centers found above, now at $0$ and $T(H)/2$.
Thus $w_H'<0$ on $(0,T(H)/2)$ and $I_{\mathrm{mon}}$ is open.

On $\mathbb R/T(H)\mathbb Z$,
\[
 w_H=\frac{n-2}{2}\,g_{T(H)}*w_H^{\frac{n}{n-2}},
\]
where $g_T$ is the one in Lemma \ref{lem:positivekernel}.
Let $H_j\in I_{\mathrm{mon}}$, $H_j\to H\in(H_c,0)$. By compactness and
periodic continuation, a subsequence satisfies $v_{H_j}\to v_H$ smoothly
and $T(H_j)\to T(H)$. The limit is even about $0$, $w_H' \le0$ on $(0,T(H)/2)$, with $w_H' \not\equiv0$. Differentiating and using odd reflection gives
\begin{equation*}
w_H' (t)=\frac n2\int_0^{T/2}
 [g_T(s-t)-g_T(t+s)]w_H(s)^{\frac{2}{n-2}}w_H' (s)\,\mathrm{d} s.
\end{equation*}
By Lemma~\ref{lem:positivekernel}, the bracket is positive for
$t,s\in(0,T/2)$, so $w_H' <0$ on $(0,T(H)/2)$. Divide by $t$, let $t\downarrow0$, and restrict
the integral to a fixed compact interval in $(0,T(H)/2)$ to obtain $w_H''(0)<0$;
similarly, $w_H''(T/2)>0$. Thus $H\in I_{\mathrm{mon}}$.
Since $I_{\mathrm{mon}}$ is nonempty, open, and closed,
$I_{\mathrm{mon}}=(H_c,0)$. 

To extend $w_H'<0$ to the interior, set $T=T(H)$.
Differentiating \eqref{eq:cylinder} and using reflection symmetry about
$0$ and $T/2$, we obtain
\[
\begin{cases}
\left(-\partial_t^2-\Delta_{\mathbb S^{n-1}}+\dfrac{(n-2)^2}{4}\right)(v_H)_t=0
&\text{in }(0,T/2)\times\mathbb S^{n-1}_+,\\
(v_H)_t(t,\pi/2)=w_H'(t)<0
&0<t<T/2,\\
(v_H)_t(0,\theta)=(v_H)_t(T/2,\theta)=0
&0\le\theta\le\pi/2.
\end{cases}
\]
We apply the maximum principles to the smooth rotational extension of
$(v_H)_t$. Its boundary values are nonpositive, so the weak maximum principle
gives $(v_H)_t\le0$. Since its equatorial trace is negative,
$(v_H)_t\not\equiv0$, and the strong maximum principle gives
$(v_H)_t<0$ in the interior.
The pole $\theta=0$ is an interior point of the hemisphere, so this
also includes $(v_H)_t(t,0)<0$. Together with the equatorial inequality,
\[
 (v_H)_t(t,\theta)<0
 \quad\mbox{for all }0<t<T(H)/2,\ 0\le\theta\le\pi/2.
\]

Fix $H_0\in(H_c,0)$ and use the Cauchy coordinates of
Proposition~\ref{prop:coordinates}. By
Proposition~\ref{prop:periodic-neighborhood}, their inverse
$(c_1,c_2)\mapsto v_{c_1,c_2}$ is real analytic near
$\mathbf c(v_{H_0})$. Write $v_{c_1,c_2}$ for
the solution with coordinates $(c_1,c_2)$ and
$\mathbf c(v_{H_0})=(c_1^0,0)$. Since
$\mathbf c(Rv)=(c_1,-c_2)$ whenever $\mathbf c(v)=(c_1,c_2)$,
after shrinking to an $R$-invariant neighborhood, injectivity gives
\[
 Rv_{c_1,c_2}=v_{c_1,-c_2},\qquad
 \mathcal H(v_{c_1,c_2})=\mathcal H(v_{c_1,-c_2}).
\]
Thus $\left.\partial_{c_2}\mathcal H(v_{c_1,c_2})\right|_{c_2=0}=0$.
Since $\mathrm d\mathcal H_{v_{H_0}}\ne0$ by \eqref{eq:dHkernel},
\[
 \left.\frac{\mathrm d}{\mathrm d c_1}\mathcal H(v_{c_1,0})
 \right|_{c_1=c_1^0}\ne0.
\]
The scalar analytic implicit function theorem gives an analytic
$c_1=c_1(H)$ with $c_1(H_0)=c_1^0$ and
$\mathcal H(v_{c_1(H),0})=H$. Set $v_H=v_{c_1(H),0}$.
Since $Rv_H=v_H$, we have $w_H'(0)=0$; by continuity,
$w_H''(0)<0$ persists, so the maximum remains at $0$.
The least period $T(H)$ is analytic by
Proposition~\ref{prop:periodic-neighborhood}.

Since $s\mapsto w_{H_0}(T(H_0)s)$ is smooth, the construction in the
proof of Proposition~\ref{prop:periodic-neighborhood} applies at this
same profile in every $H^\ell(\mathbb R/\mathbb Z)$, $\ell>5/2$.
The same kernel argument and elliptic regularity give an invertible
linearization at each Sobolev level. Thus $(c_1,c_2)\mapsto(T,W)$ is
analytic into $\mathbb R\times H^\ell$, and local uniqueness identifies
these maps at different Sobolev levels. Composing with
$H\mapsto(c_1(H),0)$ gives, with $W_H(s):=w_H(T(H)s)$,
analyticity of $H\mapsto(T(H),W_H)$ into $\mathbb R\times H^\ell$.
The analytic dependence of the Poisson operator established there gives
$v_H(T(H)\,\cdot,\cdot)=v_{T(H),W_H}$ analytically in
$H^{\ell+1/2}((\mathbb R/\mathbb Z)\times\mathbb S^{n-1}_+)$.
Choosing $\ell+1/2\ge j$, we obtain, for every integer $j\ge0$,
local real analyticity of
\[
H\longmapsto\bigl(T(H),v_H(T(H)\,\cdot,\cdot)\bigr)
\]
into $\mathbb R\times H^j((\mathbb R/\mathbb Z)\times\mathbb S^{n-1}_+)$.

At the cylinder endpoint, compactness and
$\mathcal S_{H_c}=\{v_c\}$ imply, for every fixed $R,j$,
\[
 \sup_{\tau\in\mathbb{R}}
 \|v_H(\cdot+\tau)-v_c\|_{C^j([-R,R]\times\mathbb{S}^{n-1}_+)}\longrightarrow0\quad\mbox{as } H\downarrow H_c.
\]
By Proposition~\ref{prop:cylinder-neighborhood} (ii), for $H-H_c>0$
sufficiently small, $v_H$ is a translate of a solution
$v_{\varepsilon(H)}$ on the local cylinder branch. The energy expansion
in part~(iii) gives, after restricting to a smaller branch neighborhood,
\[
\begin{aligned}
H-H_c
&=\frac{\kappa_c^2\|\phi_c\|_2^2}{2}\varepsilon(H)^2
  +O(|\varepsilon(H)|^3)\ge\frac{\kappa_c^2\|\phi_c\|_2^2}{4}\varepsilon(H)^2.
\end{aligned}
\]
Hence $\varepsilon(H)\to0$ as $H\downarrow H_c$.
Since translation preserves the least period,
\eqref{eq:cylinder-branch} yields
\[
T(H)=T(\varepsilon(H))
=\frac{2\pi}{\kappa_c}+O(|\varepsilon(H)|)
\longrightarrow\frac{2\pi}{\kappa_c}.
\]

At the bubble endpoint, for \(\bar v(\theta):=\frac1T\int_0^T v(t,\theta)\,\mathrm dt\), it satisfies the same angular equation as in \eqref{eq:f0}. Hence, \(\bar v(\theta)=\left(\frac1T\int_0^T w\,\mathrm dt\right)f_0(\theta).\) Differentiating at \(\theta=\pi/2\), using \(f_0'(\pi/2)=b_0\) and the boundary condition, yields
\[
 b_0\int_0^T w\,\mathrm{d} t
 =\frac{n-2}{2}\int_0^T w^{\frac{n}{n-2}}\,\mathrm{d} t.
\]
Since $\frac{n-2}{2}w^{\frac{n}{n-2}}-b_0w$ changes sign only at $w_c$, every nonstationary 
periodic profile satisfies
\begin{equation}\label{eq:trace-crossing}
 \min w<w_c<\max w.
\end{equation}
Centering the maximum at $t=0$ and using \eqref{eq:energy-range}, one has
\[
 v_{H_k}\rightarrow
 v_{\mathrm{bub}}^0(t,\theta):=(\cosh t+\cos\theta)^{-\frac{n-2}{2}}
 \quad\text{in }C^j([-R,R]\times\mathbb{S}^{n-1}_+)\quad\mbox{as } H_k\uparrow0
\]
for every fixed $R,j$. Indeed, any subsequential limit belongs to
$\mathcal S_0\setminus\{0\}$ by \eqref{eq:trace-crossing}, so it equals
$v_{\mathrm{bub}}^s$. The maximum at $0$ forces $s=0$.
Thus the whole family converges and $\max w_H\to1$.

If $T(H_k)\to T_\infty\in(0,\infty)$ as $H_k\uparrow0$, then $v_{\mathrm{bub}}^0$ would be
$T_\infty$-periodic. If $T(H_k)\to0$, then \eqref{eq:sharp} gives
$|w_{H_k}'|\le\frac{n-2}{2}$. Hence,
\[
 \max w_{H_k}-\min w_{H_k}
 \le\int_0^{T(H_k)}|w_{H_k}'(t)|\,\mathrm dt
 \le \frac{n-2}{2}T(H_k)\to0.
\]
Both alternatives are impossible, hence $T(H)\to\infty$ as $H\uparrow0$.
For every $R>0$,
\[
 0\le\limsup_{H\uparrow0}\min w_H
 \le\lim_{H\uparrow0}w_H(R)=(\cosh R)^{-\frac{n-2}{2}};
\]
letting $R\to\infty$ yields $\min w_H\to0$.
\end{proof}

\begin{proof}[Proof of Theorem~\ref{thm:linearized-classification}]
Set $b(t)=\frac n2w_H(t)^{\frac2{n-2}}$. By
Lemma~\ref{lem:linearized-axial-symmetry}, $h$ is axially symmetric and
$h\in\mathcal C_b$. Proposition~\ref{prop:linear-two-dimensionality}
gives $\dim\mathcal C_b=2$.

For $H_c<H<0$, differentiating the equation in $t$ and $H$ shows that
$(v_H)_t$ and $\partial_Hv_H$ solve \eqref{eq:full-linearized-fowler}.
Proposition~\ref{prop:fowler-family-properties} and differentiation of
$v_H(t+T(H),\theta)=v_H(t,\theta)$ give
\[
\partial_Hv_H(t+T(H),\theta)
=\partial_Hv_H(t,\theta)-T'(H)(v_H)_t(t,\theta).
\]
Thus $(v_H)_t$ is bounded and $\partial_Hv_H=O(1+|t|)$, so both belong
to $\mathcal C_b$ by \eqref{eq:weighted-form}. Moreover,
\[
\mathrm d\mathcal H_{v_H}[(v_H)_t]=0,\qquad
\mathrm d\mathcal H_{v_H}[\partial_Hv_H]=1,\qquad
(v_H)_t\not\equiv0.
\]
They are therefore linearly independent and span $\mathcal C_b$.
At $H=H_c$, the proof of Proposition~\ref{prop:energy-rigidity} gives
\[
\mathcal C_{b_c}
=\operatorname{span}\{\phi_c(\theta)\cos(\kappa_ct),\phi_c(\theta)\sin(\kappa_ct)\}.
\]
This proves \eqref{eq:linearized-classification}, including uniqueness
of $A,B$ and the converse.
\end{proof}

\begin{proof}[Proof of Theorem~\ref{thm:fowler-family}]
Theorem~\ref{thm:fowler-family} follows from
Propositions~\ref{prop:global-energy-classification} and~\ref{prop:fowler-family-properties}.
\end{proof}


\section{Asymptotics of local boundary singularities}\label{sec:local-singularities}
Let $U$ be as in Theorem~\ref{thm:local-asymptotics}.  Apply the Emden--Fowler transform
\[
\begin{gathered}
V(t,\omega)=\left(\frac{n-2}{2}\right)^{-\frac{n-2}{2}}e^{-\frac{n-2}{2}t}U(e^{-t}\omega),\qquad
t=-\log r,\qquad t_0=-\log r_0.
\end{gathered}
\]
Then $V$ is defined on $(t_0,\infty)\times\mathbb S^{n-1}_+$ and satisfies
\begin{equation}\label{eq:forward-cylinder}
\begin{cases}
V_{tt}+\Delta_{\mathbb S^{n-1}_+}V-\frac{(n-2)^2}{4}V=0
& \text{in }(t_0,\infty)\times\mathbb S^{n-1}_+,\\
\partial_\nu V=\frac{n-2}{2}V^{\frac{n}{n-2}}
& \text{on }(t_0,\infty)\times\partial\mathbb S^{n-1}_+.
\end{cases}
\end{equation}
By Caffarelli--Jin--Sire--Xiong~\cite{CJSX}*{Theorems~1.1--1.2} with
dimension $n-1$ and $\sigma=1/2$, the growth and Harnack estimates give
some $t_1>t_0$ such that
\begin{equation}\label{eq:local-critical-bound}
0<c_0\le V(t,\omega)\le C_0
\quad\mbox{for }t\ge t_1,\ \omega\in\overline{\mathbb S^{n-1}_+}.
\end{equation}
Set
\[
w(t,\zeta)=V(t,\pi/2,\zeta)=V|_{\partial\mathbb S^{n-1}_+},\qquad
\overline w(t)=\frac1{|\mathbb S^{n-2}|}\int_{\mathbb S^{n-2}}w(t,\zeta)\,\mathrm d\zeta.
\]
Their asymptotic symmetry estimate gives
\begin{equation}\label{eq:cjsx-angular}
\|w(t,\cdot)-\overline w(t)\|_{L^\infty(\mathbb S^{n-2})}
\le Ce^{-t}\quad\mbox{for }t\gg1.
\end{equation}
Rescaled local boundary regularity gives uniform $C^{j,\alpha}$ bounds on fixed translated cylinders.

For a positive smooth profile $W$ on a cylinder, set
\[
\begin{aligned}
 \mathcal H(W)(t)
 &=\frac{1}{|\mathbb S^{n-2}|}\Biggl[
 \frac12\int_{\mathbb S^{n-1}_+}\left(W_t^2-|\nabla_\omega W|^2-\frac{(n-2)^2}{4}W^2\right)\,\mathrm d\omega\\
 &\qquad+\frac{(n-2)^2}{4(n-1)}\int_{\partial\mathbb S^{n-1}_+}W^{\frac{2(n-1)}{n-2}}\,\mathrm d\sigma\Biggr].
\end{aligned}
\]
For axial profiles this agrees with \eqref{eq:hamiltonian}.  Integration by parts in \eqref{eq:forward-cylinder} gives
\[
\mathcal H(V)(t)\equiv H_\infty,\mbox{ which is a constant for all }t>t_0.
\]

Let $t_i\to\infty$.  By \eqref{eq:local-critical-bound} and the local uniform regularity estimates, after passing to a subsequence,
\[
V(t_i+\cdot)\longrightarrow v_*
\quad\text{smoothly on compact cylinders},
\qquad c_0\le v_*\le C_0.
\]
By Theorems~\ref{thm:conformal-end-alternative} and~\ref{thm:fowler-family}, $v_*$ is axially symmetric, and stationary or periodic.
Since $\mathcal H(v_*)=H_\infty$, Theorem~\ref{thm:fowler-family} gives
$H_\infty\in[H_c,0)$, and all translation limits lie on
\[
\mathcal O_\infty
=\{v_{H_\infty}(\,\cdot+\tau,\cdot):\tau\in\mathbb R\},
\qquad v_{H_c}:=v_c.
\]
Consequently, for every $A>0$ and $0<\alpha<1$,
\begin{equation}\label{eq:orbit-convergence}
 \operatorname{dist}_{C^{2,\alpha}([-A,A]\times\overline{\mathbb S^{n-1}_+})}
 \bigl(V(t+\cdot),\mathcal O_\infty\bigr)\rightarrow0
 \quad\mbox{as }t\to\infty.
\end{equation}

By Theorem~\ref{thm:linearized-classification}, every axially symmetric solution of
the linearized equation at $v_H$, $H\in(H_c,0)$, with at most linear
growth in $t$, has the form
\begin{equation}\label{eq:periodic-jacobi-basis}
 \phi=A(v_H)_t+B\,\partial_Hv_H,
 \qquad
 \mathrm d\mathcal H_{v_H}[(v_H)_t]=0,
 \qquad
 \mathrm d\mathcal H_{v_H}[\partial_Hv_H]=1.
\end{equation}
At $H=H_c$, every bounded axially symmetric solution of the linearized
equation at $v_c$ is of the form
\begin{equation}\label{eq:cylinder-jacobi-basis}
 \phi_c(\theta)\bigl(A\cos(\kappa_ct)+B\sin(\kappa_ct)\bigr),
\end{equation}
with positive quadratic form \eqref{eq:cylinder-quadratic} unless it is zero.

\begin{proposition}\label{prop:tail-convergence}
There exist $\beta>0$ and $\tau_\infty\in\mathbb R$ such that, for every $A>0$ and $0<\alpha<1$,
\begin{equation}\label{eq:tail-convergence}
 \bigl\|V(t+\cdot,\cdot)-v_{H_\infty}(t+\cdot+\tau_\infty,\cdot)\bigr\|_
 {C^{2,\alpha}([-A,A]\times\overline{\mathbb S^{n-1}_+})}
 \le C_{A,\alpha}e^{-\beta t}
\end{equation}
for all sufficiently large $t$. If $H_\infty\in(H_c,0)$, then $\tau_\infty$ is unique modulo $T(H_\infty)$.
\end{proposition}

\begin{proof}
The convergence argument is analogous to that of
Korevaar--Mazzeo--Pacard--Schoen \cite{KMPS}*{Section~4},
using the classification of solutions to the linearized equation.

Fix $\alpha_0\in(0,1)$ and write
$\|W\|_I:=\|W\|_{C^{2,\alpha_0}(I\times\overline{\mathbb S^{n-1}_+})}$
for compact intervals $I\subset\mathbb R$.
If $H_\infty\in(H_c,0)$, set $v_\infty=v_{H_\infty}$, $T=T(H_\infty)$,
and $J=[-2T,2T]$.
For $W\in C^{2,\alpha_0}(J\times\overline{\mathbb S^{n-1}_+})$, set
\[
\Psi(W,\varphi):=
\int_J\int_{\mathbb S^{n-1}_+}
\bigl(W(\xi,\omega)-v_\infty(\xi+\varphi,\omega)\bigr)
(v_\infty)_t(\xi+\varphi,\omega)\,\mathrm d\omega\,\mathrm d\xi.
\]
Differentiating in $\varphi$ gives
\[
\begin{aligned}
\partial_\varphi\Psi(W,\varphi)
={}&-\int_J\int_{\mathbb S^{n-1}_+}
|(v_\infty)_t(\xi+\varphi,\omega)|^2\,\mathrm d\omega\,\mathrm d\xi\\
&+\int_J\int_{\mathbb S^{n-1}_+}
\bigl(W(\xi,\omega)-v_\infty(\xi+\varphi,\omega)\bigr)
(v_\infty)_{tt}(\xi+\varphi,\omega)\,\mathrm d\omega\,\mathrm d\xi.
\end{aligned}
\]
For every $\varphi_0$, we have
$\Psi(v_\infty(\cdot+\varphi_0,\cdot),\varphi_0)=0$ and, by periodicity,
\[
\partial_\varphi\Psi
\bigl(v_\infty(\cdot+\varphi_0,\cdot),\varphi_0\bigr)
=-4\int_0^T\int_{\mathbb S^{n-1}_+}
|(v_\infty)_t(t,\omega)|^2\,\mathrm d\omega\,\mathrm dt<0,
\]
since $(v_\infty)_t\not\equiv0$.
The implicit function theorem therefore gives a unique $\varphi$ near
$\varphi_0$ with $\Psi(W,\varphi)=0$ whenever
$\|W-v_\infty(\cdot+\varphi_0,\cdot)\|_J$ is sufficiently small.
By compactness of the periodic orbit, these neighborhoods can be chosen
uniformly in $\varphi_0\in\mathbb R/T\mathbb Z$.
Taking $W=V(s+\cdot,\cdot)$ and using \eqref{eq:orbit-convergence},
for all large $s$ we obtain a unique nearby phase
$\varphi(s)\in\mathbb R/T\mathbb Z$ satisfying
\begin{equation}\label{eq:orthogonality}
\int_J\int_{\mathbb S^{n-1}_+}
\bigl(V(s+\xi,\omega)-v_\infty(\xi+\varphi(s),\omega)\bigr)
(v_\infty)_t(\xi+\varphi(s),\omega)\,\mathrm d\omega\,\mathrm d\xi=0.
\end{equation}
Define
\[
\begin{aligned}
d(s)&=\inf_{\varphi\in\mathbb R/T\mathbb Z}
\bigl\|V(s+\cdot,\cdot)-v_\infty(\cdot+\varphi,\cdot)\bigr\|_J,\\
\bar d(\tau)&=\sup_{s\ge\tau}d(s).
\end{aligned}
\]

For all sufficiently large $s$, we claim that
\begin{equation}\label{eq:phase-estimates}
\begin{gathered}
 d(s)\le\bigl\|V(s+\cdot,\cdot)-v_\infty(\cdot+\varphi(s),\cdot)\bigr\|_J\le Cd(s),\\
 \operatorname{dist}_{\mathbb R/T\mathbb Z}\bigl(\varphi(s+\delta),\varphi(s)+\delta\bigr)
 \le C\bigl(d(s)+d(s+\delta)\bigr)
\end{gathered}
\end{equation}
for $|\delta|\le T$, with $C$ independent of $s$ and $\delta$. The proof is as follows. 
Choose $\varphi_0$ attaining the infimum in $d(s)$. Use this \(\varphi_0\) as the reference phase in the implicit function theorem. Since \(\Psi(V(s+\cdot,\cdot),\varphi(s))=0\) and \(|\partial_\varphi\Psi|\) is uniformly bounded below in the relevant neighborhood, the mean value theorem gives
\[
\begin{aligned}
\operatorname{dist}_{\mathbb R/T\mathbb Z}(\varphi(s),\varphi_0)
&\le C\bigl|\Psi(V(s+\cdot,\cdot),\varphi_0)\bigr|\le C\bigl\|V(s+\cdot,\cdot)-v_\infty(\cdot+\varphi_0,\cdot)\bigr\|_J
=Cd(s).
\end{aligned}
\]
By smoothness and periodicity of $v_\infty$,
\[
\begin{aligned}
\bigl\|V(s+\cdot,\cdot)-v_\infty(\cdot+\varphi(s),\cdot)\bigr\|_J
&\le d(s)+
\bigl\|v_\infty(\cdot+\varphi_0,\cdot)-v_\infty(\cdot+\varphi(s),\cdot)\bigr\|_J\\
&\le d(s)+C\operatorname{dist}_{\mathbb R/T\mathbb Z}(\varphi(s),\varphi_0)\\
&\le Cd(s).
\end{aligned}
\]
The lower bound follows from the definition of $d(s)$. For the second estimate, take $I=[-T,T]$ and repeat the phase construction
with $J$ replaced by $I$. Then
\begin{align}\label{eq:distinv}
&\operatorname{dist}_{\mathbb R/T\mathbb Z}(\varphi_1,\varphi_2)\le C\bigl\|v_\infty(\cdot+\varphi_1,\cdot)
-v_\infty(\cdot+\varphi_2,\cdot)\bigr\|_I,
\end{align}
if the left-hand side is small, since $(v_\infty)_t\not\equiv0$. Combining this local estimate with compactness and injectivity of the phase map, we obtain \eqref{eq:distinv} for all \(\varphi_1,\varphi_2\in\mathbb R/T\mathbb Z\). For $|\delta|\le T$, we have $I,I+\delta\subset J$, so the first estimate gives
\[
\begin{aligned}
&\bigl\|v_\infty(\cdot+\varphi(s+\delta),\cdot)
-v_\infty(\cdot+\delta+\varphi(s),\cdot)\bigr\|_I\\
&\quad\le\bigl\|V(s+\delta+\cdot,\cdot)
-v_\infty(\cdot+\varphi(s+\delta),\cdot)\bigr\|_I+\bigl\|V(s+\cdot,\cdot)
-v_\infty(\cdot+\varphi(s),\cdot)\bigr\|_{I+\delta}\\
&\quad\le C\bigl(d(s+\delta)+d(s)\bigr),
\end{aligned}
\]
Applying \eqref{eq:distinv} with $\varphi_1=\varphi(s+\delta)$ and
$\varphi_2=\varphi(s)+\delta$ proves the second inequality in
\eqref{eq:phase-estimates}.

For
$s+\min\{0,\xi\}\ge\tau$,
\[
\begin{aligned}
&\bigl\|V(s+\cdot,\cdot)-v_\infty(\cdot+\varphi(s),\cdot)\bigr\|_
{C^0([\xi-2T,\xi+2T]\times\overline{\mathbb S^{n-1}_+})}\\
&\quad=
\bigl\|V(s+\xi+\cdot,\cdot)
-v_\infty(\cdot+\xi+\varphi(s),\cdot)\bigr\|_
{C^0(J\times\overline{\mathbb S^{n-1}_+})}\\
&\quad \le Cd(s+\xi)+
\bigl\|v_\infty(\cdot+\varphi(s+\xi),\cdot)
-v_\infty(\cdot+\xi+\varphi(s),\cdot)\bigr\|_
{C^0(J\times\overline{\mathbb S^{n-1}_+})}\\
&\quad\le Cd(s+\xi)
+C\operatorname{dist}_{\mathbb R/T\mathbb Z}
\bigl(\varphi(s+\xi),\varphi(s)+\xi\bigr)\\
&\quad\le C(1+|\xi|)\bar d(\tau).
\end{aligned}
\]
where in the last inequality, we used iteration over overlapping $J$-intervals and the second inequality in \eqref{eq:phase-estimates}.

If $H_\infty=H_c$, set $v_\infty=v_c$, $J=[-2,2]$,
$d(s)=\|V(s+\cdot)-v_c\|_J$, and $\bar d(\tau)=\sup_{s\ge\tau}d(s)$.
In both cases, $\bar d(\tau)\to0$ by \eqref{eq:orbit-convergence}.

We claim that there are $L,\tau_*>0$ such that
\begin{equation}\label{eq:tail-half}
 \bar d(\tau+L)\le\frac12\bar d(\tau)
 \quad\mbox{for all }\tau\ge\tau_*.
\end{equation}
Otherwise, choose $L_j\to\infty$, $\tau_j\to\infty$, and
$t_j\ge\tau_j+L_j$ with $d_j:=\bar d(\tau_j)>0$ and
$d(t_j)>d_j/2$. Since $\bar d(\tau)\to0$ as $\tau\to\infty$, we have $d_j\to0$.

In the periodic case, pass to a subsequence with
$\varphi(t_j)\to\varphi_*$ and set
$p_j(\xi,\omega)=v_\infty(\xi+\varphi(t_j),\omega)$; in the stationary
case, set $p_j=v_c$. Define
\[
U_j(\xi,\omega)=\frac{V(t_j+\xi,\omega)-p_j(\xi,\omega)}{d_j}.
\]
The preceding estimates bound $U_j$ locally, with at most linear growth
in the periodic case and a uniform bound in the stationary case.
The difference equation and local boundary estimates give, along a diagonal subsequence,
\[
U_j\longrightarrow\phi
\qquad\text{in }C^{2,\alpha_0}_{\mathrm{loc}},
\]
where $\phi$ satisfies
\[
\begin{cases}
\phi_{\xi\xi}+\Delta_{\mathbb S^{n-1}_+}\phi-\frac{(n-2)^2}{4}\phi=0
&\text{in }\mathbb R\times\mathbb S^{n-1}_+,\\
\partial_\nu\phi=\frac n2v_\infty(\xi+\varphi_*,\omega)^{\frac2{n-2}}\phi
&\text{on }\mathbb R\times\partial\mathbb S^{n-1}_+.
\end{cases}
\]
Here we set $\varphi_*=0$ in the stationary case. Since \[
\|U_j\|_J
=\frac{\|V(t_j+\cdot,\cdot)-p_j\|_J}{d_j}
\ge\frac{d(t_j)}{d_j}>\frac12,
\]
we have $\phi\not\equiv0$.

For any fixed $0<\varepsilon<\min\{1,n/4\}$, the limiting growth bounds give
\[
\sup_{\omega\in\overline{\mathbb S^{n-1}_+}}|\phi(\xi,\omega)|
\le C(1+|\xi|)\le C_\varepsilon e^{\varepsilon|\xi|}.
\]
Theorem~\ref{thm:linearized-classification}, applied after a time
translation, implies that $\phi$ is axially symmetric.
In the periodic case,
\[
\phi(\xi,\theta)
=A(v_\infty)_t(\xi+\varphi_*,\theta)
+B\left.\partial_Hv_H(\xi+\varphi_*,\theta)\right|_{H=H_\infty}.
\]
Conservation of the Hamiltonian, Taylor's formula, and
\eqref{eq:periodic-jacobi-basis} give
\[
\begin{aligned}
0
&=\frac{\mathcal H(V)(t_j)-\mathcal H(p_j)(0)}{d_j}=\mathrm d\mathcal H_{p_j}[U_j](0)+O(d_j)
\rightarrow
\mathrm d\mathcal H_{v_\infty(\cdot+\varphi_*,\cdot)}[\phi](0)=B.
\end{aligned}
\]
Thus $B=0$. Passing \eqref{eq:orthogonality} to the limit yields
\[
\begin{aligned}
0
&=\int_J\int_{\mathbb S^{n-1}_+}
\phi(\xi,\omega)(v_\infty)_t(\xi+\varphi_*,\omega)
\,\mathrm d\omega\,\mathrm d\xi=A\int_J\int_{\mathbb S^{n-1}_+}
|(v_\infty)_t(\xi+\varphi_*,\omega)|^2
\,\mathrm d\omega\,\mathrm d\xi.
\end{aligned}
\]
Since $(v_\infty)_t\not\equiv0$ and $J$ contains a full period,
$A=0$, contradicting $\phi\not\equiv0$.

In the stationary case, the same theorem gives
\[
\phi(\xi,\theta)=\phi_c(\theta)
(A\cos(\kappa_c\xi)+B\sin(\kappa_c\xi)).
\]
Since $\mathrm d\mathcal H_{v_c}=0$, conservation, Taylor's formula,
and \eqref{eq:cylinder-quadratic} yield
\[
\begin{aligned}
0
&=\frac{\mathcal H(V)(t_j)-H_c}{d_j^2}
=\frac{\mathrm d^2\mathcal H_{v_c}[U_j,U_j](0)}{2}+O(d_j)\rightarrow
\frac{\mathrm d^2\mathcal H_{v_c}[\phi,\phi]}{2}
=\frac{\kappa_c^2\|\phi_c\|_2^2}{2}(A^2+B^2).
\end{aligned}
\]
Hence $A=B=0$, again contradicting $\phi\not\equiv0$.

This proves \eqref{eq:tail-half}.

Iterating \eqref{eq:tail-half} gives
\[
\bar d(\tau_*+kL)\le2^{-k}\bar d(\tau_*)
\quad\mbox{for all }k\ge0.
\]
Since $\bar d$ is nonincreasing, taking $\beta=(\log2)/L$ gives
\begin{equation}\label{eq:tail-rate}
\bar d(t)\le Ce^{-\beta t}.
\end{equation}

For $H_\infty=H_c$, take $\tau_\infty=0$.
Covering $[-A,A]$ by finitely many overlapping translates of $J$
and applying \eqref{eq:tail-rate} gives
\eqref{eq:tail-convergence} with $\alpha=\alpha_0$.

In the periodic case, set $s_k=\tau_*+kT$. By
\eqref{eq:phase-estimates} and \eqref{eq:tail-rate},
\[
\begin{aligned}
\operatorname{dist}_{\mathbb R/T\mathbb Z}
\bigl(\varphi(s_{k+1})-s_{k+1},\varphi(s_k)-s_k\bigr)
&=\operatorname{dist}_{\mathbb R/T\mathbb Z}
\bigl(\varphi(s_{k+1}),\varphi(s_k)+T\bigr)\\
&\le C\bigl(d(s_{k+1})+d(s_k)\bigr)\\
&\le Ce^{-\beta s_k}.
\end{aligned}
\]
Thus $\{\varphi(s_k)-s_k\}$ is Cauchy in $\mathbb R/T\mathbb Z$.
Choose $\tau_\infty\in\mathbb R$ representing its limit. By applying
\eqref{eq:phase-estimates}, for $s_k\le t<s_{k+1}$ we obtain
\[
\begin{aligned}
\operatorname{dist}_{\mathbb R/T\mathbb Z}
\bigl(\varphi(t),t+\tau_\infty\bigr)&\le\operatorname{dist}_{\mathbb R/T\mathbb Z}
\bigl(\varphi(t),\varphi(s_k)+t-s_k\bigr)+\operatorname{dist}_{\mathbb R/T\mathbb Z}
\bigl(\varphi(s_k)-s_k,\tau_\infty\bigr)\\
&\le C\bigl(d(t)+d(s_k)\bigr)+Ce^{-\beta s_k}\\
&\le Ce^{-\beta t}.
\end{aligned}
\]
Combining this with \eqref{eq:phase-estimates} on overlapping $J$-intervals
proves \eqref{eq:tail-convergence} for $\alpha=\alpha_0$, with
$\tau_\infty$ unique modulo $T$.

The difference from the limiting translate satisfies a linear equation
with uniformly smooth Robin coefficient. Local boundary Schauder estimates
give \eqref{eq:tail-convergence} for every $0<\alpha<1$, with the same $\beta$.
\end{proof}

\begin{proof}[Proof of Theorem~\ref{thm:local-asymptotics}]
Define
\[
U_\infty(X)=\left(\frac{n-2}{2}\right)^{\frac{n-2}{2}}|X|^{-\frac{n-2}{2}}
 v_{H_\infty}(-\log|X|+\tau_\infty,X/|X|).
\]
Taking $r=e^{-t}$ in \eqref{eq:tail-convergence} gives, on the fixed half-annulus $\{X_n\ge0,\ 1/2\le|X|\le2\}$,
\[
\left\|r^{\frac{n-2}{2}} U(r\,\cdot)-r^{\frac{n-2}{2}}U_\infty(r\,\cdot)\right\|_{C^{2,\alpha}}
\le C_\alpha r^\beta.
\]
Since $\inf v_{H_\infty}>0$ by \eqref{eq:local-critical-bound},
$U(X)=U_\infty(X)\bigl(1+O(|X|^\beta)\bigr)$.
Uniqueness of the limiting orbit, and of the phase modulo the least period
in the periodic case, gives uniqueness of $U_\infty$.
\end{proof}

\appendix

\section{Necksize}

We show that the nonstationary Fowler family can be parametrized by
its boundary minimum or maximum. We first establish a uniqueness
lemma for even axially symmetric solutions of \eqref{eq:homoeq}.

\begin{lemma}\label{lem:even-boundary-observation}
Let $b$ satisfy \eqref{eq:coefficient-class}, with $b(-t)=b(t)$ and
$b(0)>0$. Let $h\not\equiv 0$ be an axially symmetric weak solution of
\eqref{eq:homoeq}. Suppose that $h(-t,\theta)=h(t,\theta)$ and for some $C>0$ and
$0\le\sigma<(n-1)/2$,
\[
 \|h(t)\|_2\le Ce^{\sigma|t|}\qquad\mbox{for all }t\in\mathbb R.
\]
Then
\begin{equation}\label{eq:even-boundary-sign}
 h(0,\pi/2)\neq0.
\end{equation}
\end{lemma}

\begin{proof}
Suppose that $h(0,\pi/2)=0$. Gronwall's inequality gives $b(t)\ge b(0)e^{-|t|}>0$.
By \eqref{eq:local-H2} and the time-slice trace theorem,
\[
 \|h(t)\|_{H^1_{\mathrm{ax}}}+\|h_t(t)\|_2
 \le C_n\left(\int_{t-2}^{t+2}\|h(s)\|_2^2\,\mathrm ds\right)^{1/2}
 \le C'e^{\sigma|t|}.
\]
In particular, $(h,h_t)\in C_{\mathrm{loc}}(\mathbb R;
H^1_{\mathrm{ax}}\times L^2)$ and $h_t(0,\theta)=0$.

Write $e=e_0(b)$, $\lambda=\lambda_0(b)$, and
$\delta=(n-2)^2/4-\lambda$.
All spectral quantities are evaluated at $b(t)$ unless otherwise specified.
Lemma~\ref{lem:robin-spectral-estimates} gives
\begin{equation}\label{eq:observation-spectral-range}
 0<\delta\le\frac{n(n-1)}2+\frac14<(n-1)^2,\quad \lambda_j>\frac{n^2}{4}\quad\mbox{for every }j\ge1.
\end{equation}
It follows from \eqref{eq:eigenres},
\eqref{eq:boundary-mass}, and \eqref{eq:log-resolvent-strong} that
\begin{equation}\label{eq:observation-dirichlet-resolvent}
\begin{aligned}
 \sum_{j\ge1}\frac{(e_j^\partial)^2}{\lambda_j-n^2/4}
 &=\frac{(e_0^\partial)^2}{\delta+n-1},\\
 b(e_0^\partial)^2
 &\le\frac{2\delta(\delta+n-1)}
 {2(n-1)+\delta(\delta+n-1)/(\delta+2n)}<2\delta.
\end{aligned}
\end{equation}

Decompose $h=qe+z$ and $h_t=pe+y$, with $z,y\perp e$, so
$q=\langle h,e\rangle$ and $p=\langle h_t,e\rangle$. Set
\[
 \xi_\pm=\sqrt{\frac n2}\,q\pm\sqrt{\frac2n}\,p,
 \qquad A=\frac{2\delta+n-2}{4n},\qquad B=\frac{n-1}{2}-A.
\]
By \eqref{eq:observation-spectral-range},
$0<A\le5(n-1)/16$ and $B\ge3(n-1)/16$, hence $A/B\le5/3$.
Define
\[
\begin{aligned}
 \mathscr E
 &:=Q_b[z,z]+\|y\|_2^2-(n-1)\langle z,y\rangle=\left(Q_b-\frac{(n-1)^2}{4}\right)[z,z]
   +\left\|y-\frac{n-1}{2}z\right\|_2^2\ge0,\\
F&:=\langle z,y\rangle-\frac n4q^2-\frac1n p^2.
\end{aligned}
\]
Lemma~\ref{lem:principal-eigenfunction-variation} and the weak product rule
\eqref{eq:weak-product} give
\[
\begin{aligned}
 q'&=p+\langle e',z\rangle,\qquad
 p'=\lambda q+\langle e',y\rangle,\\
 \frac{\mathrm d}{\mathrm dt}\langle z,y\rangle
 &=Q_b[z,z]+\|y\|_2^2-q\langle e',y\rangle-p\langle e',z\rangle.
\end{aligned}
\]
Consequently, almost everywhere,
\begin{equation}\label{eq:observation-moving-form}
 F'-(n-1)F
 =A\xi_+^2+B\xi_-^2+\mathscr E
 -\xi_+\left\langle e',\sqrt{\frac n2}\,z+\sqrt{\frac2n}\,y\right\rangle
 =:\mathscr D.
\end{equation}

We estimate the last term. 
Fix \(t\) and evaluate the spectral quantities at \(b(t)\). Write \(z_j=\langle z,e_j\rangle\) and \(y_j=\langle y,e_j\rangle\). Since \(z,y\perp e_0\), completing the square gives
\[
\begin{aligned}
\mathscr E
&=\sum_{j\ge1}\bigl(\lambda_jz_j^2+y_j^2-(n-1)z_jy_j\bigr)=\sum_{j\ge1}\left[
\left(\lambda_j-\frac{(n-1)^2}{4}\right)z_j^2
+\left(y_j-\frac{n-1}{2}z_j\right)^2
\right].
\end{aligned}
\]
The coefficients are positive because \(\lambda_j>n^2/4\).
Also,
\[
\begin{aligned}
\left|\left\langle e',\sqrt{\frac n2}\,z+\sqrt{\frac2n}\,y\right\rangle\right|
&=\left|\sum_{j\ge1}\langle e',e_j\rangle
\left[
\frac{2n-1}{\sqrt{2n}}z_j
+\sqrt{\frac2n}\left(y_j-\frac{n-1}{2}z_j\right)
\right]\right|\\
&\le\left(\sum_{j\ge1}|\langle e',e_j\rangle|^2
\left[
\frac{(2n-1)^2}{2n\left(\lambda_j-\frac{(n-1)^2}{4}\right)}
+\frac2n
\right]\mathscr E\right)^{1/2},
\end{aligned}
\]
where we used the Cauchy–Schwarz inequality in the last inequality. Using
$\langle e',e_j\rangle=\langle b'\partial_b e(b),e_j\rangle=b'e_0^\partial e_j^\partial/(\lambda_j-\lambda)$, we obtain
\[
\begin{aligned}
& \left|\left\langle e',\sqrt{\frac n2}\,z+\sqrt{\frac2n}\,y\right\rangle\right|^2\\&\le (b')^2(e_0^\partial)^2
 \sum_{j\ge1}
 \frac{(e_j^\partial)^2(\frac{3n}{2}-1+\frac{2\lambda_j}{n})}
 {(\lambda_j-\lambda)^2(\lambda_j-\frac{(n-1)^2}{4})}\,\mathscr E\\
 &\le (b')^2(e_0^\partial)^2
 \sum_{j\ge1}
 \frac{(2n-1)(e_j^\partial)^2}
 {(\delta+n-1)(\sqrt{\delta+n-1}+\frac12\sqrt{2n-1})^2(\lambda_j-\frac{n^2}{4})}\,\mathscr E\\
 &\le\frac{(2n-1)b^2(e_0^\partial)^4}
 {(\delta+n-1)^2(\sqrt{\delta+n-1}+\frac12\sqrt{2n-1})^2}\,\mathscr E\\
 &\le\frac{2\sqrt{2n-1}\,\delta^2}{(\delta+n-1)^{5/2}}\,\mathscr E,
\end{aligned}
\]
where we also used for $x>0$ that
\[
\begin{gathered}
 \frac{2n-1+2x/n}{x+\delta+n-1}\le\frac{2n-1}{\delta+n-1},\\
 (x+\delta+n-1)\left(x+\frac{2n-1}{4}\right)
 \ge x\left(\sqrt{\delta+n-1}+\frac{\sqrt{2n-1}}2\right)^2
\end{gathered}
\]
in the second inequality,  \eqref{eq:observation-dirichlet-resolvent} and $|b'|\le b$ in the third inequality,  \eqref{eq:observation-dirichlet-resolvent} and the Cauchy--Schwarz inequality in the fourth inequality.
With $x=\delta/(n-1)$, we have
\[
\begin{aligned}
 \left(\frac{\sqrt{2n-1}\,\delta^2}{2A(\delta+n-1)^{5/2}}\right)^2
 \le\frac{90x^4}{(4x+1)^2(1+x)^5}
 \le\frac{15}{2}\left(\frac{x}{(1+x)^2}\right)^3\le\frac{15}{128}<\left(\frac35\right)^4,
\end{aligned}
\]
where $(4x+1)^2\ge12x(1+x)$ and $(1+x)^2\ge4x$. Therefore, it follows that 
\begin{equation}\label{eq:observation-scalar-bounds}
 \mathscr D\ge A\xi_+^2+B\xi_-^2+ \mathscr E -\frac65 \sqrt{A}\sqrt{\mathscr E} |\xi_+| \ge\frac25A\xi_+^2+B\xi_-^2+\frac25\mathscr E\ge0.
\end{equation}
Moreover, \eqref{eq:observation-dirichlet-resolvent} and
Cauchy--Schwarz give directly
\[
\begin{aligned}
|z(\pi/2)|^2
=\left|\sum_{j\ge1}e_j^\partial\langle z,e_j\rangle\right|^2
&\le
\left(\sum_{j\ge1}\frac{(e_j^\partial)^2}{\lambda_j-n^2/4}\right)
\left(\sum_{j\ge1}(\lambda_j-n^2/4)|\langle z,e_j\rangle|^2\right)\\
&=\frac{(e_0^\partial)^2}{\delta+n-1}
\left(Q_b[z,z]-\frac{n^2}{4}\|z\|_2^2\right)\\
& \le\frac{(e_0^\partial)^2}{\delta+n-1}\,\mathscr E.
\end{aligned}
\]
Since 
\[
h(t,\pi/2)=q\,e_0^\partial+z(\pi/2)=\frac{\xi_++\xi_-}{\sqrt{2n}}\,e_0^\partial+z(\pi/2),
\]
it follows from the triangle inequality, 
$b(e_0^\partial)^2/(4n)<A\le 5B/3$, and
$b(e_0^\partial)^2/[2(\delta+n-1)]<1$ that
\begin{equation}\label{eq:observation-moving-trace}
\begin{aligned}
 \frac b2|h(t,\pi/2)|^2
 &\le\left(\sqrt{\frac{b(e_0^\partial)^2}{4n}}
\bigl(|\xi_+|+|\xi_-|\bigr)
+\sqrt{\frac{b(e_0^\partial)^2}{2(\delta+n-1)}}\sqrt{\mathscr E}\right)^2\\
 &\le\left(\sqrt A|\xi_+|+\sqrt{\frac{5B}{3}}|\xi_-|
              +\sqrt{\mathscr E}\right)^2\\
 &\le\frac{20}{3}\left(\frac25A\xi_+^2+B\xi_-^2
                         +\frac25\mathscr E\right)\\
& \le\frac{20}{3}\mathscr D.
\end{aligned}
\end{equation}

Since $h_t(0)=0$ and $2\sigma<n-1$,
\[
 F(0)=-\frac n4q(0)^2,\qquad
 e^{-(n-1)t}F(t)=O(e^{-(n-1-2\sigma)t})\rightarrow0
 \quad\mbox{as }t\to\infty.
\]
Integrating \eqref{eq:observation-moving-form} against $e^{-(n-1)t}$
and using \eqref{eq:observation-moving-trace}, we obtain
\begin{equation}\label{eq:observation-flux-budget}
\begin{aligned}
 \int_0^\infty e^{-(n-1)t}\mathscr D\,\mathrm dt
 &=\frac n4q(0)^2,\\
 I:=\int_0^\infty e^{-(n-1)t}b(t)|h(t,\pi/2)|^2\,\mathrm dt
 &\le\frac{10n}{3}q(0)^2.
\end{aligned}
\end{equation}
Therefore, if $q(0)=0$  then  $\mathscr D=0$, and consequently, $h\equiv0$ by
\eqref{eq:observation-scalar-bounds} and evenness.

For the remainder of the proof, freeze the spectral quantities at $b(0)$
and set $\kappa_j=\sqrt{\lambda_j}>n/2$ for $j\ge1$.
Since $h(0,\pi/2)=0$, by \eqref{eq:Dirifirst}, we have
\begin{equation}\label{eq:observation-terminal}
 (\delta+n-1)q(0)^2\le
 \sum_{j\ge1}\left(\lambda_j-\frac{n^2}{4}\right)
              |\langle h(0),e_j\rangle|^2.
\end{equation}
The weak equation gives
\[
 (-\partial_t^2+\kappa_j^2)\langle h(t),e_j\rangle
 =(b(t)-b(0))e_j^\partial h(t,\pi/2).
\]
Multiply by $e^{-\kappa_jt}$ and integrate on $[0,\infty)$.
Since $h_t(0)=0$ and $\sigma<\kappa_j$, integration by parts yields
\[
 \langle h(0),e_j\rangle
 =\frac{e_j^\partial}{\kappa_j}\int_0^\infty
 e^{-\kappa_jt}(b(t)-b(0))h(t,\pi/2)\,\mathrm dt.
\]
By $|b'|\le b$, we have $e^{-t}\le b(t)/b(0)\le e^t$ for $t\ge0$,
and hence $(b(t)-b(0))^2/b(t)\le b(0)(e^t-2+e^{-t})$. It follows from the Cauchy--Schwarz inequality, \eqref{eq:observation-dirichlet-resolvent} and
\eqref{eq:observation-flux-budget} that
\[
\begin{aligned}
 |\langle h(0),e_j\rangle|^2
 &=\left(\frac{e_j^\partial}{\kappa_j}\int_0^\infty\left[e^{-\frac{n-1}{2}t}\sqrt{b(t)}\,h(t,\pi/2)\right]
\left[e^{-(\kappa_j-\frac{n-1}{2})t}
\frac{b(t)-b(0)}{\sqrt{b(t)}}\right]\,\mathrm dt\right)^2\\
 &\le\frac{b(0) (e_j^\partial)^2 I }{\kappa_j^2}
 \int_0^\infty e^{-(2\kappa_j-n+1)t}(e^t-2+e^{-t})\,\mathrm dt\\
 &\le\frac{5  n b(0)  (e_j^\partial)^2 q(0)^2}
 {6\kappa_j^2(\kappa_j-\frac n2)(\kappa_j-\frac{n-1}{2})(\kappa_j-\frac{n-2}{2})}.
\end{aligned}
\]
For $x>n/2$, the arithmetic mean--geometric mean inequality gives
\[
 \frac{(x^2-n^2/4)(x+n/2)}
 {x^2(x-\frac{n-1}{2})(x-\frac{n-2}{2})}
 <\frac4{\frac32+\sqrt2}<\frac32.
\]
Hence, using \eqref{eq:observation-dirichlet-resolvent}, we obtain
\[
\begin{aligned}
 \sum_{j\ge1}\left(\lambda_j-\frac{n^2}{4}\right)
              |\langle h(0),e_j\rangle|^2
 &\le\frac{5n\delta}{2(\delta+n-1)}\,q(0)^2\le\frac{15}{16}(\delta+n-1)q(0)^2.
\end{aligned}
\]
Together with \eqref{eq:observation-terminal}, this gives $q(0)=0$,
and therefore $h\equiv0$.
\end{proof}

We differentiate the Fowler profiles with respect to $H$, keeping
a boundary minimum or maximum at $t=0$.
Lemma~\ref{lem:even-boundary-observation} implies that the corresponding
extremal value has a nonzero derivative. Together with the endpoint
limits, this yields the following parametrization.

\begin{proposition}\label{prop:fowler-extremal-parameters}
Let $v_H$, $w_H$, and $T(H)$ be as in
Proposition~\ref{prop:fowler-family-properties}, and define
\[
 \varepsilon(H):=\min_{t\in\mathbb R}w_H(t)
 =w_H(T(H)/2),\qquad
 M(H):=\max_{t\in\mathbb R}w_H(t)=w_H(0).
\]
Then
\begin{equation}\label{eq:fowler-extremal-diffeomorphisms}
 \varepsilon:(H_c,0)\longrightarrow(0,w_c),\qquad
 M:(H_c,0)\longrightarrow(w_c,1)
\end{equation}
are real-analytic diffeomorphisms, with
\begin{equation}\label{eq:fowler-extremal-derivatives}
 \varepsilon'(H)<0<M'(H)\qquad\mbox{for every }H\in(H_c,0).
\end{equation}
Their endpoint limits are
\[
\begin{aligned}
 &\varepsilon(H)\to w_c,\quad M(H)\to w_c
 &&\text{as }H\downarrow H_c,\\
 &\varepsilon(H)\to0,\quad M(H)\to1
 &&\text{as }H\uparrow0.
\end{aligned}
\]
Consequently, modulo translation in $t$, the nonstationary Fowler
solutions are uniquely indexed by either their boundary minimum
$\varepsilon\in(0,w_c)$ or their boundary maximum $M\in(w_c,1)$.
The representative with its minimum at $t=0$ is uniquely specified by
\[
 w(0)=\varepsilon,\qquad w'(0)=0,\qquad 0<\varepsilon<w_c;
\]
the analogous statement holds at a maximum with $w_c<M<1$.
With either normalization, the least period and the profile rescaled
to period one depend locally real analytically on the corresponding
parameter in every Sobolev space appearing in
Proposition~\ref{prop:fowler-family-properties}.
\end{proposition}

\begin{proof}
Set $W_H(s)=w_H(T(H)s)$ for $s\in\mathbb R/\mathbb Z$.
Proposition~\ref{prop:fowler-family-properties} gives local real
analyticity of $H\mapsto W_H$ in every $H^\ell(\mathbb R/\mathbb Z)$,
$\ell>5/2$. By the boundedness of the trace operator $\varepsilon(H)=W_H(1/2)$ and $M(H)=W_H(0)$ are real analytic. Their endpoint limits follow from the same
proposition and the last paragraph of its proof. In particular, at the
cylinder endpoint the convergence is uniform over all translates.

For $\sigma\in\{0,1/2\}$, set
\[
 u_H^\sigma(t,\theta)=v_H(t+\sigma T(H),\theta),\qquad
 h_H^\sigma=\partial_Hu_H^\sigma,
\]
where the derivative is taken at fixed $t$.
Reflection symmetry makes $u_H^\sigma$ and $h_H^\sigma$ even in $t$.
Differentiating \eqref{eq:cylinder} shows that
\[
 (h_H^\sigma)_{tt}-A_{b_H^\sigma(t)}h_H^\sigma=0,
 \qquad
 b_H^\sigma(t)=\frac n2
 \bigl(u_H^\sigma(t,\pi/2)\bigr)^{\frac2{n-2}}.
\]
The coefficient is positive, even, and satisfies
\eqref{eq:coefficient-class} by Lemma~\ref{lem:secant}.
Differentiation of periodicity gives
\[
 h_H^\sigma(t+T(H),\theta)
 =h_H^\sigma(t,\theta)-T'(H)(u_H^\sigma)_t(t,\theta).
\]
Thus $\|h_H^\sigma(t)\|_2\le C_H(1+|t|)$. Moreover,
\[
 \mathrm d\mathcal H_{u_H^\sigma}[h_H^\sigma]
 =\frac{\mathrm d}{\mathrm dH}\mathcal H(u_H^\sigma)=1,
\]
so $h_H^\sigma\not\equiv0$. Then it follows from
Lemma~\ref{lem:even-boundary-observation} that $M'(H)=h_H^0(0,\pi/2)\ne0$ and $\varepsilon'(H)=h_H^{1/2}(0,\pi/2)\ne0$. Each derivative has a fixed sign on the connected interval $(H_c,0)$.
The endpoint limits determine these signs and give
\eqref{eq:fowler-extremal-derivatives} and the exact ranges in
\eqref{eq:fowler-extremal-diffeomorphisms}.
The analytic inverse function theorem gives real-analytic inverses.

There is exactly one translation orbit at each energy by
Proposition~\ref{prop:global-energy-classification}, and each orbit
has exactly one minimum and one maximum in a least period by
Proposition~\ref{prop:fowler-family-properties}.
Together with \eqref{eq:trace-crossing}, this proves the stated
uniqueness under the two normalizations.
Finally, composition with the analytic inverses proves analytic
dependence on either extremal parameter. Centering the minimum at $t=0$ amounts, after rescaling, to the fixed translation $s\mapsto s+1/2$, which preserves this analyticity.
\end{proof}

\medskip

\textbf{AI usage statement}: The authors proposed using local
bifurcation as the starting point for the global classification.
They used GPT-6 Astra as an interactive research-assistance tool
to explore and develop mathematical arguments, guided by their
questions and directions. In particular, GPT-6 Astra suggested
the local uniqueness argument, which was developed and refined
through subsequent exchanges and revisions. The authors also used
it for language editing and the organization of the exposition.
They take full responsibility for all mathematical statements,
arguments, citations, and conclusions in this paper.

\bigskip

\noindent T. Jin

\noindent Department of Mathematics, The Hong Kong University of Science and Technology\\
Clear Water Bay, Kowloon, Hong Kong\\
Email: \textsf{tianlingjin@ust.hk}

\bigskip

\noindent Y. Liao

\noindent School of Mathematical Sciences, Laboratory of Mathematics and Complex Systems, MOE\\ Beijing Normal University, 
Beijing 100875, China\\
Email: \textsf{yxliao@mail.bnu.edu.cn}

\bigskip

\noindent J. Xiong

\noindent School of Mathematical Sciences, Laboratory of Mathematics and Complex Systems, MOE\\ Beijing Normal University, 
Beijing 100875, China\\
Email: \textsf{jx@bnu.edu.cn}

\end{document}